%% file: main.tex
\documentclass{article}

\usepackage{arxiv}

\usepackage[utf8]{inputenc} 
\usepackage[T1]{fontenc}    
\usepackage[numbers,sort&compress]{natbib} 
\usepackage{hyperref}       
\usepackage{url}            
\usepackage{booktabs}       
\usepackage{longtable}
\usepackage{placeins}  
\usepackage{pgfplots}
\pgfplotsset{compat=1.18}
\usepgfplotslibrary{groupplots}
\usepackage{amsfonts}       
\usepackage{amsmath}
\usepackage{amssymb}
\usepackage{amsthm}
\usepackage{nicefrac}       
\usepackage{microtype}      
\usepackage{graphicx}
\usepackage{multirow}
\usepackage{algorithm}
\usepackage[noend]{algpseudocode}
\graphicspath{ {./} }

\title{Presolve heuristics in HiGHS: implementation and computational study}

\author{
 Simon Spoorendonk \\
  Denmark \\
  \texttt{simon@spoorendonk.dk} \\
  \href{https://orcid.org/0009-0007-4304-6956}{ORCID:~0009-0007-4304-6956} \\
}

\input{claims}

\begin{document}
\maketitle
\begin{abstract}
\input{abstract}
\end{abstract}

\keywords{Mixed-integer programming \and Primal heuristics \and Presolve heuristics \and HiGHS \and Feasibility Jump \and Primal integral}

\input{body}

\bibliographystyle{plainnat}
\bibliography{references}

\input{appendix}

\end{document}

%% file: claims.tex
\newcommand{\ConstNInstances}{233}
\newcommand{\ConstTimeLimitS}{600}
\newcommand{\ConstHighsVersion}{1.15.1}
\newcommand{\ConstShiftHeadline}{0.001}
\newcommand{\ConstEffortfj}{0.3317}
\newcommand{\ConstEffortfpr}{3.161}
\newcommand{\ConstEffortlocalmip}{3.287}
\newcommand{\ConstEffortscylla}{0}
\newcommand{\ConstFjupstreamtotaleffortshift}{10}
\newcommand{\ConstLocalmiprecheckperiod}{100}
\newcommand{\ConstPatiencefj}{0}
\newcommand{\ConstPatiencefpr}{0.3372}
\newcommand{\ConstPatiencelocalmip}{3.194}
\newcommand{\ConstPatiencescylla}{0}
\newcommand{\HeadlineN}{233}
\newcommand{\HeadlineRatio}{0.763}
\newcommand{\HeadlineCI}{[0.684, 0.852]}
\newcommand{\HeadlinePct}{23.7}
\newcommand{\HeadlineFeasiblePatched}{214}
\newcommand{\HeadlineFeasibleVanilla}{211}
\newcommand{\HeadlineFeasibleOnlyPatched}{4}
\newcommand{\HeadlineFeasibleOnlyVanilla}{1}
\newcommand{\HeadlineWinPatched}{49}
\newcommand{\HeadlineWinVanilla}{36}
\newcommand{\HeadlineSgmGapPatched}{0.0054}
\newcommand{\HeadlineSgmGapVanilla}{0.0066}
\newcommand{\HeadlineSgmHeadlinePatched}{0.0484}
\newcommand{\HeadlineSgmHeadlineVanilla}{0.0637}
\newcommand{\HeadlinePairedGapRatio}{0.836}
\newcommand{\HeadlinePairedGapCI}{[0.718, 0.974]}
\newcommand{\HeadlinePairedGapN}{210}
\newcommand{\HeadlineGapP}{0.022}
\newcommand{\HeadlineGapBetter}{42}
\newcommand{\HeadlineGapWorse}{34}
\newcommand{\HeldoutN}{143}
\newcommand{\HeldoutRatio}{0.871}
\newcommand{\HeldoutCI}{[0.782, 0.970]}
\newcommand{\HeldoutP}{0.012}
\newcommand{\HeldoutPct}{12.9}
\newcommand{\HeldoutBetter}{63}
\newcommand{\HeldoutWorse}{61}
\newcommand{\HeldoutPairedGapRatio}{0.871}
\newcommand{\HeldoutPairedGapCI}{[0.720, 1.053]}
\newcommand{\HeldoutPairedGapN}{125}
\newcommand{\TuningN}{90}
\newcommand{\HeldoutUnseenN}{95}
\newcommand{\HeldoutUnseenRatio}{0.861}
\newcommand{\HeldoutUnseenCI}{[0.746, 0.994]}
\newcommand{\HeldoutUnseenP}{0.041}
\newcommand{\SecondArmRatio}{0.965}
\newcommand{\SecondArmCI}{[0.906, 1.027]}
\newcommand{\SecondArmP}{0.260}
\newcommand{\SecondArmBetter}{140}
\newcommand{\SecondArmWorse}{67}
\newcommand{\ReplicationRatio}{0.915}
\newcommand{\ReplicationCI}{[0.841, 0.995]}
\newcommand{\ProbePoolN}{191}
\newcommand{\ProbeFjSoleCracker}{7}
\newcommand{\ProbeFprSoleCracker}{3}
\newcommand{\ProbeLocalMipSoleCracker}{18}
\newcommand{\ProbeScyllaSoleCracker}{2}
\newcommand{\ProbeScyllaProduces}{129}
\newcommand{\ProbeFjVsFprFirstBetter}{81}
\newcommand{\ProbeFjVsFprSecondBetter}{24}
\newcommand{\ProbeFjVsLocalMipFirstBetter}{65}
\newcommand{\ProbeFjVsLocalMipSecondBetter}{46}
\newcommand{\ProbeLocalMipVsFprFirstBetter}{81}
\newcommand{\ProbeLocalMipVsFprSecondBetter}{24}
\newcommand{\PatienceFprRatio}{1844}
\newcommand{\SearchConfirmInstancesN}{49}
\newcommand{\SearchHeldoutInstancesN}{48}
\newcommand{\SearchCandidatesN}{1447}
\newcommand{\SearchCapabilityTimeLimitS}{120}
\newcommand{\SearchCapabilityDiveN}{10}
\newcommand{\SearchMixSurvivorsN}{11}
\newcommand{\SearchEffortTotalMeasured}{29.85}
\newcommand{\SearchHeldoutCheapVsControlBetter}{31}
\newcommand{\SearchHeldoutCheapVsControlWorse}{10}
\newcommand{\SearchConfirmFjOnlyVsMixP}{0.050}
\newcommand{\SearchConfirmFprLpOnSelectedN}{49}
\newcommand{\SearchConfirmFprLpOnMeasuredP}{0.049}
\newcommand{\SearchArmTotalBPrime}{6.779}
\newcommand{\AttrFeasibleN}{214}
\newcommand{\AttrHighsOtherBest}{188}
\newcommand{\AttrChainFjFirstN}{114}
\newcommand{\AttrChainFprImprovesN}{16}
\newcommand{\AttrChainLocalMipImprovesN}{62}
\newcommand{\AttrChainEitherImprovesN}{71}
\newcommand{\AttrChainGapFjMed}{0.624}
\newcommand{\AttrChainGapAfterMed}{0.382}
\newcommand{\AttrOnlyPatchedOursN}{4}
\newcommand{\ProfTotal}{233}
\newcommand{\ProfNever}{22}
\newcommand{\ProfUnderTenS}{154}
\newcommand{\ProfUnderTenSPct}{66.1}

%% file: abstract.tex
We present reference implementations of four published primal heuristics for mixed-integer
programming inside the open-source solver HiGHS~\cite{huangfu2018highs}: Feasibility
Jump~\cite{luteberget2023fj}, fix-propagate-repair (FPR)~\cite{salvagnin2025fpr} with its
LP-guided dive-time variant, LocalMIP~\cite{lin2024localmip}, and the first-order-LP feasibility
pump Scylla~\cite{mexi2023scylla}. They share one integration interface with per-heuristic
budgets and patiences, and each runs in parallel on independently seeded
workers. The evaluation is on the
\texttt{mipfeas} benchmark~\cite{bussieck2026mipfeas}, \ConstNInstances{} instances at
\ConstTimeLimitS\,s, against a separately built unpatched HiGHS on the same machine. The
primal integral improves by \HeadlinePct\% over the full set and by \HeldoutPct\% on the
instances not used for tuning: the patched solver reaches a good solution sooner. At the time
limit it is slightly ahead. It finds a solution on \HeadlineFeasibleOnlyPatched{} instances
where the unpatched solver finds none, against \HeadlineFeasibleOnlyVanilla{} the other way, and
its final gap is smaller, significantly so only over the full set. Scylla and the LP-guided dive produced no accepted incumbent
in any full-limit run that carried them and ship disabled; their implementations remain. The
implementation is released as open source~\cite{mip_heuristics_software}.

%% file: body.tex
\section{Introduction}
\label{sec:intro}

Finding a feasible, and quickly a good, incumbent early is decisive for the practical performance
of a branch-and-bound mixed-integer programming (MIP) solver: a strong incumbent drives bound-based
fixing, pruning, and node selection. Primal heuristics are the components that produce such
incumbents, and a substantial literature has developed and surveyed
them~\cite{berthold2006primal,fischetti2005feaspump,berthold2014rens,%
danna2005rins,luteberget2023fj,salvagnin2025fpr,lin2024localmip,mexi2023scylla};
Section~\ref{sec:related} reviews the strands relevant here.

In this paper we provide
\emph{reference implementations of four published primal heuristics (Feasibility Jump, FPR,
LocalMIP and Scylla) inside one open solver, and an empirical evaluation of them}. The solver is
HiGHS~\cite{huangfu2018highs}. The heuristics are not new. Our contribution is that they now
live in the same solver behind one integration interface, share one budgeting, patience and
submission infrastructure, report their effort usage per dispatch, and are evaluated together on
the \texttt{mipfeas} benchmark~\cite{bussieck2026mipfeas}. The evaluation compares the patched
solver with the unpatched one, which we call \emph{vanilla} HiGHS, and asks whether the added
heuristics improve the full solve. They do, modestly. The primal integral improves: the patched solver reaches good solutions
\emph{sooner}. At the time limit it is slightly ahead. It finds a solution on
\HeadlineFeasibleOnlyPatched{} instances where vanilla HiGHS finds none, against
\HeadlineFeasibleOnlyVanilla{} the other way, and its final gap is smaller, a difference that
separates over the full set but not on the instances held out from tuning. We measure this net effect
only; what the added heuristics displace inside the host is not measured. Our contributions are:

\begin{itemize}
\item \textbf{Four published heuristics, re-implemented on one open stack.} A
  Feasibility Jump~\cite{luteberget2023fj} wrapper over HiGHS's own implementation, aligned with its paper in two places,
  FPR~\cite{salvagnin2025fpr} (with its LP-guided dive-time variant, which we count as a fifth
  mechanism), LocalMIP~\cite{lin2024localmip}, and Scylla~\cite{mexi2023scylla} are integrated
  natively into HiGHS, through a presolve hook and one call in the dive, and released as open
  software at \mbox{\href{https://github.com/spoorendonk/mip-heuristics}{\nolinkurl{github.com/spoorendonk/mip-heuristics}}}~\cite{mip_heuristics_software}.
  Section~\ref{sec:heuristics} documents each at the level of algorithmic detail needed to
  reproduce it, distinguishing what follows the source paper from where our implementation departs
  from it.
\item \textbf{A common infrastructure for integration, budgeting, parallel execution and attribution.} One runner contract,
  a fixed heuristic call chain with a per-heuristic effort budget and an absolute, instance-scaled patience
  each, one report of effort usage, immediate incumbent submission with per-heuristic source tags, and
  a record of accepted solutions from which productive versus stale effort per heuristic is read. All heuristics run in parallel on independently seeded workers, and the infrastructure is documented in Section~\ref{sec:infra}.
\item \textbf{A reproducible evaluation.} A five-stage campaign on the \ConstNInstances-instance
  \texttt{mipfeas} set against a separately built unpatched HiGHS at the same tag on the same
  machine (Section~\ref{sec:results}): (1) a baseline run of that vanilla HiGHS, (2) a
  presolve-only probe that measures each heuristic's productive versus stale effort, (3) a joint
  racing search over heuristic mix, budgets and patiences, (4) a separate measurement of the
  LP-guided dive, which runs after the root LP, and (5) the headline comparison over the full
  set.
\item \textbf{Documented negative results.} Two of the five mechanisms produced no accepted incumbent in any
  full-limit run that carried them, and both ship at effort zero. Scylla finds solutions on
  \ProbeScyllaProduces{} instances when run alone with an unbounded budget, but on none in full
  solves, including a run over all \HeadlineN{} instances (Section~\ref{sec:search}). The
  LP-guided dive \texttt{fpr\_lp} produces nothing in paired full-limit solves on
  \SearchConfirmFprLpOnSelectedN{} instances on each of two backgrounds
  (Section~\ref{sec:fprlp-stage}). The selected
  configuration enables the other three, and no dive-time heuristic of ours. All four
  implementations remain, documented and runnable~\cite{mip_heuristics_software}.
\end{itemize}

\paragraph{Scope of the claim.}
The \texttt{mipfeas} benchmark~\cite{bussieck2026mipfeas} is the feasible part of the
MIPLIB~2017~\cite{miplib2017} benchmark set, scored by the primal integral under a fixed
wall-clock limit, and its results are hosted on Mittelmann's
site~\cite{mittelmann_benchmarks}. Our only quantitative claim is a controlled comparison:
the selected configuration of the suite versus vanilla HiGHS, on one machine, under that
instance list, limit and metric. Our absolute numbers are not comparable to the rankings posted
on the site, which are measured on other hardware.
Stronger solvers exist. On the benchmark's own
published table the parallel solver ReXi~\cite{mexi2026rexi} leads the open-source field,
followed by the concurrent SCIP variant SCIPConc~\cite{mexi2026rexi,scipconc} and, just behind
it, the B200 run of the GPU solver cuOpt~\cite{cuopt2025}; the virtual mean of the three
commercial solvers is ahead of all of them. None is the comparison this paper makes, which is a
patched HiGHS against the HiGHS it patches.

The remainder of the paper is organized as follows. Section~\ref{sec:heuristics} walks through the
four heuristics and their integration into HiGHS, one heuristic at a time.
Section~\ref{sec:infra} describes the shared infrastructure: the presolve hook, effort and
patience, parallel execution and the solution pool, and instrumentation. Section~\ref{sec:results} presents the computational study stage by stage, each with its
results. Section~\ref{sec:conclusion} concludes.

\subsection{Related work}
\label{sec:related}

\paragraph{Primal heuristics in MIP.}
Primal heuristics are the solver components that find and improve incumbents inside
branch-and-bound. Berthold's thesis~\cite{berthold2006primal} catalogues the rounding, diving
and neighborhood heuristics of an early SCIP, Fischetti and Lodi~\cite{fischetti2011heuristics}
give a compact classification, and Berthold, Lodi and Salvagnin~\cite{berthold2025book} are the
current book-length treatment, whose grouping the paragraphs below follow. Berthold, Hendel and
Koch~\cite{berthold2018threephases} split a solve into feasibility, improvement and proof
phases; our chain runs entirely inside the first, before the root LP. What a heuristic costs
the solver it is added to depends on a second distinction, what it needs from the LP relaxation.
\emph{LP-free} heuristics never solve one and work on the constraint
matrix, the bounds and the integrality directly. \emph{LP-guided} heuristics need a fractional
point, and differ in who pays for it: one kind consumes a point the host has already computed by
simplex, while another solves its own relaxation by a first-order method, which is
\emph{matrix-free}, never factorizes, and returns an approximate point at low cost.

The primal integral, the area under the primal gap over
time, was introduced by Achterberg, Berthold and Hendel~\cite{achterberg2012rounding} and
studied as a measure by Berthold~\cite{berthold2013integral}; we use it in the benchmark's
normalized form (Section~\ref{sec:setup}). It rewards early incumbents, which is what a
presolve-window heuristic can offer, and it does not say whether the solve ends anywhere
different, which is why we report the final gap beside it.

\paragraph{Feasibility-first heuristics.}
Constructive heuristics start without an incumbent and aim for a first feasible point. The
feasibility pump~\cite{fischetti2005feaspump} alternates an LP projection with a rounding, the
objective pump~\cite{achterberg2007objfp} blends the original cost into the projection, and
\emph{Feasibility Pump 2.0}~\cite{fischetti2009fp2} replaces plain rounding by propagation-aware
rounding; Berthold, Lodi and Salvagnin~\cite{berthold2019tenyears} survey the decade of
descendants. Rounding and propagation heuristics work without any LP solve:
Shift-and-Propagate~\cite{berthold2015shiftpropagate} fixes variables one at a time, propagates,
and shifts to repair, and the structure-driven fix-and-propagate heuristics of Gamrath et
al.~\cite{gamrath2019fixpropagate} choose the fixings by their predicted propagation impact and
finish with an LP or a sub-MIP. FPR~\cite{salvagnin2025fpr} adds to that scheme backtracking and
a repair phase in the manner of WalkSAT~\cite{selman1994walksat}, the stochastic local search for
satisfiability that flips a variable in a violated clause, with noise to escape local minima, and
it runs a portfolio of ranking and value rules. RENS~\cite{berthold2014rens} is the exact
counterpart, a sub-MIP over the roundings of an LP point. Feasibility
Jump~\cite{luteberget2023fj} needs neither
an LP nor a propagation engine, only a weighted violation and one-variable moves, and it won the
MIP 2022 computational competition, whose brief was LP-free heuristics. HiGHS ships an
implementation, which we wrap and align with the paper in two places (Section~\ref{sec:fj}).

\paragraph{Neighborhood and improvement heuristics.}
Once an incumbent exists, large-neighborhood search (LNS) improves it by solving a sub-MIP.
RINS~\cite{danna2005rins} fixes the variables on which the incumbent and the LP solution agree;
local branching~\cite{fischetti2003localbranching} adds a Hamming-ball constraint around the
incumbent; DINS~\cite{ghosh2007dins} defines the neighborhood by distance from the LP solution;
solution polishing~\cite{rothberg2007polishing} recombines a population of incumbents; proximity
search~\cite{fischetti2014proximity} swaps the objective for a distance to the incumbent; and
Hendel's adaptive LNS~\cite{hendel2022alns} runs the LNS neighborhoods of SCIP under one
controller. All of them need an incumbent and most need an LP solution. HiGHS runs its own RINS and RENS in the dive, which is why our dive-time
heuristic (FPR LP-guided) shares its effort with them.

\paragraph{Local search from SAT to MIP.}
LocalMIP descends from the stochastic local search of satisfiability solving, with linear
constraints in place of clauses. Constraint weighting, in the additive form of PAWS~\cite{thornton2004paws}, raises
the weight of constraints that stay violated so that the landscape changes under the search, and
both FJ and LocalMIP use it. Walser~\cite{walser1999book} carried these ideas to integer
programs early on, with a domain-independent local search. The line from
Cai's group applies them to general MIP with tight and lift moves:
Local-MIP~\cite{lin2024localmip,lin2025localmipaij} is a standalone solver; on seven public
benchmarks including MIPLIB its authors report that it outperforms CPLEX, HiGHS, SCIP and FJ
and is competitive with Gurobi~\cite{lin2024localmip}.
ViolationLS~\cite{davies2024violationls} takes FJ's violation-driven moves into a
constraint-programming solver and reports gains from running it inside that solver's portfolio.

\paragraph{First-order LP solves inside heuristics.}
PDLP~\cite{applegate2021pdlp,applegate2026pdlp} made the primal-dual hybrid gradient practical
for LP through restarts, preconditioning and adaptive steps, and the GPU implementations
cuPDLP.jl~\cite{lu2025cupdlpjl} and cuPDLP-C~\cite{lu2023cupdlpc} made it fast. A first-order
solve is matrix-free and cheap per iteration but slow to converge, so it suits heuristics that
want many rough LP points rather than one exact one. Scylla~\cite{mexi2023scylla} builds a
feasibility pump on such solves and follows the objective pump's schedule. Kempke and Koch~\cite{kempke2026fixpropagate} show that low-precision PDLP
solutions lose nothing as the input to fix-and-propagate and use them on models of up to 243
million nonzeros; cuOpt~\cite{cuopt2025} fuses GPU PDLP with a feasibility pump, FJ and
fix-and-propagate; and CHAP~\cite{tjusila2026chap} coordinates the same ingredients across GPU
and CPU through a shared solution pool.

\section{The heuristic suite}
\label{sec:heuristics}

The suite integrates four published heuristics into HiGHS: Feasibility Jump, FPR, LocalMIP and
Scylla, with FPR dispatched at two points of the solve, so five mechanisms in all. They divide
along the LP axis of Section~\ref{sec:related}. Three are \emph{LP-free}: Feasibility Jump
(Section~\ref{sec:fj}), FPR (Section~\ref{sec:fpr}) and LocalMIP (Section~\ref{sec:localmip}).
Two are \emph{LP-guided}, and both round a fractional point with the same kernel. LP-guided FPR
(\texttt{fpr\_lp}, Section~\ref{sec:fprlp}) takes its points computed by simplex, which is why it runs during the
branch-and-bound dive rather than at presolve. Scylla (Section~\ref{sec:scylla}) solves its own
relaxations with the first-order method PDLP~\cite{applegate2021pdlp} and pumps, alternating a
solve with a rounding and moving the objective between them; it is matrix-free.

Three of the five share a rounding kernel, \texttt{fpr\_attempt} (Section~\ref{sec:fpr}): FPR
itself, LP-guided FPR and Scylla, which differ in the reference point they hand it. Feasibility
Jump and LocalMIP have their own move loops. The four presolve heuristics are dispatched as one
fixed chain through a single presolve hook, and all run under the shared budgeting, patience,
submission and instrumentation layer of Section~\ref{sec:infra}. The LP-guided FPR runs within the dive,
and its budget is shared with the host's own RINS and RENS.

\subsection{Feasibility Jump}
\label{sec:fj}

Feasibility Jump (FJ)~\cite{luteberget2023fj} is an LP-free Lagrangian heuristic. It keeps a weight
on every row and moves one variable at a time to its \emph{jump value}, the value that minimizes the
weighted violation of the variable's rows with the other variables held fixed. When no move helps,
the walk is at a local minimum, and it raises the weights of the violated rows until some move helps
again. Once the point is feasible, the same step raises the weight of the objective instead, so the
walk goes on toward better solutions. HiGHS ships an implementation of the core, and we run that core
rather than a new one. Algorithm~\ref{alg:fj} gives one worker.

\begin{algorithm}[t]
\caption{Feasibility Jump: one worker (HiGHS's core inside our wrapper)}
\label{alg:fj}
\begin{algorithmic}[1]
\State \textbf{Input:} model; start point $x^0$ (pool best, else the incumbent, else a finite bound); seed; per-worker budget $B = e\cdot(\text{nnz}\ll\ConstFjupstreamtotaleffortshift)$; patience $P$
\State build the core's problem: integer bounds rounded inward, each finite row side its own $\le$ or $\ge$ row, costs multiplied by the model sense \Comment{not charged to $B$, not interruptible}
\State $x \gets x^0$ clamped to the bounds;\; $w_i \gets 1$ for every row;\; $w_{\text{obj}} \gets 0$;\; $\Delta \gets 1$ with decay $1.0$ \Comment{\texttt{weightUpdateIncrement}, \texttt{weightUpdateDecay}}
\State for every $j$: $v_j \gets \textsc{Jump}(j)$, $s_j \gets \textsc{Score}(j, v_j)$;\; $G \gets \{j : s_j > 0\}$ \Comment{charges $\text{nnz}$ plus one column length per variable}
\Loop
  \If{effort since the last callback $> 500\,000$} \Comment{\texttt{CALLBACK\_EFFORT}}
    \State \textbf{callback}: \textbf{stop} if past the deadline, or the attempt's effort $\ge B/10$, or the total effort $> B$, or the effort since the core's last improvement $> P$ \label{ln:fj-callback}
  \EndIf
  \If{no row is violated \textbf{and} $c^\top x$ beats the best by a relative margin $> 10^{-4}$} \Comment{\texttt{kMinRelativeObjectiveImprovement}}
    \State record the best; \textbf{callback} with $x$: offer it to the pool at once, then the same stop tests \label{ln:fj-offer}
  \EndIf
  \If{$G \ne \varnothing$}
    \State with probability $0.001$ pick $j \in G$ uniformly; otherwise sample $\min(25, |G|)$ members of $G$ with replacement and take the largest $s_j$ \Comment{\texttt{randomVarProbability}, \texttt{maxMovesToEvaluate}}
  \Else \Comment{local minimum: bump the weights}
    \If{some row is violated}
      \State $w_i \gets w_i + \Delta$ for every violated row $i$; add the weighted change to $s_j$ for every $j$ in those rows; update $G$
    \Else
      \State $w_{\text{obj}} \gets w_{\text{obj}} + \Delta$;\; $s_j \gets s_j - \Delta\, c_j\,(v_j - x_j)$ for every $j$ \Comment{objective term subtracted; see text} \label{ln:fj-wobj}
    \EndIf
    \State \textbf{if} any weight $> 10^{20}$: scale every $w_i$, $w_{\text{obj}}$ and $\Delta$ by $10^{-20}$ and rescore every variable
    \If{some row is violated}
      \State pick a violated row $i$ uniformly; with probability $0.01$ pick $j$ in row $i$ uniformly, otherwise the $j$ in row $i$ with the largest $s_j$ \Comment{\texttt{randomCellProbability}}
    \Else
      \State pick $j$ uniformly among all variables
    \EndIf
  \EndIf
  \State $x_j \gets v_j$; update the row activities, the scores of every variable sharing a row with $j$, and $G$ \Comment{charges the rows' lengths}
  \State $v_j \gets \textsc{Jump}(j)$;\; $s_j \gets \textsc{Score}(j, v_j)$ \Comment{\texttt{resetMoves}}
\EndLoop
\Statex
\Function{Jump}{$j$} \Comment{\texttt{JumpMove::updateValue}}
  \For{each row $i$ containing $j$, the other variables held at $x$}
    \State $[l, u] \gets$ the values of $x_j$ that satisfy row $i$; \textbf{if} $a_{ij} < 0$ \textbf{then} swap $l$ and $u$ \label{ln:fj-swap} \Comment{our change; see text}
    \State round $[l, u]$ inward if $j$ is integer; skip the row if it is empty
    \State start the slope at $-w_i$ if $l$ lies above the lower bound, at $+w_i$ if $u$ lies at or below it; add $l$ and $u$ (when below the upper bound) as breakpoints, each raising the slope by $w_i$
  \EndFor
  \State add the finite bounds of $x_j$ as breakpoints; sort; walk the piecewise-linear weighted violation from the first breakpoint
  \State \textbf{return} the breakpoint of least violation other than $x_j$ itself, stopping once the slope is $\ge 0$
\EndFunction
\Function{Score}{$j, v$} \Comment{\texttt{resetMoves}}
  \State \textbf{return} $-\,w_{\text{obj}}\, c_j\,(v - x_j) \;+\; \sum_{i \ni j} w_i\,\big[\text{viol}_i(x) - \text{viol}_i(x \text{ with } x_j = v)\big]$ \label{ln:fj-score}
\EndFunction
\end{algorithmic}
\end{algorithm}

We align the implementation in HiGHS with the paper in two places. The first is the jump value
of a variable with a negative row coefficient (line~\ref{ln:fj-swap}). Dividing the row's sides
by a negative coefficient reverses the interval of satisfying values, and without the swap that
interval comes out empty and the row contributes nothing to the jump value. The second is the
sign of the objective term in the move score (lines~\ref{ln:fj-wobj} and~\ref{ln:fj-score}): the
score is an improvement, so an increase of the objective enters it negated.

Vanilla HiGHS calls its Feasibility Jump once, on one thread, with a budget of
$\text{nnz}\ll\ConstFjupstreamtotaleffortshift$ effort units. The patched solver calls the same
core, with the two changes above, from our concurrent runner in place of that single call: on
$N$ workers, each with its own budget.

\subsection{Fix, Propagate, and Repair (FPR)}
\label{sec:fpr}

Fix, Propagate, and Repair (FPR)~\cite{salvagnin2025fpr} is an LP-free constructive heuristic. A
depth-first search fixes the integer variables one at a time, in a chosen order and to a value
chosen by a rule, and propagates each fixing to tighten the remaining domains. When a fixing leaves
some row unsatisfiable, FPR repairs the partial assignment by shifting domains, and backtracks if the
repair fails. Once every integer is fixed, the other variables are filled in, and a local search,
WalkSAT~\cite{selman1994walksat} or the paper's RepairSearch, repairs whatever violation the complete
point still carries. A \emph{configuration} sets the variable order $\sigma$, the value rule $\nu$
and a mode $m$ that turns propagation, repair and backtracking on or off, so one kernel spans the
paper's family from a plain dive to a repairing tree search. The orders range from formulation order
and lock counts to clique covers from the HiGHS clique table; the value rules from the bounds, chosen
by objective direction or by locks, to an LP-based rule that rounds a reference point up with
probability equal to its fractional part. Algorithm~\ref{alg:fpr} gives one attempt of the kernel,
\texttt{fpr\_attempt}, and Algorithm~\ref{alg:fprworker} the presolve worker that drives it. The
kernel is shared: LP-guided FPR (Section~\ref{sec:fprlp}) and Scylla (Section~\ref{sec:scylla}) call
it with a reference point for the LP-based rule.

\begin{algorithm}[t]
\caption{FPR: one attempt of the kernel under one configuration}
\label{alg:fpr}
\begin{algorithmic}[1]
\State \textbf{Input:} model $(A,b,\ell,u,\mathcal{I})$ with $n$ columns; configuration $(\sigma,\nu,m)$; attempt index $a$; per-call effort slice $S$; optional reference point $x^\ast$
\State $\textit{prop} \gets (m \ne \textsf{dive})$;\; $\textit{rep} \gets m\in\{\textsf{dfsrep},\textsf{dive},\textsf{diveprop}\}$;\; $\textit{bt} \gets m\in\{\textsf{dfs},\textsf{dfsrep},\textsf{repairsearch}\}$
\State order $\gets$ precomputed order for $\sigma$; \textbf{if} $a>0$ \textbf{then} shuffle its first $\max(1,\lfloor 3n/10\rfloor)$ entries \Comment{Phase~1}
\State \textbf{if} $\nu=\textsf{loosedyn}$ \textbf{or} $\textit{rep}$ \textbf{then} initialise incremental row activity ranges
\ForAll{unfixed integer $j$} \Comment{trivially roundable, from HiGHS lock counts}
  \State \textbf{if} $\text{uplocks}_j=0$ \textbf{and} $\text{downlocks}_j>0$ \textbf{then} fix $x_j\gets u_j$
  \State \textbf{if} $\text{downlocks}_j=0$ \textbf{and} $\text{uplocks}_j>0$ \textbf{then} fix $x_j\gets \ell_j$
\EndFor
\State propagate to a fixpoint from every fixed column \Comment{each call capped at $100\cdot\text{nnz}$ accesses, \texttt{kPropagateBudgetPerNnz}}
\State $k \gets$ first unfixed integer (static order, or top of the domain-width heap for \textsf{domainsize})
\If{none} \State $\textit{complete}\gets\textbf{true}$
\Else \State $v\gets\nu(k)$; \textbf{if} $\textit{bt}$ \textbf{then} push $(k,\bar v)$; push $(k,v)$ \Comment{$\bar v$: the other bound; $1-v$ for binaries}
\EndIf
\While{stack nonempty \textbf{and} nodes $< n+1$ \textbf{and} \textbf{not} \textit{complete} \textbf{and} slice effort $< S$} \Comment{Phase~2; deadline polled every 16 nodes}
  \State pop $(j,v)$; restore bounds, values, activities and heap to the node's undo marks
  \State $\textit{infeas} \gets \neg\,\textsc{Fix}(x_j\gets v)$ \Comment{false only if $v$ is outside $[\ell_j,u_j]$}
  \State \textbf{if} $\neg\textit{infeas}$ \textbf{and} $\textit{rep}$ \textbf{then} $\textit{infeas}\gets$ some row of column $j$ is unsatisfiable by every completion of the current domains
  \State \textbf{if} $\neg\textit{infeas}$ \textbf{and} $\textit{prop}$ \textbf{then} $\textit{infeas}\gets(\textsc{Propagate}(j)=\textsf{infeasible})$ \Comment{a truncated fixpoint is not a refutation}
  \State \textbf{if} $\textit{infeas}$ \textbf{and} $\textit{rep}$ \textbf{then} $\textit{infeas}\gets\neg\,\textsc{RepairWalk}$ on the partial assignment \Comment{$\le 200$ shifts, noise $0.75$, tabu $3$, soft restart every $10$, $\le 100\cdot\text{nnz}+\text{rows}$ accesses}
  \State \textbf{if} $\textit{infeas}$ \textbf{and} $\textit{bt}$ \textbf{then} \textbf{continue} \Comment{the sibling pushed below this node is the backtrack}
  \State $k \gets$ next unfixed integer;\; \textbf{if} none \textbf{then} $\textit{complete}\gets\textbf{true}$; \textbf{break} \Comment{set even if $\textit{infeas}$ in a non-backtracking mode}
  \State $v\gets\nu(k)$; \textbf{if} $\textit{bt}$ \textbf{then} push $(k,\bar v)$; push $(k,v)$
\EndWhile
\State \textbf{if} slice exhausted with stack nonempty \textbf{then return} \textsf{paused} \Comment{state kept; the next call resumes here}
\ForAll{unfixed $j$} \Comment{Phase~2.5, charged $\text{nnz}$ for the row rebuild; the fill itself is not counted}
  \State continuous with cost: bound in the improving direction; costless: $0$; both clamped to a $\pm 10^{5}$ box when a bound is infinite
  \State integer: $\nu(j)$, rounded and clamped to $[\ell_j,u_j]$
\EndFor
\State compute all row activities; \textbf{if not} \textit{complete} \textbf{then return} the point as infeasible \Comment{no Phase~3}
\If{some row violated \textbf{and} $m=\textsf{repairsearch}$} \State \textsc{RepairSearch}: $\le 50$ nodes, stall threshold $10$, noise $0.75$ \Comment{Phase~3}
\ElsIf{some row violated \textbf{and} $\textit{rep}$} \State \textsc{WalkSAT}: $\le 200$ steps, noise $0.75$, $\le 100\cdot\text{nnz}$ accesses, restore the best state seen
\EndIf
\State re-check every row \Comment{the one verdict site; Phase~3's own answer is not trusted}
\State \textbf{if} feasible \textbf{then} greedy 1-opt: shift each costed integer by $\pm1$ toward the objective where no row breaks; \textbf{return} feasible point and objective
\State \textbf{else return} the point as infeasible
\end{algorithmic}
\end{algorithm}

\begin{algorithm}[t]
\caption{Presolve FPR: one call of worker $w$}
\label{alg:fprworker}
\begin{algorithmic}[1]
\State \textbf{Input:} slice $S=\max(\lfloor \text{total}/(10N)\rfloor,1)$; patience room $P$ (patience minus effort since this worker last improved the incumbent)
\State $C\gets\min(S,\max(P,1))$;\; spent $\gets 0$;\; started $\gets 0$
\While{spent $< C$ \textbf{and} deadline not passed}
  \If{no attempt is live}
    \State \textbf{if} started $=32$ \textbf{then break} \Comment{\texttt{kMaxAttemptsPerCall}}
    \State started $\gets$ started $+1$;\; config $\gets$ \texttt{kInitialFprConfigs}$[(w+a)\bmod 8]$;\; begin Algorithm~\ref{alg:fpr} with $(\text{config}, a)$
  \EndIf
  \State step Algorithm~\ref{alg:fpr} with slice $C-\text{spent}$;\; add its effort to spent
  \State \textbf{if} it returned \textsf{paused} \textbf{then} charge spent to the patience; \textbf{return}
  \State finish Algorithm~\ref{alg:fpr};\; add its effort to spent;\; offer a feasible point to the pool
  \State $a\gets a+1$ \Comment{next attempt rotates the configuration and reshuffles the order}
\EndWhile
\State charge spent to the patience (reset if the incumbent improved)
\end{algorithmic}
\end{algorithm}

\subsection{LP-guided FPR}
\label{sec:fprlp}

LP-guided FPR (\texttt{fpr\_lp}) realizes Classes~2--3 of~\cite{salvagnin2025fpr}: the kernel of
Algorithm~\ref{alg:fpr} with the order and the value rule driven by an LP point. It needs LP
solutions, so it runs inside branch-and-bound rather than at presolve: the patch adds one call to
HiGHS's dive-time heuristics, after RENS~\cite{berthold2014rens} and RINS~\cite{danna2005rins}. Each
dispatch uses three reference points: the LP solution at the current node, which the host has already
paid for, and the analytic center and a vertex of the relaxation with the objective set to zero, from
two fresh LP solves. Each worker is bound to one pairing of a configuration with one of the
three points, and rounds that point with the LP-based value rule.

\subsection{LocalMIP}
\label{sec:localmip}

LocalMIP~\cite{lin2024localmip,lin2025localmipaij} is an LP-free weighted local search. Like FJ it
moves one variable at a time and raises row weights at local optima, but it chooses among named
operators. The \emph{tight move} shifts a variable so that one of its rows becomes tight. The
\emph{breakthrough move} pushes a costed variable just past the best objective found so far. On a
feasible point the \emph{lift move} improves the objective as much as it can without breaking a row.
Candidates are ranked by weighted progress on the rows and the objective, with a bonus for moves that
make rows strictly satisfied or beat the best objective as the tie-break. A tabu list forbids undoing
a move for a few steps. The weights follow the probabilistic version of PAWS, the pure additive
weighting scheme from SAT~\cite{thornton2004paws}. At a local optimum, the weight of every violated
row rises by one, or the objective's weight when the point is feasible. With a small smoothing
probability the weights fall by one instead: those of the satisfied rows, and the objective's when
the point beats the best objective (line~\ref{ln:lm-paws}). Violated rows thus gain weight
over time without the weights growing apart without bound. Our implementation follows the paper's
Algorithms~1 and~2 and, where the paper is silent, the authors' public code~\cite{lin2024localmipcode}.
Algorithm~\ref{alg:localmip} gives one worker and Algorithm~\ref{alg:localmipcand} the
candidate generation it calls.

\begin{algorithm}[t]
\caption{LocalMIP: one worker}
\label{alg:localmip}
\begin{algorithmic}[1]
\State \textbf{Input:} model; start $x$; per-worker effort budget; a tolerance $\epsilon$
\State $w_i \gets 1$ for every row $i$;\; $w_{\text{obj}} \gets 1$;\; $x^\star \gets \varnothing$;\; $s \gets 0$ (moves since last improvement);\; $r \gets 0$ (random walks)
\State compute all row activities; $V \gets$ violated rows; $S \gets$ satisfied inequality rows
\While{budget remains \textbf{and} not stale \textbf{and} deadline not passed (polled every $65{,}536$ effort units)} \label{ln:lm-budget}
  \If{$V = \varnothing$} \Comment{feasible mode}
    \State on entry from infeasible mode and every \ConstLocalmiprecheckperiod{}th feasible step: recompute all activities and $V,S$ ($\text{nnz}$ effort); \textbf{if} $V \neq \varnothing$ \textbf{then continue}
    \If{$c^\top x < c^\top x^\star - \epsilon$} \State $x^\star \gets x$; offer $x$ to the pool; $s \gets 0$; $r \gets 0$ \EndIf
    \State recompute the lift bounds of the costed columns marked dirty (all of them after an infeasible episode)
    \State $m \gets$ the lift move of largest objective gain, tabu list not consulted
    \If{no such $m$} \State $m \gets \textsc{Candidates}(x, V, S)$ \Comment{Algorithm~\ref{alg:localmipcand} with $V = \varnothing$} \EndIf
  \Else \Comment{infeasible mode}
    \State $m \gets \textsc{Candidates}(x, V, S)$
  \EndIf
  \If{$m$ exists}
    \State apply $m$ to $x_j$, updating the activities of column $j$ and $V,S$; forbid the reverse direction of $x_j$ for $3 + \text{rand}(10)$ steps
    \State $s \gets s+1$; \textbf{if} the step was infeasible \textbf{and} now $V = \varnothing$ \textbf{then} $s \gets 0$
  \EndIf
  \If{the step was feasible \textbf{and} $s \ge 5000$} \Comment{feasible plateau} \label{ln:lm-plateau}
    \If{$r < 20$} \State perturb $x$ (each integer with prob.\ $0.2$); $w \gets 1$; rebuild activities; $r \gets r+1$; $s \gets 0$
    \Else \State retire \Comment{the runner rebuilds the worker from the pool}
    \EndIf
  \EndIf
  \If{$s \ge 200{,}000$} \State restart from $x^\star$ (odd restarts, when it exists) or from a random assignment; clear both tabu lists; $s \gets 0$ \label{ln:lm-restart} \EndIf
  \State every $100{,}000$ steps: rebuild activities and $V,S$ from scratch
\EndWhile
\State \textbf{return} $x^\star$
\end{algorithmic}
\end{algorithm}

\begin{algorithm}[t]
\caption{LocalMIP: candidate generation $\textsc{Candidates}(x, V, S)$ (paper Algorithm~2 with the authors' sampling caps; Phases~5 and~6 are ours)}
\label{alg:localmipcand}
\begin{algorithmic}[1]
\State \textbf{Input:} $x$, $V$, $S$, weights $w$, best objective $z^\star$ if one exists
\State $B \gets$ tight moves of the variables in the $12$ heaviest rows among $36$ rows sampled from $V$ (all of $V$ when $|V| \le 36$), at most $2250$ row entries examined \Comment{Phase~1, BMS}
\If{$z^\star$ exists \textbf{and} $c^\top x \ge z^\star - \epsilon$} \State add to $B$ the breakthrough move $x_j \gets x_j - (c^\top x - z^\star + \epsilon)/c_j$ of every costed variable, rounded and clamped \Comment{Phase~1b} \EndIf
\State $m \gets$ best of $B$ (tabu moves admitted only when they beat $z^\star$); \textbf{if} $\text{progress}(m) > 0$ \textbf{then return} $m$
\State $B \gets$ tight moves of the variables in one row sampled from $S$, at most $80$ row entries examined; $m \gets$ better of $m$ and best of $B$ (no aspiration); \textbf{if} $\text{progress}(m) > 0$ \textbf{then return} $m$ \Comment{Phase~2}
\State $B \gets$ flips of up to $5000$ binaries, a window at a random offset; $m \gets$ better of $m$ and best of $B$ (aspiration as in Phase~1); \textbf{if} $\text{progress}(m) > 0$ \textbf{then return} $m$ \Comment{Phase~3}
\State with prob.\ $1 - 3\times10^{-4}$: $w_i \gets w_i + 1$ for $i \in V$, or $w_{\text{obj}} \gets w_{\text{obj}} + 1$ when $V = \varnothing$;\; with prob.\ $3\times10^{-4}$: $w_i \gets \max(1, w_i - 1)$ for $i \in S$, and $w_{\text{obj}} \gets \max(1, w_{\text{obj}} - 1)$ when $c^\top x < z^\star$ \Comment{Phase~4, PAWS} \label{ln:lm-paws}
\State $B \gets$ tight moves of the variables in one row sampled from $V$ (no cap); $m \gets$ better of $m$ and best of $B$ (no aspiration); \textbf{if} $m$ exists \textbf{then return} $m$, whatever its sign
\State pick one variable of one row sampled from $V$: flip a binary, shift a general integer by $\pm 1$, or move a continuous variable by up to $\pm 10\%$ of its range; \textbf{if} it moves \textbf{then return} it \Comment{Phase~5, ours} \label{ln:lm-phase5}
\State $B \gets$ for $5$ random variables: the value nearest zero, each finite bound, and the midpoint (continuous only); \textbf{return} best of $B$ (no aspiration; may be none) \Comment{Phase~6, ours} \label{ln:lm-phase6}
\end{algorithmic}
\end{algorithm}

Our additions follow from one change of setting. The paper runs a single search to a time cutoff; a worker
here runs under an effort budget (line~\ref{ln:lm-budget}) and should not spend it stuck. On a cold
start, with an empty pool, the paper's zero assignment is refined by a greedy sweep: in order of
decreasing column length, each variable moves to whichever of a few candidate values (its bounds,
zero, tight moves for its violated rows) satisfies the most rows net of those it breaks. When
candidate generation still has no move after the weight update, two fallback phases supply one: a
random step on a variable of a violated row (line~\ref{ln:lm-phase5}), and failing that the best of a
few simple values for a handful of random variables (line~\ref{ln:lm-phase6}). A worker stuck on a
feasible plateau perturbs its point and resets its weights, and after a bounded number of such walks
it retires so that the runner rebuilds it from the pool (line~\ref{ln:lm-plateau}). A longer run
without improvement restarts the worker, alternating between the best point found and a random
assignment (line~\ref{ln:lm-restart}). The authors' code also restarts after a run without
improvement, from the best point or a random one; the alternation is ours.

\subsection{Scylla}
\label{sec:scylla}

Scylla~\cite{mexi2023scylla} is a feasibility pump~\cite{fischetti2005feaspump,achterberg2007objfp,fischetti2009fp2}
whose projection step is an approximate LP solve by the first-order method PDLP~\cite{applegate2021pdlp}
and whose rounding step is fix-and-propagate. Each round solves the relaxation under an objective
that blends the original cost with the $\ell_1$ distance to the last rounded point, rounds the LP
point with the FPR kernel, and shifts the blend further toward distance, as in the objective
feasibility pump~\cite{achterberg2007objfp}. The LP tolerance tightens from round to round, so early
solves are cheap. When a rounding repeats a recent one, the pump perturbs it to break the cycle.
Scylla is matrix-free in its reference's sense, since PDLP never factorizes, and LP-guided in this
paper's split, since it solves relaxations. Algorithm~\ref{alg:scylla} gives one worker.

\begin{algorithm}[t]
\caption{Scylla: one worker sharing a PDLP instance}
\label{alg:scylla}
\begin{algorithmic}[1]
\State \textbf{Input:} model $(A,b,\ell,u,\mathcal{I})$; dispatch budget $B$; worker $w$ of $N$; shared PDLP instance $P$ (mutex; latest snapshot with generation $g$)
\State $K \gets 0$;\; $\alpha \gets 1$;\; $\varepsilon \gets 10^{-2}$;\; $\tilde c \gets c$;\; $\sigma \gets \sqrt{|\mathcal{I}|}/\lVert c\rVert_2$;\; history $H \gets \varnothing$;\; $s \gets 0$ \Comment{$s$: consecutive stale rounds}
\State $\bar s \gets \min\{16,\, 4 + \lfloor \text{nnz}/83\,000 \rfloor\}$ \Comment{\texttt{compute\_max\_stale\_rounds}}
\State (order, mode) $\gets$ entry $w$ of $\{$(clique, DFS), (locks, DFS), (locks, dive), (formulation, DFS+repair)$\}$ for $w<4$, seed-chosen otherwise; value rule: LP-based \Comment{\texttt{kFprConfigs}}
\While{attempt slice and $B$ remain \textbf{and} deadline not passed}
  \If{a peer improved the incumbent since last check} \State reset this worker's staleness counter \EndIf
  \If{effort since last incumbent improvement $> k_{\text{scylla}}$} \State retire; the runner rebuilds the worker with a fresh seed; \textbf{stop} \EndIf
  \If{$s \ge \bar s$} \State $\bar x \gets \textsc{Solve}(P, \tilde c, \varepsilon)$, blocking on the mutex \Comment{forced fresh}
  \ElsIf{try-lock on $P$ succeeds} \State $\bar x \gets \textsc{Solve}(P, \tilde c, \varepsilon)$ \Comment{fresh}
  \Else \State $\bar x \gets$ latest snapshot of $P$; \textbf{if} none, or $g$ unchanged since this worker last read it \textbf{then} $s \gets s+1$; charge $1$; \textbf{continue} \Comment{stale}
  \EndIf
  \State \textsc{Solve}: set costs $\tilde c$; warm-start from this worker's own last primal--dual pair; \texttt{kkt\_tolerance} $\gets \varepsilon$; iteration limit $\max\{100, \lfloor B/(4\,\text{nnz}) \rfloor\}$; time limit $\gets$ remaining deadline; publish the result as a new snapshot; retire the worker on an error, an infeasible LP, or $3$ consecutive zero-iteration solves \Comment{\texttt{kMaxPdlpStalls}}
  \State charge $\text{iters}\times\text{nnz}$ to the attempt and $\text{iters}\times\text{nnz}/N$ to this worker's counters \Comment{$0$ on a stale round}
  \State $s \gets 0$ if fresh, else $s \gets s+1$
  \If{$\bar x$ is integral within the MIP tolerance \textbf{and} satisfies every row} \State offer $\bar x$ to the pool; \textbf{continue} \Comment{fast path}
  \EndIf
  \State $\hat x \gets \textsc{FprAttempt}(\bar x)$ with reference point $\bar x$: integer $x_j$ rounds up with probability $\mathrm{frac}(\bar x_j)$, zero-cost continuous $x_j \gets \bar x_j$; charge its effort \Comment{Section~\ref{sec:fpr}}
  \If{$\hat x$ is feasible} \State offer $\hat x$ to the pool \Comment{the pump continues either way}
  \EndIf
  \If{fresh} \Comment{pump state advances only on fresh rounds}
    \If{$\hat x$ equals some point in $H$ on every integer column} \State perturb $\hat x$: each integer column, with probability $0.2$, moves to a uniformly random other value in its domain
    \EndIf
    \State write $\hat x$ into slot $K \bmod 3$ of $H$ \Comment{$|H| \le 3$, \texttt{kCycleWindow}}
    \State $\alpha \gets 0.9\,\alpha$;\; $\tilde c_j \gets \alpha\,\sigma\,c_j + (1-\alpha)\,\delta_j$ for $j\in\mathcal{I}$, $\tilde c_j \gets \alpha\,\sigma\,c_j$ otherwise \Comment{\texttt{kAlpha}}
    \Statex \hspace{\algorithmicindent}\hspace{\algorithmicindent} with $\delta_j = 1 - 2\hat x_j$ if $[\ell_j,u_j]=[0,1]$, else $\delta_j = \operatorname{sign}(\bar x_j - \hat x_j)$
    \State $\varepsilon \gets \max\{0.98\,\varepsilon,\, 10^{-8}\}$;\; $K \gets K+1$ \Comment{\texttt{kBeta}, \texttt{kEpsilonFloor}}
  \EndIf
\EndWhile
\end{algorithmic}
\end{algorithm}

We use HiGHS's own PDLP solver through its \texttt{solver=pdlp} option, which in
\ConstHighsVersion{} runs cuPDLP-C~\cite{lu2023cupdlpc} (not HiPDLP). 
The $N$ workers share this
\emph{single} instance behind a lock, since a solve is expensive. A worker that takes the lock solves its own objective and publishes the result. A worker
that does not rounds the latest published solution instead of waiting, and after too many such stale
rounds it waits for a fresh solve of its own. Pump state advances only on fresh rounds, and any
worker's improvement of the incumbent resets the staleness of all of them. Effort counts PDLP
iterations times nonzeros plus the rounding's coefficient accesses. A worker whose effort since the
last improvement exceeds the patience retires and is rebuilt with a fresh seed.

In our implementation a feasible rounding does not end the pump: the point is offered to the pool 
and the pump goes on, 
since the pool keeps improving solutions. The rounding orders the variables once per dispatch rather
than by the current LP point; only the value rule reads it. The shared solver and the stale-snapshot
rounding are ours.

PDLP is well suited to a GPU, and HiGHS can build cuPDLP-C for one. In a preliminary test a GPU
build gave no advantage on this benchmark, so the campaign runs Scylla on the CPU
(Section~\ref{sec:machine}); we did not study the GPU build further.

\section{Integration infrastructure}
\label{sec:infra}

The heuristics share a thin layer between HiGHS and their kernels. It decides where
they are called, how much they may spend, how their workers run in parallel, and how their solutions
reach the solver. This section describes that layer at the level the experiments need; the code is
the reference for the rest.

\subsection{The presolve hook}
\label{sec:chain}

HiGHS calls the layer once presolve is complete and before the root LP is solved
(\textsc{PresolveHook} in Algorithm~\ref{alg:runner}). The hook runs every presolve heuristic whose
effort option is positive, one after another in the fixed order FJ $\to$ FPR $\to$ LocalMIP $\to$
Scylla, and each heuristic runs its own workers in parallel. A zero effort is the only way to exclude
a heuristic, so any subset of the five is a zero pattern of five continuous options. At the shipped
defaults Scylla's effort is \ConstEffortscylla{} (Section~\ref{sec:search}), so the default chain is
FJ $\to$ FPR $\to$ LocalMIP. One solution pool, seeded with the incumbent, serves the whole chain, so
a later heuristic can start from what an earlier one found. Because the chain is sequential, the time
limit truncates its tail: if an early heuristic uses the whole limit, the later ones never run.
LP-guided FPR is called from the dive instead (Section~\ref{sec:fprlp}), through the same runner. The
layer reads HiGHS options and never writes them; it changes solver state only by submitting solutions
and, for \texttt{fpr\_lp}, by charging LP iterations.

\subsection{Effort and patience}
\label{sec:budget}

Each presolve heuristic has two options that control its spend, an effort and a patience. Both
are multiples of $\text{nnz}\ll\ConstFjupstreamtotaleffortshift$, the budget vanilla HiGHS gives its
own Feasibility Jump, so they scale with model size and an effort of one is one such budget.

The \emph{effort} $e_h$ sets the budget $T_h = e_h\,(\text{nnz}\ll\ConstFjupstreamtotaleffortshift)$
that the heuristic's workers share in one dispatch; only FJ gives each worker the full budget
(Algorithm~\ref{alg:runner}). Each heuristic counts effort in its own unit (FJ in the effort its
core charges, FPR and LocalMIP in coefficient accesses, Scylla in PDLP iterations times nonzeros),
so equal values for two heuristics are not equal work.

The \emph{patience} $k_h$ stops a heuristic that spends without improving the best objective. A
counter of the effort charged since the last improvement is reset to zero at every improvement,
and the dispatch ends as soon as it reaches $\min\{k_h,\ T_h/4\}$; each worker of FJ, FPR and
LocalMIP also stops once its own counter reaches its share of that limit. A heuristic that keeps improving therefore runs until
its budget is spent, and one that stops improving runs for at most the patience more. The cap at a
quarter of the budget, the ratio HiGHS's own Feasibility Jump uses, lets the patience end a
dispatch before the budget does. $k_h=0$ means no patience limit.

The shipped values are those of the selected configuration (Section~\ref{sec:search}): efforts \ConstEffortfj,
\ConstEffortfpr, \ConstEffortlocalmip{} and \ConstEffortscylla{} and patiences \ConstPatiencefj,
\ConstPatiencefpr, \ConstPatiencelocalmip{} and \ConstPatiencescylla{} for FJ, FPR, LocalMIP and
Scylla. At these values FJ has no patience limit, FPR's lies below the cap, and LocalMIP's is the
cap. Table~\ref{tab:arms} lists them beside the other configurations the search compared.

\subsection{Parallel execution and the solution pool}
\label{sec:modes}
\label{sec:pool}

\begin{algorithm}[t]
\caption{Presolve hook, runner, and the heuristic's callbacks}
\label{alg:runner}
\begin{algorithmic}[1]
\Procedure{PresolveHook}{} \Comment{after presolve, before the root LP}
  \State $\mathcal{P} \gets$ solution pool holding the incumbent, if any \Comment{one pool for the chain}
  \For{$h$ in FJ, FPR, LocalMIP, Scylla \textbf{with} $e_h > 0$}
    \State \textbf{if} HiGHS has terminated \textbf{then break}
    \State $T \gets e_h\,(\text{nnz}\ll 10)$;\; $K \gets \min\{k_h\,(\text{nnz}\ll 10),\ T/4\}$, or $\infty$ if $k_h = 0$ \Comment{both times $N$ for FJ}
    \State $h$'s setup;\; \textsc{Runner}$(h,\ T,\ C = \max\{\lfloor T/(10N)\rfloor, 1\},\ K)$ \label{ln:runner-hook}
  \EndFor
\EndProcedure
\Statex
\Procedure{Runner}{$h, T, C, K$} \Comment{\texttt{fpr\_lp} calls it from the dive}
  \State shared atomics: $E \gets 0$;\; $E_{\text{stale}} \gets 0$;\; $\textit{stop} \gets \textbf{false}$;\; $\textit{seat} \gets 0$
  \ForAll{slots $w = 0, \dots, N-1$, \textbf{in parallel}}
    \State $\textit{rng} \gets$ base seed $+\,w$;\; $a \gets 0$;\; $S \gets h.\textsc{MakeState}(w)$
    \While{\textbf{not} \textsc{ShouldStop}$(w, a)$}
      \State $c \gets \min\{C,\ T - E\}$;\; \textbf{if} $c = 0$ \textbf{then} $\textit{stop} \gets \textbf{true}$; \textbf{break} \label{ln:runner-cap}
      \State $(\delta, \textit{improved}) \gets h.\textsc{Attempt}(S, c)$;\; $a \gets a + 1$
      \State \textbf{if} $\delta = 0$ \textbf{then break} \label{ln:runner-guard} \Comment{this slot retires; the others run on}
      \State $E_{\text{stale}} \gets 0$ \textbf{if} \textit{improved} \textbf{else} $E_{\text{stale}} + \delta$;\; $E \gets E + \delta$;\; \textbf{if} $E_{\text{stale}} \ge K$ \textbf{or} $E \ge T$ \textbf{then} $\textit{stop} \gets \textbf{true}$
    \EndWhile
    \State \textbf{if} $\textit{seat} = w$ \textbf{then} $\textit{seat} \gets \textsf{vacant}$
  \EndFor
\EndProcedure
\Statex
\Function{ShouldStop}{$w, a$}
  \State \textbf{if} past the deadline \textbf{then} $\textit{stop} \gets \textbf{true}$ \Comment{any slot}
  \If{$a$ is even \textbf{and} ($\textit{seat} = w$ \textbf{or} CAS $\textit{seat}\!: \textsf{vacant} \to w$)} \Comment{seat holder only}
    \State \textbf{if} HiGHS has terminated \textbf{then} $\textit{stop} \gets \textbf{true}$
  \EndIf
  \State \textbf{return} \textit{stop}
\EndFunction
\Statex
\Function{$h$.Attempt}{$S, c$} \Comment{the shape all five callbacks share}
  \State \textbf{if} $S$ has no live worker \textbf{then} $S.\textit{worker} \gets$ new worker from \textsc{StartPoint}$_h$, with a fresh share of $T$ \Comment{not FPR}
  \State run $S.\textit{worker}$'s own loop (Section~\ref{sec:heuristics}) until it charges $c$ or retires; pass each feasible point to \textsc{Offer} \label{ln:runner-call}
  \State \textbf{if} it charged $0$ and retired \textbf{then} rebuild and run again \Comment{Scylla, \texttt{fpr\_lp}} \label{ln:runner-retry}
  \State \textbf{return} (effort charged, whether any offered point improved the best objective)
\EndFunction
\Statex
\Function{Offer}{$x$} \Comment{thread-safe}
  \State admit $x$ if $\mathcal{P}$ holds fewer than ten points, or $x$ beats the worst, or $x$ is within a tenth of the best objective and differs from every stored point on at least a twentieth of the integer columns; evict the worst, or in the last case the most similar
  \State \textbf{if} admitted \textbf{then} submit $x$ to HiGHS now, tagged with $h$
  \State \textbf{return} whether $x$ beat every point $\mathcal{P}$ has admitted
\EndFunction
\Statex
\Function{StartPoint$_h$}{} \Comment{empty pool: the incumbent, else $h$'s own construction}
  \State FJ, and LocalMIP's first build: the best point of $\mathcal{P}$;\; FPR, Scylla, \texttt{fpr\_lp}: none
  \State LocalMIP's rebuild: pick two distinct points of $\mathcal{P}$ and draw $u$ uniformly, then perturb the result:
  \State \hspace{\algorithmicindent} $u < 2/5$: keep the integer columns they agree on; take every other column from either at random
  \State \hspace{\algorithmicindent} $u < 7/10$: keep the better point; take the integer columns where they differ from either at random
  \State \hspace{\algorithmicindent} otherwise: copy one point, from the better half of $\mathcal{P}$ with probability $1/2$
\EndFunction
\end{algorithmic}
\end{algorithm}

The runner is \emph{opportunistic}: the $N$ workers of a heuristic run continuously and share state
without ever synchronizing, so the order in which their work lands, and with it the result, can
differ from one run to the next. Algorithm~\ref{alg:runner} gives the whole layer.
\textsc{PresolveHook} runs the chain and \textsc{Runner} runs one heuristic's $N$ slots in parallel.
A heuristic plugs into the runner with two callbacks. \textsc{MakeState} sets up a slot, and
\textsc{Attempt} runs the heuristic: it advances the slot's worker, one of the
algorithms of Section~\ref{sec:heuristics}, by at most $c$ effort units, and replaces the worker with
a new one from \textsc{StartPoint} once it has retired. Between attempts the runner checks the
deadline and the solver's terminator, and charges the effort against the budget and the patience.

The pool (\textsc{Offer}) passes every admitted point to HiGHS immediately, so an incumbent's
timestamp is its discovery time, and tags it with the heuristic that found it (\texttt{J}=FJ,
\texttt{A}=FPR, \texttt{D}=LP-guided FPR, \texttt{M}=LocalMIP, \texttt{G}=Scylla), from which
attribution is read. The pool also gives rebuilt workers their start points (\textsc{StartPoint}).
Only FJ and LocalMIP read one: the FPR kernel, which Scylla and \texttt{fpr\_lp} also use, writes every
column before reading it (Section~\ref{sec:fpr}).

\subsection{Instrumentation: effort usage and accepted solutions}
\label{sec:instr}

Two log records carry everything the evaluation in Section~\ref{sec:results} is built on:
\begin{itemize}
\item \texttt{[Heur]}, one per dispatch with its effort usage: heuristic, phase (presolve or
  dive), start and end on the solver's clock, charged effort, wall time, and whether the pool
  accepted anything;
\item \texttt{[HeurSol]}, one per \emph{accepted} solution: the charged effort and
  wall time at which it happened, from which productive effort (effort at the last accepted
  solution), stale effort (the remainder) and the inter-acceptance effort gaps per nonzero are
  derived. These are the quantities the presolve probe of Section~\ref{sec:probe} reports.
\end{itemize}
A patched build also gains an option, \texttt{mip\_heuristic\_presolve\_only}, that exits the
solve after the chain and before the root LP with the incumbent preserved.

\section{Computational study}
\label{sec:results}

The evaluation runs in five stages on one machine (Table~\ref{tab:stages}). It follows the
good practice Beiranvand, Hare and Lucet~\cite{beiranvand2017bestpractices} set out for comparing
optimization algorithms: the metric is the benchmark's own and not one of our choosing, tuning and
evaluation are kept apart where the benchmark allows it, and every comparison is paired per
instance and reported with its confidence interval and the numbers of instances won and lost. We
call each solver configuration in a comparison an \emph{arm}. The stages use
a tuning set, a confirmation set drawn from it, and a held-out complement with a held-out check
drawn from it, all stratified by the baseline's time to first feasible solution
(Sections~\ref{sec:profile} and~\ref{sec:search}).

Stage~1 runs vanilla HiGHS over the benchmark as the baseline. Stage~2, the probe, runs each
heuristic alone to see what it reaches and how it spends its effort. Stage~3 searches the effort
budgets and patiences of the four heuristics jointly and confirms the finalists in full solves.
Stage~4 measures the LP-guided dive, which the earlier stages cannot see, and Stage~5 compares
the selected configuration with the baseline over the whole benchmark. Stages~2 and~3 score
heuristics in \emph{presolve-only} runs: the solve stops when the presolve heuristics return,
before the root LP, so a run is short and sees the presolve heuristics alone. Every accepted
solution carries a source tag naming the heuristic or HiGHS component that found it, and every
count of what a heuristic produced comes from these tags.

\subsection{Benchmark, metrics and machine}
\label{sec:setup}

We evaluate on the \texttt{mipfeas} benchmark~\cite{bussieck2026mipfeas}, whose results are
hosted on Mittelmann's benchmark site~\cite{mittelmann_benchmarks}: the 233 feasible
MIPLIB~2017~\cite{miplib2017} benchmark instances (the 240-instance set without its 7 infeasible
instances), each with a wall-clock limit of $600$\,s. We use the instance list and the time limit
unchanged. The reference objective $z^\ast$ is the benchmark's own, the file of MIPLIB~2017
objective values its published scoring script reads.

We score runs with the benchmark's own formula, computed as its published script computes it. The primal gap of an incumbent $z$ is
\begin{equation}
p \;=\; \begin{cases}
1 & \text{if } z\,z^\ast < 0,\\[2pt]
\dfrac{|z - z^\ast|}{\max(|z|,\,|z^\ast|,\,1)} & \text{otherwise,}
\end{cases}
\label{eq:gap}
\end{equation}
with $p=2$ before the first incumbent. The primal integral~\cite{berthold2013integral}, in the
benchmark's normalized form~\cite{bussieck2026mipfeas}, is $\frac{1}{T}\int_0^T p(t)\,dt$
over the time limit $T$, so it lies in $[0,2]$. For values
$x_1,\dots,x_n$ and shift $s>0$, the shifted geometric mean~\cite{achterberg2007thesis} is
\begin{equation}
\operatorname{SGM}_s(x) \;=\; \exp\!\left(\frac{1}{n}\sum_{i=1}^{n}\ln(x_i+s)\right) - s .
\label{eq:sgm}
\end{equation}
We report four metrics.
\begin{itemize}
\item \textbf{Primal integral}, $\operatorname{SGM}_{\ConstShiftHeadline}$ over all
  \ConstNInstances{} instances, the aggregate the benchmark's results page reports. This is the
  headline. It rewards reaching good solutions early.
  An instance with no feasible solution scores $2$, so failures stay in the set.
\item \textbf{Final gap}: the gap $p$ to the reference best known solution at the time limit,
  $\operatorname{SGM}_{\ConstShiftHeadline}$, with the number of instances on which each arm ends
  with the strictly better objective. This is the solution a user who waits for the limit gets,
  so reaching a good solution sooner is an improvement only if the final solution is not worse.
  Paired, it is compared over the instances both arms made feasible. A sub-MIP heuristic of
  HiGHS that the time limit cuts off hands its best solution back as the solve closes, so that
  solution can be stamped a fraction of a second past the limit. It is the run's final solution: we count it, up to one
  second past the limit, in the final gap and in the number of feasible instances, and we
  allow the same second at the caps of the presolve-only runs. The primal
  integral runs to the limit and is unaffected.
\item \textbf{Time to first feasible solution}, $\operatorname{SGM}_1$ in seconds over all
  instances, a run with no solution counting at the limit, $600$\,s.
\item \textbf{Number of instances with a feasible solution} at $600$\,s.
\end{itemize}

\paragraph{Paired comparisons.}
\label{sec:stats}
For two arms with values $x_i$ and $y_i$ of a metric with shift $s$ on instances
$i=1,\dots,n$, let $d_i=\ln\bigl((x_i+s)/(y_i+s)\bigr)$, with mean $\bar d$ and standard deviation
$\sigma_d$. We report the ratio $r=e^{\bar d}$, its 95\% confidence interval
$e^{\bar d\pm1.96\,\sigma_d/\sqrt{n}}$, the two-sided $p$-value of $\bar d\sqrt{n}/\sigma_d$
under the normal distribution, the decrease $1-r$, and the numbers of instances with $d_i<0$ and
$d_i>0$. Two arms \emph{separate} when this interval excludes $1$. There is one run per arm and
instance, and the result of a run can differ from one run to the next
(Section~\ref{sec:modes}). That noise is part of $\sigma_d$, so it widens the interval, and we
read aggregates only, never a single instance. The $p$-values are not adjusted for the number of
comparisons, and we treat those near the five per cent level as borderline.

\paragraph{Machine and build.}
\label{sec:machine}
All runs use one benchmark machine (AMD Ryzen~9 3950X, 16~cores / 32~threads, so HiGHS's default
of 16 workers), with threads unpinned, and every run used seed $0$. The patched solver~\cite{mip_heuristics_software} is HiGHS
\texttt{v1.15.1} built with the CPU PDLP. The baseline is a \emph{separately built, unpatched}
HiGHS at the same tag, identified by the absence of the patch marker line since the version
banners are identical.

\input{stages}

\subsection{Stage 1: baseline and the benchmark's difficulty profile}
\label{sec:profile}

Stage~1 runs vanilla HiGHS once over all \ConstNInstances{} instances. That run is the baseline
for every later comparison, and Table~\ref{tab:profile} shows when it first found a feasible
solution: \ProfUnderTenS{} of the \ProfTotal{} instances (\ProfUnderTenSPct\%) within ten
seconds, and \ProfNever{} never. This bounds what a presolve heuristic can
contribute: on most instances the most it can buy is a fraction of the first ten seconds of the
integral. The easy instances matter as much as the hard ones. Presolve effort buys
feasibility where feasibility is hard and is overhead where branch-and-bound has an
incumbent in the first second, so a set drawn from the hard end alone would see the benefit
without the cost. Every instance set of the later stages is therefore a draw stratified on the
strata of Table~\ref{tab:profile}, so that it samples the whole difficulty range.

\input{profile}

\subsection{Stage 2: presolve probe}
\label{sec:probe}

Stage~2 sets the default effort and patience of each heuristic from what it does alone. Each
runs by itself, presolve-only, over all \ProfTotal{} instances, with unlimited effort, no
patience limit and a $30$\,s wall-clock cap, so the clock is the only stopping rule. The four
run separately, because in a chain a slow first heuristic uses the whole cap.

From these runs, the effort $e_h$ is the median, over the runs that stopped improving before the
cap, of the effort at the last improvement. The patience is $k_h=\min\{q_{95},\ e_h/4\}$, where
$q_{95}$ is the 95th percentile of the effort between improvements. For all four heuristics
$q_{95}$ lies above $e_h/4$, FPR's by a factor of \PatienceFprRatio{}, so every $k_h$ is
$e_h/4$. Table~\ref{tab:effort} gives the resulting \emph{measured vector} with what each
heuristic reached, and Figure~\ref{fig:probe-time} shows when its improvements arrive. Every
heuristic spends far more effort after its last improvement than before it, which is the tail
the patience is meant to cut. Stage~3 centres its search ranges on the measured vector and
keeps it as the control.

\paragraph{Speed and quality.}
Running alone, LocalMIP solves the most instances, and FJ, FPR and LocalMIP are about equally
often the first to find a solution (Table~\ref{tab:effort}). Scylla is first far less often,
which supports placing it last but does not order the other three. They differ in quality. On
the instances two of them both solve, FJ ends better than FPR on \ProbeFjVsFprFirstBetter{} and
worse on \ProbeFjVsFprSecondBetter{}, better than LocalMIP on \ProbeFjVsLocalMipFirstBetter{}
and worse on \ProbeFjVsLocalMipSecondBetter{}, and LocalMIP ends better than FPR on
\ProbeLocalMipVsFprFirstBetter{} and worse on \ProbeLocalMipVsFprSecondBetter{}. FJ does not
trade quality for speed. Each heuristic also solves instances no other one does: LocalMIP
\ProbeLocalMipSoleCracker{}, FJ \ProbeFjSoleCracker{}, FPR \ProbeFprSoleCracker{} and Scylla
\ProbeScyllaSoleCracker{}.

Their order FJ, FPR, LocalMIP, Scylla was initially fixed on the rule \emph{cheapest likely
success first}: our best guess from the speed and solution quality each heuristic's own paper
reports. The probe gives no reason to change it. Among heuristics that are equally fast, the
first to run can be expected to find most first solutions, so what a later heuristic adds is
what it does after FJ, which Section~\ref{sec:headline} measures.

\input{effort}
\input{fig/probe-time}

\subsection{Stage 3: joint search}
\label{sec:search}

\paragraph{Race.}
Stage~3 tunes eight numbers at once, an effort $e_h$ and a patience $k_h$ for each heuristic.
Effort zero switches a heuristic off, so the choice of heuristics is part of the same search.
Iterated racing (\texttt{irace}~\cite{lopezibanez2016irace}) runs presolve-only over the tuning
set, each run stopping at the presolve exit or at a $60$\,s cap, and scores it by its
presolve-exit gap plus $\lambda\tau$, where $\tau$ is the heuristics' own wall time. The weight
is of the order of what a second of delay costs the primal integral, between $1/T$ and $2/T$
since the gap is at most $2$. The race ran at $\lambda=1/600$ and at half and twice it, and once
more with all four heuristics forced on; \SearchCandidatesN{}
candidates in all. Racing ends with a set of survivors, which are re-scored by the SGM of the
presolve-exit gap; ties go to the simpler configuration, with fewer heuristics and then less
total effort. Finalists then run full solves on a confirmation set and a held-out check.

\paragraph{Instance sets.}
The \emph{tuning set} is \TuningN{} instances, drawn after the probe from the \ProbePoolN{}
instances on which some heuristic running alone found a solution: on the others the
presolve-only race scores every configuration alike. Its size is set by the race, which ranks
many configurations on few instances. The race's selection noise falls with the square root of the
number of instances, and \TuningN{} was chosen to bring it below the differences between
configurations the race has to rank; an earlier draw of 25 did not. The other \HeldoutN{}
instances form the \emph{held-out complement}, which therefore holds every instance no
heuristic solved alone. The confirmation set is \SearchConfirmInstancesN{} instances of the
tuning set, and the held-out check \SearchHeldoutInstancesN{} of the complement, the only
complement instances any stage sees before the headline. Both sizes come from a power
calculation on the paired spread of the first confirmation runs: at that spread about fifty
instances detect a twenty per cent difference four times in five. The Set column of
Table~\ref{tab:per-instance} marks the instances of every set.

\paragraph{Race outcome.}
Table~\ref{tab:arms} lists the configurations that ran full solves, under the names of the code
repository~\cite{mip_heuristics_software}; the text uses the descriptions beside them. The
tie-break selected FJ alone at $\lambda=1/600$ and FJ, FPR and LocalMIP at the other two
weights. The selection is not monotone in $\lambda$, so the cost weight does not decide it. Of the free
races' survivors, \SearchMixSurvivorsN{} run FJ, FPR and LocalMIP, the rest FJ alone, and none
Scylla. The mix at $\lambda=1/1200$ nearly duplicates the one at $1/300$ and was not run; the
tie-break applied once to the \SearchMixSurvivorsN{} pooled three-heuristic survivors selects a
cheaper mix. Both mixes went to confirmation, with FJ alone, the measured vector and the winner
of the forced race.

\paragraph{Confirmation and held-out check.}
The five finalists ran \ConstTimeLimitS\,s solves on the confirmation set
(Figure~\ref{fig:forest}). FJ alone was nominally worst, borderline against the generous mix
($p=\SearchConfirmFjOnlyVsMixP$), and was dropped. The other four did not separate, so the
tie-break decided: fewer heuristics keeps the two mixes, and less total effort picks the cheap
one. On the held-out check the cheap mix separated from the measured vector, better on
\SearchHeldoutCheapVsControlBetter{} instances and worse on
\SearchHeldoutCheapVsControlWorse{}, and re-measured on the same instances in the headline run
it is \ReplicationRatio{}, \ReplicationCI{}. The tie-break picks among finalists the data do not
separate; the held-out check shows that the cheaper pick is not worse.

\paragraph{The selected configuration.}
The cheap mix ships as the released solver's option defaults~\cite{mip_heuristics_software}.
Against the measured vector it cuts total effort from \SearchEffortTotalMeasured{} to
\SearchArmTotalBPrime{} (Table~\ref{tab:arms}) by dropping Scylla and taking FPR and LocalMIP to
a fraction of their measured budgets.

\paragraph{Scylla ships disabled.}
Scylla produced no accepted incumbent in any full-limit run that carried it, in this stage or
with the measured vector in Stage~5: no incumbent in those logs carries its source tag. It runs
last in the chain, so a point of its own counts only if it beats the incumbent that FJ, FPR and
LocalMIP have left. It ships at effort zero, and its implementation remains.

\input{arms}
\input{fig/forest}

\subsection{Stage 4: the LP-guided dive}
\label{sec:fprlp-stage}

The LP-guided dive \texttt{fpr\_lp} runs after the root LP, so the presolve-only runs of
Stages~2 and~3 cannot see it. Stage~4 measures it in two steps. Running alone for
\SearchCapabilityTimeLimitS\,s on the confirmation instances, with HiGHS's own RENS and RINS
switched off so that it has their shared dive budget to itself, it produces accepted incumbents
on \SearchCapabilityDiveN{} of \SearchConfirmInstancesN{}, so the dive can find solutions. The second step adds
it, with RENS and RINS back at their defaults, to the selected configuration and to the measured
vector, and changes nothing else. Paired against their backgrounds in \ConstTimeLimitS\,s solves
on the confirmation instances, neither produced a single accepted incumbent from the dive. On the
selected configuration it changes nothing; on the measured vector it is a cost at the edge of
significance ($p=\SearchConfirmFprLpOnMeasuredP$; Figure~\ref{fig:forest}). It ships at effort
zero, and its implementation remains.

\subsection{Stage 5: the headline comparison}
\label{sec:headline}

Stage~5 runs the selected configuration and the measured vector over all \ConstNInstances{}
instances. Table~\ref{tab:headline} compares the selected configuration with the Stage-1
vanilla run on the full set, the tuning set and the held-out complement;
Table~\ref{tab:per-instance} in Appendix~\ref{app:per-instance} gives both on every instance,
and the Stage-5 block of Figure~\ref{fig:forest} shows the paired ratios. Over the full set the
primal-integral SGM falls from \HeadlineSgmHeadlineVanilla{} to \HeadlineSgmHeadlinePatched{}
(PI column), a paired decrease of \HeadlinePct\%. The tuning set is drawn from instances some
heuristic solves alone, so the held-out complement, at \HeldoutRatio{}, \HeldoutCI{}, is the
fairer estimate. There the selected configuration is better on \HeldoutBetter{} instances and
worse on \HeldoutWorse{}: it does not win more often, and the decrease comes from a minority of
instances with large gains. The \SearchHeldoutInstancesN{} held-out-check instances of the
complement were already seen in Stage~3; on the other \HeldoutUnseenN{} the ratio is
\HeldoutUnseenRatio{}, \HeldoutUnseenCI{} ($p=\HeldoutUnseenP$; Figure~\ref{fig:forest}). Over
the full set the selected configuration is better than the measured vector on
\SecondArmBetter{} instances and worse on \SecondArmWorse{}, a clear majority, but a few large
losses offset the wins and the ratio does not separate (\SecondArmRatio{}, \SecondArmCI{},
$p=\SecondArmP$; Figure~\ref{fig:forest}). The headline therefore does not depend on which of
the two was picked.

\paragraph{Over time and at the limit.}
Figure~\ref{fig:gap-time} shows the mean primal gap of both arms over time. Over the full set
the selected configuration's curve lies below vanilla's from the first tenth of a second to the
limit; on the held-out complement the two curves nearly meet after about twenty seconds. At the
limit (Table~\ref{tab:headline}) the selected configuration is feasible on
\HeadlineFeasiblePatched{} instances against \HeadlineFeasibleVanilla{}, holds the strictly
better objective on \HeadlineWinPatched{} against \HeadlineWinVanilla{}, and has the smaller
final-gap SGM, \HeadlineSgmGapPatched{} against \HeadlineSgmGapVanilla{}. It finds a solution on
\HeadlineFeasibleOnlyPatched{} instances where vanilla finds none (the rows of
Table~\ref{tab:per-instance} with a gap for the selected configuration only), and the reverse
happens on \HeadlineFeasibleOnlyVanilla{}, where vanilla's solution arrives as the solve closes;
on all \AttrOnlyPatchedOursN{} the first
solution comes from one of our heuristics. Paired over the
\HeadlinePairedGapN{} instances both made feasible, its final gap is smaller on
\HeadlineGapBetter{} and larger on \HeadlineGapWorse{}, a ratio of \HeadlinePairedGapRatio{},
\HeadlinePairedGapCI{} ($p=\HeadlineGapP$); on the complement it is \HeldoutPairedGapRatio{},
\HeldoutPairedGapCI{} over \HeldoutPairedGapN{} instances, which does not separate. The
heuristics buy a good solution \emph{sooner}; at the limit the patched solver is slightly
ahead, and the difference separates only over the full set.

One property of the host shows in these numbers. HiGHS passes the solution of a sub-MIP
heuristic to the main solve only when the sub-MIP returns. On five runs, two of the selected
configuration and three of vanilla, the final solution comes from a sub-MIP that the time limit
cut off (Section~\ref{sec:setup}), so it enters the primal integral at the limit and earns
nothing there, whenever the sub-MIP found it. On \texttt{neos-4647030-tutaki} the selected
configuration ends within a fraction of a per cent of the reference and still scores a primal
integral close to one: it holds a poor incumbent from the presolve chain until a sub-MIP
delivers the good one at the limit. The effect is the same in both arms and does not bias the
comparison.

\paragraph{Attribution.}
Table~\ref{tab:attribution} counts, over the \AttrFeasibleN{} instances the selected
configuration made feasible, the source of the first solution and of the best at
\ConstTimeLimitS\,s. Our heuristics find the first solution on most of them, and HiGHS's own
components hold the best solution at the limit on \AttrHighsOtherBest{}. Inside the chain, FPR
and LocalMIP add to what FJ finds. On \AttrChainEitherImprovesN{} of the
\AttrChainFjFirstN{} instances where FJ found the first solution, a later heuristic found a
better one, LocalMIP on \AttrChainLocalMipImprovesN{} and FPR on \AttrChainFprImprovesN{}; on
those, the median gap went from \AttrChainGapFjMed{} at FJ's best to \AttrChainGapAfterMed{}.
LocalMIP restarts from the incumbent, so part of its gain is refining FJ's solution.

\paragraph{The cost of budgeting.}
On six instances (\texttt{app1-2}, \texttt{fhnw-binpack4-48}, \texttt{germanrr},
\texttt{neos-1354092}, \texttt{neos-3656078-kumeu} and \texttt{neos-5104907-jarama}) a heuristic
running alone in the probe produced a solution within its $30$\,s cap, and neither the selected
configuration nor vanilla found one in \ConstTimeLimitS\,s. Neither did the measured vector,
with several times the budget and Scylla enabled. The heuristics can reach these instances and
a budgeted chain does not, which is the visible cost of budgeting. A presolve window
delays the root LP, so it has to be budgeted; a heuristic running alongside branch-and-bound
would not hold up the root and need not be cut short (Section~\ref{sec:conclusion}).

\input{headline}
\input{attribution}
\input{fig/gap-time}

\FloatBarrier

\section{Conclusion}
\label{sec:conclusion}

We have implemented four published primal heuristics (Feasibility Jump, FPR, LocalMIP and
Scylla) inside HiGHS behind one integration interface, and evaluated them on \texttt{mipfeas}
against a separately built unpatched HiGHS on the same machine. Section~\ref{sec:intro} asked
whether the added heuristics improve the full solve. They do, modestly. The
patched solver reaches a good solution sooner: the primal integral decreases by
\HeadlinePct\% over all \HeadlineN{} instances (\HeadlineRatio{}, \HeadlineCI{}) and by
\HeldoutPct\% on the \HeldoutN{} held-out instances (\HeldoutRatio{}, \HeldoutCI{},
$p=\HeldoutP$), the fairer estimate. At the limit it is slightly ahead. It finds a
solution on \HeadlineFeasibleOnlyPatched{} instances where vanilla HiGHS finds none, each first
found by one of our heuristics, against \HeadlineFeasibleOnlyVanilla{} the other way. Its final
gap is smaller on \HeadlineGapBetter{} instances and larger on \HeadlineGapWorse{}, which
separates over the full set (\HeadlinePairedGapRatio{}, \HeadlinePairedGapCI{}) but not on the
held-out instances (\HeldoutPairedGapRatio{}, \HeldoutPairedGapCI{}). HiGHS's own components hold the best
solution at the limit on \AttrHighsOtherBest{} of the \AttrFeasibleN{} instances the patched
solver makes feasible: the heuristics supply an earlier incumbent, and branch-and-bound improves
it from there.

Two of the five mechanisms contribute nothing in full solves. Scylla and the LP-guided dive \texttt{fpr\_lp} produced no
accepted incumbent in any full-limit run that carried them, Scylla although it produces on
\ProbeScyllaProduces{} instances when it runs alone, and both ship at effort zero; their
implementations remain, documented and runnable. The configuration that ships enables the other
three. It was picked by a tie-break among finalists the confirmation runs did not separate, and
over the full set the measured vector, with over four times the total effort and Scylla enabled,
does not separate from it either (\SecondArmRatio{}, \SecondArmCI{}). The comparison with vanilla
HiGHS is a net effect: what the heuristics displace inside the host was not measured.

We see two directions for further work. A presolve window delays the root LP and so has to be budgeted; the six
instances of Section~\ref{sec:headline} that a heuristic reaches alone and the budgeted chain
does not are what that costs. Running the heuristics alongside branch-and-bound removes that
constraint; ReXi~\cite{mexi2026rexi}, a parallel portfolio whose workers exchange information,
is a recent example of that design, and the heuristics here could serve as workers in such a portfolio. The second direction is selection: the productive-versus-stale effort
measured per heuristic is a natural input to per-instance heuristic selection. An online bandit we tried
before this campaign did not beat the fixed chain, but a single failed selector does not rule
selection out.

%% file: stages.tex

\begin{table}[t]
\centering
\small
\setlength{\tabcolsep}{4pt}
\resizebox{\ifdim\width>\linewidth\linewidth\else\width\fi}{!}{%
\begin{tabular}{llrlrr}
\toprule
Stage & What runs & Instances & Limit (s) & Configs & Runs \\
\midrule
1\enspace Baseline and the benchmark's difficulty profile & unpatched HiGHS & 233 & 600 & 1 & 233 \\
2\enspace Presolve probe & each heuristic alone & 233 & 30 (presolve) & 4 & 932 \\
3\enspace Joint search &  &  &  &  &  \\
\qquad Race & \texttt{irace} over $e_h, k_h$ & 90 & 60 (presolve) & 1447 & 12\,000 \\
\qquad Confirmation & race finalists & 49 & 600 & 5 & 245 \\
\qquad Held-out check & last finalists & 48 & 600 & 3 & 144 \\
4\enspace The LP-guided dive &  &  &  &  &  \\
\qquad Capability & \texttt{fpr\_lp} alone, RENS/RINS off & 49 & 120 & 1 & 49 \\
\qquad Contribution & \texttt{fpr\_lp} added, two backgrounds & 49 & 600 & 2 & 98 \\
5\enspace The headline comparison & selected, measured vector & 233 & 600 & 2 & 466 \\
\bottomrule
\end{tabular}
}
\par\vspace{7pt}
\caption{The five stages, one per subsection, in the order the paper reports them; indented rows are the separate runs within a stage. \emph{Presolve} marks the presolve-only runs; the limit beside it is the wall clock they still carry. The race runs on the tuning set; the confirmation and the \texttt{fpr\_lp} stages on a stratified draw from it; the held-out check on a stratified draw from the complement.}
\label{tab:stages}
\end{table}

%% file: profile.tex

\begin{table}[t]
\centering
\begin{tabular}{lrrr}
\toprule
Stratum (vanilla T1st) & \#Instances & Share (\%) & \#Tuning \\
\midrule
$<1$\,s & 92 & 39.5 & 42 \\
$1$--$10$\,s & 62 & 26.6 & 27 \\
$10$--$100$\,s & 40 & 17.2 & 10 \\
$100$--$600$\,s & 17 & 7.3 & 7 \\
never feasible & 22 & 9.4 & 4 \\
\midrule
total & 233 & 100.0 & 90 \\
\bottomrule
\end{tabular}
\par\vspace{7pt}
\caption{Feasibility-difficulty profile of \texttt{mipfeas}: vanilla time-to-first-feasible, HiGHS \texttt{v1.15.1}, 233 instances, $600$\,s, one seed, 16 workers. The last column is the tuning set: seats per stratum by largest remainder, at least one per stratum, drawn from the instances on which at least one heuristic produced a solution in the presolve probe.}
\label{tab:profile}
\end{table}

%% file: effort.tex

\begin{table}[t]
\centering
\small
\setlength{\tabcolsep}{4pt}
\resizebox{\ifdim\width>\linewidth\linewidth\else\width\fi}{!}{%
\begin{tabular}{lrrrrrrrr}
\toprule
Heuristic & Feasible & First & Incumbents & Gap (med.) & Productive (med.) & Stale (med.) & Measured $e_h$ & Measured $k_h$ \\
\midrule
FJ & 146 & 107 & 8125 & 0.099 & 1.565e+07 & 7.117e+10 & 0.5665 & 0.1416 \\
FPR & 136 & 109 & 724 & 0.448 & 1.387e+06 & 6.757e+10 & 12.26 & 3.064 \\
LocalMIP & 160 & 116 & 23\,168 & 0.182 & 1.232e+08 & 4.596e+10 & 13.96 & 3.49 \\
Scylla & 129 & 40 & 604 & 0.304 & 2.845e+06 & 8.845e+09 & 3.068 & 0.767 \\
\bottomrule
\end{tabular}
}
\par\vspace{7pt}
\caption{The presolve probe: each heuristic alone, presolve-only, effort unbounded, no patience, 233 instances. \emph{Feasible}: instances on which it found a feasible solution. \emph{First}: instances on which its first solution came no later than any other heuristic's; incumbents are logged at $0.1$\,s resolution, and a tie, on 94 instances, counts for every heuristic in it. \emph{Incumbents}: the solutions it found, each an improvement of its run's best objective, summed over instances. \emph{Gap}: the median gap $p$ of its best solution, over the instances it solved. \emph{Productive} and \emph{stale}: the median effort it spent before its last improvement and after it. Measured $e_h$ and $k_h$: the effort and patience the probe derives, the defaults Stage~3 starts from. Effort is in each heuristic's own unit and is not comparable across rows.}
\label{tab:effort}
\end{table}

%% file: fig/probe-time.tex

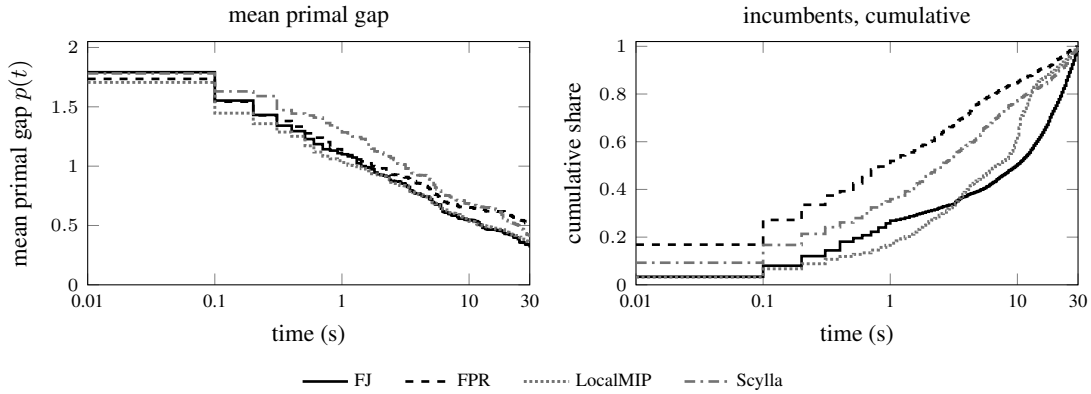
\begin{figure}[t]
\centering
\begin{tikzpicture}
\begin{groupplot}[
  group style={group size=2 by 1, horizontal sep=40pt},
  width=0.45\linewidth, height=4.8cm,
  xmode=log, xmin=0.01, xmax=30,
  xtick={0.01,0.1,1,10,30}, xticklabels={0.01,0.1,1,10,30},
  tick label style={font=\scriptsize}, label style={font=\small},
  title style={font=\small, yshift=-4pt}, xlabel={time (s)},
]
\nextgroupplot[title={mean primal gap}, ymin=0, ymax=2.05, ylabel={mean primal gap $p(t)$}]
\addplot[const plot, black, line width=1pt] coordinates {(0.01,1.79154) (0.0102714,1.79154) (0.0105501,1.79154) (0.0108365,1.79154) (0.0111306,1.79154) (0.0114326,1.79154) (0.0117429,1.79154) (0.0120616,1.79154) (0.0123889,1.79154) (0.0127251,1.79154) (0.0130705,1.79154) (0.0134252,1.79154) (0.0137895,1.79154) (0.0141638,1.79154) (0.0145482,1.79154) (0.014943,1.79154) (0.0153485,1.79154) (0.0157651,1.79154) (0.0161929,1.79154) (0.0166324,1.79154) (0.0170838,1.79154) (0.0175474,1.79154) (0.0180236,1.79154) (0.0185127,1.79154) (0.0190152,1.79154) (0.0195312,1.79154) (0.0200613,1.79154) (0.0206057,1.79154) (0.0211649,1.79154) (0.0217393,1.79154) (0.0223293,1.79154) (0.0229353,1.79154) (0.0235577,1.79154) (0.0241971,1.79154) (0.0248537,1.79154) (0.0255282,1.79154) (0.026221,1.79154) (0.0269327,1.79154) (0.0276636,1.79154) (0.0284143,1.79154) (0.0291855,1.79154) (0.0299775,1.79154) (0.0307911,1.79154) (0.0316267,1.79154) (0.032485,1.79154) (0.0333666,1.79154) (0.0342722,1.79154) (0.0352023,1.79154) (0.0361576,1.79154) (0.0371389,1.79154) (0.0381468,1.79154) (0.0391821,1.79154) (0.0402454,1.79154) (0.0413377,1.79154) (0.0424595,1.79154) (0.0436118,1.79154) (0.0447954,1.79154) (0.0460111,1.79154) (0.0472598,1.79154) (0.0485424,1.79154) (0.0498597,1.79154) (0.0512129,1.79154) (0.0526027,1.79154) (0.0540303,1.79154) (0.0554966,1.79154) (0.0570028,1.79154) (0.0585498,1.79154) (0.0601387,1.79154) (0.0617708,1.79154) (0.0634472,1.79154) (0.0651691,1.79154) (0.0669377,1.79154) (0.0687543,1.79154) (0.0706203,1.79154) (0.0725368,1.79154) (0.0745054,1.79154) (0.0765274,1.79154) (0.0786042,1.79154) (0.0807375,1.79154) (0.0829286,1.79154) (0.0851792,1.79154) (0.0874908,1.79154) (0.0898652,1.79154) (0.0923041,1.79154) (0.0948091,1.79154) (0.0973821,1.79154) (0.100025,1.55268) (0.10274,1.55268) (0.105528,1.55268) (0.108392,1.55268) (0.111333,1.55268) (0.114355,1.55268) (0.117458,1.55268) (0.120646,1.55268) (0.12392,1.55268) (0.127283,1.55268) (0.130737,1.55268) (0.134286,1.55268) (0.13793,1.55268) (0.141673,1.55268) (0.145518,1.55268) (0.149467,1.55268) (0.153524,1.55268) (0.15769,1.55268) (0.16197,1.55268) (0.166365,1.55268) (0.17088,1.55268) (0.175518,1.55268) (0.180281,1.55268) (0.185174,1.55268) (0.190199,1.55268) (0.195361,1.55268) (0.200663,1.43164) (0.206109,1.43164) (0.211702,1.43164) (0.217447,1.43164) (0.223349,1.43164) (0.22941,1.43164) (0.235636,1.43164) (0.242031,1.43164) (0.248599,1.43164) (0.255346,1.43164) (0.262276,1.43164) (0.269394,1.43164) (0.276705,1.43164) (0.284214,1.43164) (0.291928,1.43164) (0.29985,1.43164) (0.307988,1.34139) (0.316346,1.34139) (0.324932,1.34139) (0.33375,1.34139) (0.342807,1.34139) (0.352111,1.34139) (0.361667,1.34139) (0.371482,1.34139) (0.381564,1.34139) (0.391919,1.34139) (0.402555,1.29504) (0.41348,1.29504) (0.424701,1.29504) (0.436227,1.29504) (0.448066,1.29504) (0.460226,1.29504) (0.472716,1.29504) (0.485545,1.29504) (0.498722,1.29504) (0.512257,1.22819) (0.526159,1.22819) (0.540438,1.22819) (0.555105,1.22819) (0.57017,1.22819) (0.585644,1.22819) (0.601538,1.18551) (0.617863,1.18551) (0.634631,1.18551) (0.651854,1.18551) (0.669544,1.18551) (0.687715,1.18551) (0.706379,1.14248) (0.725549,1.14248) (0.74524,1.14248) (0.765465,1.14248) (0.786239,1.14248) (0.807576,1.13111) (0.829493,1.13111) (0.852005,1.13111) (0.875127,1.13111) (0.898877,1.13111) (0.923271,1.10951) (0.948328,1.10951) (0.974065,1.10951) (1.0005,1.10011) (1.02765,1.10011) (1.05554,1.10011) (1.08419,1.10011) (1.11361,1.0825) (1.14383,1.0825) (1.17488,1.0825) (1.20676,1.07286) (1.23951,1.07286) (1.27315,1.07286) (1.3077,1.0244) (1.34319,1.0244) (1.37964,1.0244) (1.41709,1.00023) (1.45554,1.00023) (1.49505,1.00023) (1.53562,0.991125) (1.57729,0.991125) (1.6201,0.988396) (1.66407,0.988396) (1.70923,0.973667) (1.75562,0.973667) (1.80326,0.946737) (1.8522,0.946737) (1.90247,0.931686) (1.9541,0.931686) (2.00713,0.926889) (2.0616,0.926889) (2.11755,0.919902) (2.17502,0.919902) (2.23405,0.912787) (2.29467,0.912787) (2.35695,0.884808) (2.42091,0.875392) (2.48662,0.875392) (2.5541,0.874505) (2.62341,0.870775) (2.69461,0.870775) (2.76774,0.869771) (2.84285,0.868701) (2.92001,0.854276) (2.99925,0.854276) (3.08065,0.85282) (3.16425,0.839441) (3.25013,0.825666) (3.33833,0.824907) (3.42893,0.795492) (3.52199,0.793728) (3.61757,0.782844) (3.71575,0.772626) (3.81659,0.770344) (3.92017,0.769827) (4.02656,0.768857) (4.13583,0.761124) (4.24807,0.749697) (4.36336,0.745719) (4.48178,0.745192) (4.60341,0.7451) (4.72834,0.744092) (4.85666,0.722618) (4.98847,0.721388) (5.12385,0.697328) (5.2629,0.690366) (5.40573,0.676983) (5.55244,0.676913) (5.70312,0.676009) (5.8579,0.675964) (6.01688,0.641606) (6.18017,0.618125) (6.34789,0.611504) (6.52017,0.611294) (6.69712,0.6077) (6.87887,0.59762) (7.06555,0.597107) (7.2573,0.586523) (7.45426,0.585497) (7.65656,0.585396) (7.86435,0.584811) (8.07778,0.571481) (8.297,0.557656) (8.52217,0.557327) (8.75346,0.557247) (8.99101,0.556922) (9.23502,0.548531) (9.48565,0.548007) (9.74308,0.546316) (10.0075,0.546081) (10.2791,0.546028) (10.558,0.545471) (10.8446,0.543866) (11.1389,0.535825) (11.4412,0.535641) (11.7517,0.526571) (12.0706,0.505972) (12.3982,0.504421) (12.7347,0.481047) (13.0803,0.479711) (13.4353,0.471664) (13.7999,0.467108) (14.1744,0.466791) (14.5591,0.466728) (14.9542,0.466664) (15.36,0.465743) (15.7769,0.465318) (16.2051,0.462485) (16.6448,0.457431) (17.0966,0.455957) (17.5605,0.454953) (18.0371,0.454094) (18.5266,0.441226) (19.0294,0.429768) (19.5459,0.429756) (20.0763,0.429714) (20.6212,0.422213) (21.1808,0.420602) (21.7556,0.420277) (22.346,0.400295) (22.9525,0.399704) (23.5754,0.399674) (24.2152,0.385463) (24.8724,0.359685) (25.5474,0.355953) (26.2407,0.355488) (26.9528,0.35442) (27.6843,0.340506) (28.4356,0.338207) (29.2073,0.337714) (30,0.314353)};
\addplot[const plot, black, line width=1pt, dashed] coordinates {(0.01,1.73527) (0.0102714,1.73527) (0.0105501,1.73527) (0.0108365,1.73527) (0.0111306,1.73527) (0.0114326,1.73527) (0.0117429,1.73527) (0.0120616,1.73527) (0.0123889,1.73527) (0.0127251,1.73527) (0.0130705,1.73527) (0.0134252,1.73527) (0.0137895,1.73527) (0.0141638,1.73527) (0.0145482,1.73527) (0.014943,1.73527) (0.0153485,1.73527) (0.0157651,1.73527) (0.0161929,1.73527) (0.0166324,1.73527) (0.0170838,1.73527) (0.0175474,1.73527) (0.0180236,1.73527) (0.0185127,1.73527) (0.0190152,1.73527) (0.0195312,1.73527) (0.0200613,1.73527) (0.0206057,1.73527) (0.0211649,1.73527) (0.0217393,1.73527) (0.0223293,1.73527) (0.0229353,1.73527) (0.0235577,1.73527) (0.0241971,1.73527) (0.0248537,1.73527) (0.0255282,1.73527) (0.026221,1.73527) (0.0269327,1.73527) (0.0276636,1.73527) (0.0284143,1.73527) (0.0291855,1.73527) (0.0299775,1.73527) (0.0307911,1.73527) (0.0316267,1.73527) (0.032485,1.73527) (0.0333666,1.73527) (0.0342722,1.73527) (0.0352023,1.73527) (0.0361576,1.73527) (0.0371389,1.73527) (0.0381468,1.73527) (0.0391821,1.73527) (0.0402454,1.73527) (0.0413377,1.73527) (0.0424595,1.73527) (0.0436118,1.73527) (0.0447954,1.73527) (0.0460111,1.73527) (0.0472598,1.73527) (0.0485424,1.73527) (0.0498597,1.73527) (0.0512129,1.73527) (0.0526027,1.73527) (0.0540303,1.73527) (0.0554966,1.73527) (0.0570028,1.73527) (0.0585498,1.73527) (0.0601387,1.73527) (0.0617708,1.73527) (0.0634472,1.73527) (0.0651691,1.73527) (0.0669377,1.73527) (0.0687543,1.73527) (0.0706203,1.73527) (0.0725368,1.73527) (0.0745054,1.73527) (0.0765274,1.73527) (0.0786042,1.73527) (0.0807375,1.73527) (0.0829286,1.73527) (0.0851792,1.73527) (0.0874908,1.73527) (0.0898652,1.73527) (0.0923041,1.73527) (0.0948091,1.73527) (0.0973821,1.73527) (0.100025,1.54327) (0.10274,1.54327) (0.105528,1.54327) (0.108392,1.54327) (0.111333,1.54327) (0.114355,1.54327) (0.117458,1.54327) (0.120646,1.54327) (0.12392,1.54327) (0.127283,1.54327) (0.130737,1.54327) (0.134286,1.54327) (0.13793,1.54327) (0.141673,1.54327) (0.145518,1.54327) (0.149467,1.54327) (0.153524,1.54327) (0.15769,1.54327) (0.16197,1.54327) (0.166365,1.54327) (0.17088,1.54327) (0.175518,1.54327) (0.180281,1.54327) (0.185174,1.54327) (0.190199,1.54327) (0.195361,1.54327) (0.200663,1.42762) (0.206109,1.42762) (0.211702,1.42762) (0.217447,1.42762) (0.223349,1.42762) (0.22941,1.42762) (0.235636,1.42762) (0.242031,1.42762) (0.248599,1.42762) (0.255346,1.42762) (0.262276,1.42762) (0.269394,1.42762) (0.276705,1.42762) (0.284214,1.42762) (0.291928,1.42762) (0.29985,1.42762) (0.307988,1.38025) (0.316346,1.38025) (0.324932,1.38025) (0.33375,1.38025) (0.342807,1.38025) (0.352111,1.38025) (0.361667,1.38025) (0.371482,1.38025) (0.381564,1.38025) (0.391919,1.38025) (0.402555,1.33335) (0.41348,1.33335) (0.424701,1.33335) (0.436227,1.33335) (0.448066,1.33335) (0.460226,1.33335) (0.472716,1.33335) (0.485545,1.33335) (0.498722,1.33335) (0.512257,1.2767) (0.526159,1.2767) (0.540438,1.2767) (0.555105,1.2767) (0.57017,1.2767) (0.585644,1.2767) (0.601538,1.23989) (0.617863,1.23989) (0.634631,1.23989) (0.651854,1.23989) (0.669544,1.23989) (0.687715,1.23989) (0.706379,1.20279) (0.725549,1.20279) (0.74524,1.20279) (0.765465,1.20279) (0.786239,1.20279) (0.807576,1.16246) (0.829493,1.16246) (0.852005,1.16246) (0.875127,1.16246) (0.898877,1.16246) (0.923271,1.14273) (0.948328,1.14273) (0.974065,1.14273) (1.0005,1.12065) (1.02765,1.12065) (1.05554,1.12065) (1.08419,1.12065) (1.11361,1.09297) (1.14383,1.09297) (1.17488,1.09297) (1.20676,1.06903) (1.23951,1.06903) (1.27315,1.06903) (1.3077,1.05686) (1.34319,1.05686) (1.37964,1.05686) (1.41709,1.04077) (1.45554,1.04077) (1.49505,1.04077) (1.53562,1.02426) (1.57729,1.02426) (1.6201,1.02411) (1.66407,1.02411) (1.70923,0.999962) (1.75562,0.999962) (1.80326,0.999586) (1.8522,0.999586) (1.90247,0.990408) (1.9541,0.990408) (2.00713,0.981943) (2.0616,0.981943) (2.11755,0.981808) (2.17502,0.981808) (2.23405,0.943515) (2.29467,0.943515) (2.35695,0.9355) (2.42091,0.928082) (2.48662,0.928082) (2.5541,0.927041) (2.62341,0.92685) (2.69461,0.92685) (2.76774,0.92491) (2.84285,0.923927) (2.92001,0.916524) (2.99925,0.916524) (3.08065,0.907242) (3.16425,0.907235) (3.25013,0.907226) (3.33833,0.895402) (3.42893,0.875395) (3.52199,0.874572) (3.61757,0.874345) (3.71575,0.873545) (3.81659,0.857054) (3.92017,0.857053) (4.02656,0.856158) (4.13583,0.85292) (4.24807,0.843673) (4.36336,0.840771) (4.48178,0.831915) (4.60341,0.831791) (4.72834,0.819558) (4.85666,0.807277) (4.98847,0.794833) (5.12385,0.794423) (5.2629,0.767276) (5.40573,0.759923) (5.55244,0.758742) (5.70312,0.743888) (5.8579,0.714735) (6.01688,0.7137) (6.18017,0.700415) (6.34789,0.698632) (6.52017,0.698632) (6.69712,0.688577) (6.87887,0.687732) (7.06555,0.687497) (7.2573,0.687488) (7.45426,0.687028) (7.65656,0.687028) (7.86435,0.679065) (8.07778,0.665677) (8.297,0.656263) (8.52217,0.656162) (8.75346,0.656162) (8.99101,0.655455) (9.23502,0.655286) (9.48565,0.655092) (9.74308,0.653157) (10.0075,0.653053) (10.2791,0.652812) (10.558,0.652252) (10.8446,0.643717) (11.1389,0.643181) (11.4412,0.643168) (11.7517,0.642859) (12.0706,0.633915) (12.3982,0.633915) (12.7347,0.633915) (13.0803,0.633896) (13.4353,0.621256) (13.7999,0.620265) (14.1744,0.620265) (14.5591,0.619469) (14.9542,0.619224) (15.36,0.619217) (15.7769,0.618903) (16.2051,0.618903) (16.6448,0.618889) (17.0966,0.618814) (17.5605,0.618618) (18.0371,0.61857) (18.5266,0.618108) (19.0294,0.600715) (19.5459,0.598928) (20.0763,0.598088) (20.6212,0.586511) (21.1808,0.586175) (21.7556,0.574642) (22.346,0.566383) (22.9525,0.556474) (23.5754,0.556462) (24.2152,0.556233) (24.8724,0.54446) (25.5474,0.54394) (26.2407,0.543924) (26.9528,0.536494) (27.6843,0.520695) (28.4356,0.520695) (29.2073,0.511888) (30,0.51184)};
\addplot[const plot, black!55, line width=1.1pt, densely dotted] coordinates {(0.01,1.70568) (0.0102714,1.70568) (0.0105501,1.70568) (0.0108365,1.70568) (0.0111306,1.70568) (0.0114326,1.70568) (0.0117429,1.70568) (0.0120616,1.70568) (0.0123889,1.70568) (0.0127251,1.70568) (0.0130705,1.70568) (0.0134252,1.70568) (0.0137895,1.70568) (0.0141638,1.70568) (0.0145482,1.70568) (0.014943,1.70568) (0.0153485,1.70568) (0.0157651,1.70568) (0.0161929,1.70568) (0.0166324,1.70568) (0.0170838,1.70568) (0.0175474,1.70568) (0.0180236,1.70568) (0.0185127,1.70568) (0.0190152,1.70568) (0.0195312,1.70568) (0.0200613,1.70568) (0.0206057,1.70568) (0.0211649,1.70568) (0.0217393,1.70568) (0.0223293,1.70568) (0.0229353,1.70568) (0.0235577,1.70568) (0.0241971,1.70568) (0.0248537,1.70568) (0.0255282,1.70568) (0.026221,1.70568) (0.0269327,1.70568) (0.0276636,1.70568) (0.0284143,1.70568) (0.0291855,1.70568) (0.0299775,1.70568) (0.0307911,1.70568) (0.0316267,1.70568) (0.032485,1.70568) (0.0333666,1.70568) (0.0342722,1.70568) (0.0352023,1.70568) (0.0361576,1.70568) (0.0371389,1.70568) (0.0381468,1.70568) (0.0391821,1.70568) (0.0402454,1.70568) (0.0413377,1.70568) (0.0424595,1.70568) (0.0436118,1.70568) (0.0447954,1.70568) (0.0460111,1.70568) (0.0472598,1.70568) (0.0485424,1.70568) (0.0498597,1.70568) (0.0512129,1.70568) (0.0526027,1.70568) (0.0540303,1.70568) (0.0554966,1.70568) (0.0570028,1.70568) (0.0585498,1.70568) (0.0601387,1.70568) (0.0617708,1.70568) (0.0634472,1.70568) (0.0651691,1.70568) (0.0669377,1.70568) (0.0687543,1.70568) (0.0706203,1.70568) (0.0725368,1.70568) (0.0745054,1.70568) (0.0765274,1.70568) (0.0786042,1.70568) (0.0807375,1.70568) (0.0829286,1.70568) (0.0851792,1.70568) (0.0874908,1.70568) (0.0898652,1.70568) (0.0923041,1.70568) (0.0948091,1.70568) (0.0973821,1.70568) (0.100025,1.4466) (0.10274,1.4466) (0.105528,1.4466) (0.108392,1.4466) (0.111333,1.4466) (0.114355,1.4466) (0.117458,1.4466) (0.120646,1.4466) (0.12392,1.4466) (0.127283,1.4466) (0.130737,1.4466) (0.134286,1.4466) (0.13793,1.4466) (0.141673,1.4466) (0.145518,1.4466) (0.149467,1.4466) (0.153524,1.4466) (0.15769,1.4466) (0.16197,1.4466) (0.166365,1.4466) (0.17088,1.4466) (0.175518,1.4466) (0.180281,1.4466) (0.185174,1.4466) (0.190199,1.4466) (0.195361,1.4466) (0.200663,1.35797) (0.206109,1.35797) (0.211702,1.35797) (0.217447,1.35797) (0.223349,1.35797) (0.22941,1.35797) (0.235636,1.35797) (0.242031,1.35797) (0.248599,1.35797) (0.255346,1.35797) (0.262276,1.35797) (0.269394,1.35797) (0.276705,1.35797) (0.284214,1.35797) (0.291928,1.35797) (0.29985,1.35797) (0.307988,1.28727) (0.316346,1.28727) (0.324932,1.28727) (0.33375,1.28727) (0.342807,1.28727) (0.352111,1.28727) (0.361667,1.28727) (0.371482,1.28727) (0.381564,1.28727) (0.391919,1.28727) (0.402555,1.25157) (0.41348,1.25157) (0.424701,1.25157) (0.436227,1.25157) (0.448066,1.25157) (0.460226,1.25157) (0.472716,1.25157) (0.485545,1.25157) (0.498722,1.25157) (0.512257,1.17158) (0.526159,1.17158) (0.540438,1.17158) (0.555105,1.17158) (0.57017,1.17158) (0.585644,1.17158) (0.601538,1.11894) (0.617863,1.11894) (0.634631,1.11894) (0.651854,1.11894) (0.669544,1.11894) (0.687715,1.11894) (0.706379,1.08389) (0.725549,1.08389) (0.74524,1.08389) (0.765465,1.08389) (0.786239,1.08389) (0.807576,1.06155) (0.829493,1.06155) (0.852005,1.06155) (0.875127,1.06155) (0.898877,1.06155) (0.923271,1.04619) (0.948328,1.04619) (0.974065,1.04619) (1.0005,1.03069) (1.02765,1.03069) (1.05554,1.03069) (1.08419,1.03069) (1.11361,1.00713) (1.14383,1.00713) (1.17488,1.00713) (1.20676,1.00543) (1.23951,1.00543) (1.27315,1.00543) (1.3077,0.995059) (1.34319,0.995059) (1.37964,0.995059) (1.41709,0.978703) (1.45554,0.978703) (1.49505,0.978703) (1.53562,0.975916) (1.57729,0.975916) (1.6201,0.960386) (1.66407,0.960386) (1.70923,0.936037) (1.75562,0.936037) (1.80326,0.934337) (1.8522,0.934337) (1.90247,0.926454) (1.9541,0.926454) (2.00713,0.924245) (2.0616,0.924245) (2.11755,0.912256) (2.17502,0.912256) (2.23405,0.898746) (2.29467,0.898746) (2.35695,0.88778) (2.42091,0.863087) (2.48662,0.863087) (2.5541,0.855977) (2.62341,0.853419) (2.69461,0.853419) (2.76774,0.852589) (2.84285,0.836917) (2.92001,0.835923) (2.99925,0.835923) (3.08065,0.834842) (3.16425,0.834521) (3.25013,0.834019) (3.33833,0.802876) (3.42893,0.785948) (3.52199,0.778867) (3.61757,0.77717) (3.71575,0.775978) (3.81659,0.768959) (3.92017,0.768363) (4.02656,0.760205) (4.13583,0.758973) (4.24807,0.756273) (4.36336,0.749107) (4.48178,0.746355) (4.60341,0.741525) (4.72834,0.740776) (4.85666,0.733574) (4.98847,0.708636) (5.12385,0.691914) (5.2629,0.684478) (5.40573,0.683891) (5.55244,0.682756) (5.70312,0.671539) (5.8579,0.64654) (6.01688,0.644472) (6.18017,0.637617) (6.34789,0.632391) (6.52017,0.6313) (6.69712,0.621951) (6.87887,0.621445) (7.06555,0.60697) (7.2573,0.604188) (7.45426,0.603467) (7.65656,0.601949) (7.86435,0.592851) (8.07778,0.592055) (8.297,0.59122) (8.52217,0.57292) (8.75346,0.572533) (8.99101,0.571912) (9.23502,0.570081) (9.48565,0.556817) (9.74308,0.5443) (10.0075,0.534169) (10.2791,0.533865) (10.558,0.532936) (10.8446,0.53209) (11.1389,0.520002) (11.4412,0.512924) (11.7517,0.505052) (12.0706,0.502305) (12.3982,0.498864) (12.7347,0.497082) (13.0803,0.496172) (13.4353,0.495074) (13.7999,0.487339) (14.1744,0.486608) (14.5591,0.486302) (14.9542,0.485422) (15.36,0.477854) (15.7769,0.477694) (16.2051,0.470152) (16.6448,0.470099) (17.0966,0.468513) (17.5605,0.466881) (18.0371,0.464896) (18.5266,0.464304) (19.0294,0.457981) (19.5459,0.443525) (20.0763,0.442652) (20.6212,0.427017) (21.1808,0.426721) (21.7556,0.426349) (22.346,0.41903) (22.9525,0.415377) (23.5754,0.415035) (24.2152,0.408034) (24.8724,0.40159) (25.5474,0.398815) (26.2407,0.398266) (26.9528,0.382002) (27.6843,0.368301) (28.4356,0.36772) (29.2073,0.367619) (30,0.367209)};
\addplot[const plot, black!55, line width=1pt, dash pattern=on 4pt off 2pt on 1pt off 2pt] coordinates {(0.01,1.78135) (0.0102714,1.78135) (0.0105501,1.78135) (0.0108365,1.78135) (0.0111306,1.78135) (0.0114326,1.78135) (0.0117429,1.78135) (0.0120616,1.78135) (0.0123889,1.78135) (0.0127251,1.78135) (0.0130705,1.78135) (0.0134252,1.78135) (0.0137895,1.78135) (0.0141638,1.78135) (0.0145482,1.78135) (0.014943,1.78135) (0.0153485,1.78135) (0.0157651,1.78135) (0.0161929,1.78135) (0.0166324,1.78135) (0.0170838,1.78135) (0.0175474,1.78135) (0.0180236,1.78135) (0.0185127,1.78135) (0.0190152,1.78135) (0.0195312,1.78135) (0.0200613,1.78135) (0.0206057,1.78135) (0.0211649,1.78135) (0.0217393,1.78135) (0.0223293,1.78135) (0.0229353,1.78135) (0.0235577,1.78135) (0.0241971,1.78135) (0.0248537,1.78135) (0.0255282,1.78135) (0.026221,1.78135) (0.0269327,1.78135) (0.0276636,1.78135) (0.0284143,1.78135) (0.0291855,1.78135) (0.0299775,1.78135) (0.0307911,1.78135) (0.0316267,1.78135) (0.032485,1.78135) (0.0333666,1.78135) (0.0342722,1.78135) (0.0352023,1.78135) (0.0361576,1.78135) (0.0371389,1.78135) (0.0381468,1.78135) (0.0391821,1.78135) (0.0402454,1.78135) (0.0413377,1.78135) (0.0424595,1.78135) (0.0436118,1.78135) (0.0447954,1.78135) (0.0460111,1.78135) (0.0472598,1.78135) (0.0485424,1.78135) (0.0498597,1.78135) (0.0512129,1.78135) (0.0526027,1.78135) (0.0540303,1.78135) (0.0554966,1.78135) (0.0570028,1.78135) (0.0585498,1.78135) (0.0601387,1.78135) (0.0617708,1.78135) (0.0634472,1.78135) (0.0651691,1.78135) (0.0669377,1.78135) (0.0687543,1.78135) (0.0706203,1.78135) (0.0725368,1.78135) (0.0745054,1.78135) (0.0765274,1.78135) (0.0786042,1.78135) (0.0807375,1.78135) (0.0829286,1.78135) (0.0851792,1.78135) (0.0874908,1.78135) (0.0898652,1.78135) (0.0923041,1.78135) (0.0948091,1.78135) (0.0973821,1.78135) (0.100025,1.62985) (0.10274,1.62985) (0.105528,1.62985) (0.108392,1.62985) (0.111333,1.62985) (0.114355,1.62985) (0.117458,1.62985) (0.120646,1.62985) (0.12392,1.62985) (0.127283,1.62985) (0.130737,1.62985) (0.134286,1.62985) (0.13793,1.62985) (0.141673,1.62985) (0.145518,1.62985) (0.149467,1.62985) (0.153524,1.62985) (0.15769,1.62985) (0.16197,1.62985) (0.166365,1.62985) (0.17088,1.62985) (0.175518,1.62985) (0.180281,1.62985) (0.185174,1.62985) (0.190199,1.62985) (0.195361,1.62985) (0.200663,1.59029) (0.206109,1.59029) (0.211702,1.59029) (0.217447,1.59029) (0.223349,1.59029) (0.22941,1.59029) (0.235636,1.59029) (0.242031,1.59029) (0.248599,1.59029) (0.255346,1.59029) (0.262276,1.59029) (0.269394,1.59029) (0.276705,1.59029) (0.284214,1.59029) (0.291928,1.59029) (0.29985,1.59029) (0.307988,1.47219) (0.316346,1.47219) (0.324932,1.47219) (0.33375,1.47219) (0.342807,1.47219) (0.352111,1.47219) (0.361667,1.47219) (0.371482,1.47219) (0.381564,1.47219) (0.391919,1.47219) (0.402555,1.44301) (0.41348,1.44301) (0.424701,1.44301) (0.436227,1.44301) (0.448066,1.44301) (0.460226,1.44301) (0.472716,1.44301) (0.485545,1.44301) (0.498722,1.44301) (0.512257,1.427) (0.526159,1.427) (0.540438,1.427) (0.555105,1.427) (0.57017,1.427) (0.585644,1.427) (0.601538,1.40912) (0.617863,1.40912) (0.634631,1.40912) (0.651854,1.40912) (0.669544,1.40912) (0.687715,1.40912) (0.706379,1.36314) (0.725549,1.36314) (0.74524,1.36314) (0.765465,1.36314) (0.786239,1.36314) (0.807576,1.34397) (0.829493,1.34397) (0.852005,1.34397) (0.875127,1.34397) (0.898877,1.34397) (0.923271,1.30505) (0.948328,1.30505) (0.974065,1.30505) (1.0005,1.28467) (1.02765,1.28467) (1.05554,1.28467) (1.08419,1.28467) (1.11361,1.25853) (1.14383,1.25853) (1.17488,1.25853) (1.20676,1.24968) (1.23951,1.24968) (1.27315,1.24968) (1.3077,1.24442) (1.34319,1.24442) (1.37964,1.24442) (1.41709,1.22907) (1.45554,1.22907) (1.49505,1.22907) (1.53562,1.19632) (1.57729,1.19632) (1.6201,1.17043) (1.66407,1.17043) (1.70923,1.14722) (1.75562,1.14722) (1.80326,1.12474) (1.8522,1.12474) (1.90247,1.12329) (1.9541,1.12329) (2.00713,1.11458) (2.0616,1.11458) (2.11755,1.10024) (2.17502,1.10024) (2.23405,1.0839) (2.29467,1.0839) (2.35695,1.08289) (2.42091,1.06823) (2.48662,1.06823) (2.5541,1.05314) (2.62341,1.022) (2.69461,1.022) (2.76774,1.02112) (2.84285,1.01915) (2.92001,1.01369) (2.99925,1.01369) (3.08065,1.01329) (3.16425,0.979317) (3.25013,0.979175) (3.33833,0.978735) (3.42893,0.977282) (3.52199,0.942408) (3.61757,0.941715) (3.71575,0.940364) (3.81659,0.938807) (3.92017,0.92462) (4.02656,0.924565) (4.13583,0.923691) (4.24807,0.92353) (4.36336,0.92353) (4.48178,0.923359) (4.60341,0.921775) (4.72834,0.891182) (4.85666,0.890209) (4.98847,0.890205) (5.12385,0.874523) (5.2629,0.842429) (5.40573,0.818309) (5.55244,0.804378) (5.70312,0.794023) (5.8579,0.793567) (6.01688,0.767403) (6.18017,0.757734) (6.34789,0.755983) (6.52017,0.752675) (6.69712,0.751846) (6.87887,0.751378) (7.06555,0.741026) (7.2573,0.725255) (7.45426,0.725021) (7.65656,0.72405) (7.86435,0.714045) (8.07778,0.714036) (8.297,0.713834) (8.52217,0.710943) (8.75346,0.710409) (8.99101,0.710091) (9.23502,0.687004) (9.48565,0.686791) (9.74308,0.686777) (10.0075,0.681669) (10.2791,0.673235) (10.558,0.672815) (10.8446,0.671209) (11.1389,0.670623) (11.4412,0.670518) (11.7517,0.656651) (12.0706,0.656451) (12.3982,0.65634) (12.7347,0.656318) (13.0803,0.653323) (13.4353,0.652137) (13.7999,0.652111) (14.1744,0.651942) (14.5591,0.650029) (14.9542,0.649565) (15.36,0.640968) (15.7769,0.640968) (16.2051,0.631195) (16.6448,0.630841) (17.0966,0.614019) (17.5605,0.613381) (18.0371,0.591015) (18.5266,0.561894) (19.0294,0.56125) (19.5459,0.545788) (20.0763,0.536734) (20.6212,0.525413) (21.1808,0.521129) (21.7556,0.520592) (22.346,0.519966) (22.9525,0.510979) (23.5754,0.509218) (24.2152,0.499791) (24.8724,0.485821) (25.5474,0.470277) (26.2407,0.468816) (26.9528,0.449379) (27.6843,0.433674) (28.4356,0.421812) (29.2073,0.421533) (30,0.401425)};
\nextgroupplot[title={incumbents, cumulative}, ymin=0, ymax=1.02, ylabel={cumulative share}, legend to name=probetimelegend, legend columns=4, legend style={font=\scriptsize, draw=none, /tikz/every even column/.append style={column sep=10pt}}]
\addplot[const plot, black, line width=1pt] coordinates {(0.01,0.0338462) (0.0102714,0.0338462) (0.0105501,0.0338462) (0.0108365,0.0338462) (0.0111306,0.0338462) (0.0114326,0.0338462) (0.0117429,0.0338462) (0.0120616,0.0338462) (0.0123889,0.0338462) (0.0127251,0.0338462) (0.0130705,0.0338462) (0.0134252,0.0338462) (0.0137895,0.0338462) (0.0141638,0.0338462) (0.0145482,0.0338462) (0.014943,0.0338462) (0.0153485,0.0338462) (0.0157651,0.0338462) (0.0161929,0.0338462) (0.0166324,0.0338462) (0.0170838,0.0338462) (0.0175474,0.0338462) (0.0180236,0.0338462) (0.0185127,0.0338462) (0.0190152,0.0338462) (0.0195312,0.0338462) (0.0200613,0.0338462) (0.0206057,0.0338462) (0.0211649,0.0338462) (0.0217393,0.0338462) (0.0223293,0.0338462) (0.0229353,0.0338462) (0.0235577,0.0338462) (0.0241971,0.0338462) (0.0248537,0.0338462) (0.0255282,0.0338462) (0.026221,0.0338462) (0.0269327,0.0338462) (0.0276636,0.0338462) (0.0284143,0.0338462) (0.0291855,0.0338462) (0.0299775,0.0338462) (0.0307911,0.0338462) (0.0316267,0.0338462) (0.032485,0.0338462) (0.0333666,0.0338462) (0.0342722,0.0338462) (0.0352023,0.0338462) (0.0361576,0.0338462) (0.0371389,0.0338462) (0.0381468,0.0338462) (0.0391821,0.0338462) (0.0402454,0.0338462) (0.0413377,0.0338462) (0.0424595,0.0338462) (0.0436118,0.0338462) (0.0447954,0.0338462) (0.0460111,0.0338462) (0.0472598,0.0338462) (0.0485424,0.0338462) (0.0498597,0.0338462) (0.0512129,0.0338462) (0.0526027,0.0338462) (0.0540303,0.0338462) (0.0554966,0.0338462) (0.0570028,0.0338462) (0.0585498,0.0338462) (0.0601387,0.0338462) (0.0617708,0.0338462) (0.0634472,0.0338462) (0.0651691,0.0338462) (0.0669377,0.0338462) (0.0687543,0.0338462) (0.0706203,0.0338462) (0.0725368,0.0338462) (0.0745054,0.0338462) (0.0765274,0.0338462) (0.0786042,0.0338462) (0.0807375,0.0338462) (0.0829286,0.0338462) (0.0851792,0.0338462) (0.0874908,0.0338462) (0.0898652,0.0338462) (0.0923041,0.0338462) (0.0948091,0.0338462) (0.0973821,0.0338462) (0.100025,0.0801231) (0.10274,0.0801231) (0.105528,0.0801231) (0.108392,0.0801231) (0.111333,0.0801231) (0.114355,0.0801231) (0.117458,0.0801231) (0.120646,0.0801231) (0.12392,0.0801231) (0.127283,0.0801231) (0.130737,0.0801231) (0.134286,0.0801231) (0.13793,0.0801231) (0.141673,0.0801231) (0.145518,0.0801231) (0.149467,0.0801231) (0.153524,0.0801231) (0.15769,0.0801231) (0.16197,0.0801231) (0.166365,0.0801231) (0.17088,0.0801231) (0.175518,0.0801231) (0.180281,0.0801231) (0.185174,0.0801231) (0.190199,0.0801231) (0.195361,0.0801231) (0.200663,0.120369) (0.206109,0.120369) (0.211702,0.120369) (0.217447,0.120369) (0.223349,0.120369) (0.22941,0.120369) (0.235636,0.120369) (0.242031,0.120369) (0.248599,0.120369) (0.255346,0.120369) (0.262276,0.120369) (0.269394,0.120369) (0.276705,0.120369) (0.284214,0.120369) (0.291928,0.120369) (0.29985,0.120369) (0.307988,0.144369) (0.316346,0.144369) (0.324932,0.144369) (0.33375,0.144369) (0.342807,0.144369) (0.352111,0.144369) (0.361667,0.144369) (0.371482,0.144369) (0.381564,0.144369) (0.391919,0.144369) (0.402555,0.180923) (0.41348,0.180923) (0.424701,0.180923) (0.436227,0.180923) (0.448066,0.180923) (0.460226,0.180923) (0.472716,0.180923) (0.485545,0.180923) (0.498722,0.180923) (0.512257,0.195815) (0.526159,0.195815) (0.540438,0.195815) (0.555105,0.195815) (0.57017,0.195815) (0.585644,0.195815) (0.601538,0.207877) (0.617863,0.207877) (0.634631,0.207877) (0.651854,0.207877) (0.669544,0.207877) (0.687715,0.207877) (0.706379,0.2224) (0.725549,0.2224) (0.74524,0.2224) (0.765465,0.2224) (0.786239,0.2224) (0.807576,0.240246) (0.829493,0.240246) (0.852005,0.240246) (0.875127,0.240246) (0.898877,0.240246) (0.923271,0.2576) (0.948328,0.2576) (0.974065,0.2576) (1.0005,0.269169) (1.02765,0.269169) (1.05554,0.269169) (1.08419,0.269169) (1.11361,0.272985) (1.14383,0.272985) (1.17488,0.272985) (1.20676,0.2768) (1.23951,0.2768) (1.27315,0.2768) (1.3077,0.2816) (1.34319,0.2816) (1.37964,0.2816) (1.41709,0.287385) (1.45554,0.287385) (1.49505,0.287385) (1.53562,0.290954) (1.57729,0.290954) (1.6201,0.296246) (1.66407,0.296246) (1.70923,0.299077) (1.75562,0.299077) (1.80326,0.302646) (1.8522,0.302646) (1.90247,0.306585) (1.9541,0.306585) (2.00713,0.309538) (2.0616,0.309538) (2.11755,0.312615) (2.17502,0.312615) (2.23405,0.315569) (2.29467,0.315569) (2.35695,0.319262) (2.42091,0.322954) (2.48662,0.322954) (2.5541,0.326523) (2.62341,0.329231) (2.69461,0.329231) (2.76774,0.3312) (2.84285,0.333785) (2.92001,0.335385) (2.99925,0.335385) (3.08065,0.337231) (3.16425,0.339692) (3.25013,0.347446) (3.33833,0.352492) (3.42893,0.356308) (3.52199,0.360862) (3.61757,0.365785) (3.71575,0.367754) (3.81659,0.369231) (3.92017,0.370585) (4.02656,0.376492) (4.13583,0.378708) (4.24807,0.381785) (4.36336,0.383508) (4.48178,0.385231) (4.60341,0.387938) (4.72834,0.390031) (4.85666,0.392369) (4.98847,0.394215) (5.12385,0.398646) (5.2629,0.400123) (5.40573,0.404062) (5.55244,0.406277) (5.70312,0.411323) (5.8579,0.413538) (6.01688,0.418954) (6.18017,0.423508) (6.34789,0.427692) (6.52017,0.431385) (6.69712,0.434708) (6.87887,0.442092) (7.06555,0.446277) (7.2573,0.451569) (7.45426,0.455754) (7.65656,0.459446) (7.86435,0.463631) (8.07778,0.467569) (8.297,0.472) (8.52217,0.478646) (8.75346,0.4816) (8.99101,0.485415) (9.23502,0.490585) (9.48565,0.494154) (9.74308,0.499323) (10.0075,0.505354) (10.2791,0.509292) (10.558,0.515938) (10.8446,0.521477) (11.1389,0.528369) (11.4412,0.534769) (11.7517,0.540923) (12.0706,0.547815) (12.3982,0.554585) (12.7347,0.565538) (13.0803,0.572923) (13.4353,0.583877) (13.7999,0.592) (14.1744,0.600492) (14.5591,0.610215) (14.9542,0.618708) (15.36,0.628308) (15.7769,0.6384) (16.2051,0.649846) (16.6448,0.661046) (17.0966,0.670892) (17.5605,0.686031) (18.0371,0.702769) (18.5266,0.715938) (19.0294,0.728369) (19.5459,0.739569) (20.0763,0.751262) (20.6212,0.764923) (21.1808,0.777846) (21.7556,0.791138) (22.346,0.807385) (22.9525,0.822892) (23.5754,0.8352) (24.2152,0.8576) (24.8724,0.873231) (25.5474,0.8912) (26.2407,0.905969) (26.9528,0.924185) (27.6843,0.940923) (28.4356,0.963446) (29.2073,0.981662) (30,1)};
\addlegendentry{FJ}
\addplot[const plot, black, line width=1pt, dashed] coordinates {(0.01,0.168508) (0.0102714,0.168508) (0.0105501,0.168508) (0.0108365,0.168508) (0.0111306,0.168508) (0.0114326,0.168508) (0.0117429,0.168508) (0.0120616,0.168508) (0.0123889,0.168508) (0.0127251,0.168508) (0.0130705,0.168508) (0.0134252,0.168508) (0.0137895,0.168508) (0.0141638,0.168508) (0.0145482,0.168508) (0.014943,0.168508) (0.0153485,0.168508) (0.0157651,0.168508) (0.0161929,0.168508) (0.0166324,0.168508) (0.0170838,0.168508) (0.0175474,0.168508) (0.0180236,0.168508) (0.0185127,0.168508) (0.0190152,0.168508) (0.0195312,0.168508) (0.0200613,0.168508) (0.0206057,0.168508) (0.0211649,0.168508) (0.0217393,0.168508) (0.0223293,0.168508) (0.0229353,0.168508) (0.0235577,0.168508) (0.0241971,0.168508) (0.0248537,0.168508) (0.0255282,0.168508) (0.026221,0.168508) (0.0269327,0.168508) (0.0276636,0.168508) (0.0284143,0.168508) (0.0291855,0.168508) (0.0299775,0.168508) (0.0307911,0.168508) (0.0316267,0.168508) (0.032485,0.168508) (0.0333666,0.168508) (0.0342722,0.168508) (0.0352023,0.168508) (0.0361576,0.168508) (0.0371389,0.168508) (0.0381468,0.168508) (0.0391821,0.168508) (0.0402454,0.168508) (0.0413377,0.168508) (0.0424595,0.168508) (0.0436118,0.168508) (0.0447954,0.168508) (0.0460111,0.168508) (0.0472598,0.168508) (0.0485424,0.168508) (0.0498597,0.168508) (0.0512129,0.168508) (0.0526027,0.168508) (0.0540303,0.168508) (0.0554966,0.168508) (0.0570028,0.168508) (0.0585498,0.168508) (0.0601387,0.168508) (0.0617708,0.168508) (0.0634472,0.168508) (0.0651691,0.168508) (0.0669377,0.168508) (0.0687543,0.168508) (0.0706203,0.168508) (0.0725368,0.168508) (0.0745054,0.168508) (0.0765274,0.168508) (0.0786042,0.168508) (0.0807375,0.168508) (0.0829286,0.168508) (0.0851792,0.168508) (0.0874908,0.168508) (0.0898652,0.168508) (0.0923041,0.168508) (0.0948091,0.168508) (0.0973821,0.168508) (0.100025,0.272099) (0.10274,0.272099) (0.105528,0.272099) (0.108392,0.272099) (0.111333,0.272099) (0.114355,0.272099) (0.117458,0.272099) (0.120646,0.272099) (0.12392,0.272099) (0.127283,0.272099) (0.130737,0.272099) (0.134286,0.272099) (0.13793,0.272099) (0.141673,0.272099) (0.145518,0.272099) (0.149467,0.272099) (0.153524,0.272099) (0.15769,0.272099) (0.16197,0.272099) (0.166365,0.272099) (0.17088,0.272099) (0.175518,0.272099) (0.180281,0.272099) (0.185174,0.272099) (0.190199,0.272099) (0.195361,0.272099) (0.200663,0.335635) (0.206109,0.335635) (0.211702,0.335635) (0.217447,0.335635) (0.223349,0.335635) (0.22941,0.335635) (0.235636,0.335635) (0.242031,0.335635) (0.248599,0.335635) (0.255346,0.335635) (0.262276,0.335635) (0.269394,0.335635) (0.276705,0.335635) (0.284214,0.335635) (0.291928,0.335635) (0.29985,0.335635) (0.307988,0.374309) (0.316346,0.374309) (0.324932,0.374309) (0.33375,0.374309) (0.342807,0.374309) (0.352111,0.374309) (0.361667,0.374309) (0.371482,0.374309) (0.381564,0.374309) (0.391919,0.374309) (0.402555,0.392265) (0.41348,0.392265) (0.424701,0.392265) (0.436227,0.392265) (0.448066,0.392265) (0.460226,0.392265) (0.472716,0.392265) (0.485545,0.392265) (0.498722,0.392265) (0.512257,0.425414) (0.526159,0.425414) (0.540438,0.425414) (0.555105,0.425414) (0.57017,0.425414) (0.585644,0.425414) (0.601538,0.451657) (0.617863,0.451657) (0.634631,0.451657) (0.651854,0.451657) (0.669544,0.451657) (0.687715,0.451657) (0.706379,0.483425) (0.725549,0.483425) (0.74524,0.483425) (0.765465,0.483425) (0.786239,0.483425) (0.807576,0.495856) (0.829493,0.495856) (0.852005,0.495856) (0.875127,0.495856) (0.898877,0.495856) (0.923271,0.505525) (0.948328,0.505525) (0.974065,0.505525) (1.0005,0.519337) (1.02765,0.519337) (1.05554,0.519337) (1.08419,0.519337) (1.11361,0.527624) (1.14383,0.527624) (1.17488,0.527624) (1.20676,0.542818) (1.23951,0.542818) (1.27315,0.542818) (1.3077,0.549724) (1.34319,0.549724) (1.37964,0.549724) (1.41709,0.555249) (1.45554,0.555249) (1.49505,0.555249) (1.53562,0.563536) (1.57729,0.563536) (1.6201,0.566298) (1.66407,0.566298) (1.70923,0.581492) (1.75562,0.581492) (1.80326,0.587017) (1.8522,0.587017) (1.90247,0.592541) (1.9541,0.592541) (2.00713,0.600829) (2.0616,0.600829) (2.11755,0.603591) (2.17502,0.603591) (2.23405,0.617403) (2.29467,0.617403) (2.35695,0.624309) (2.42091,0.631215) (2.48662,0.631215) (2.5541,0.639503) (2.62341,0.645028) (2.69461,0.645028) (2.76774,0.650552) (2.84285,0.657459) (2.92001,0.661602) (2.99925,0.661602) (3.08065,0.665746) (3.16425,0.668508) (3.25013,0.66989) (3.33833,0.674033) (3.42893,0.689227) (3.52199,0.694751) (3.61757,0.696133) (3.71575,0.700276) (3.81659,0.708564) (3.92017,0.709945) (4.02656,0.712707) (4.13583,0.718232) (4.24807,0.729282) (4.36336,0.737569) (4.48178,0.741713) (4.60341,0.744475) (4.72834,0.755525) (4.85666,0.758287) (4.98847,0.765193) (5.12385,0.767956) (5.2629,0.772099) (5.40573,0.776243) (5.55244,0.779006) (5.70312,0.78453) (5.8579,0.787293) (6.01688,0.791436) (6.18017,0.792818) (6.34789,0.798343) (6.52017,0.798343) (6.69712,0.802486) (6.87887,0.80663) (7.06555,0.813536) (7.2573,0.814917) (7.45426,0.819061) (7.65656,0.819061) (7.86435,0.824586) (8.07778,0.827348) (8.297,0.832873) (8.52217,0.835635) (8.75346,0.835635) (8.99101,0.838398) (9.23502,0.839779) (9.48565,0.84116) (9.74308,0.845304) (10.0075,0.849448) (10.2791,0.854972) (10.558,0.859116) (10.8446,0.861878) (11.1389,0.864641) (11.4412,0.866022) (11.7517,0.870166) (12.0706,0.875691) (12.3982,0.875691) (12.7347,0.877072) (13.0803,0.881215) (13.4353,0.889503) (13.7999,0.893646) (14.1744,0.893646) (14.5591,0.899171) (14.9542,0.903315) (15.36,0.906077) (15.7769,0.90884) (16.2051,0.90884) (16.6448,0.910221) (17.0966,0.911602) (17.5605,0.912983) (18.0371,0.91989) (18.5266,0.922652) (19.0294,0.929558) (19.5459,0.935083) (20.0763,0.94337) (20.6212,0.946133) (21.1808,0.947514) (21.7556,0.953039) (22.346,0.955801) (22.9525,0.961326) (23.5754,0.962707) (24.2152,0.96547) (24.8724,0.975138) (25.5474,0.980663) (26.2407,0.983425) (26.9528,0.987569) (27.6843,0.995856) (28.4356,0.995856) (29.2073,0.998619) (30,1)};
\addlegendentry{FPR}
\addplot[const plot, black!55, line width=1.1pt, densely dotted] coordinates {(0.01,0.0323722) (0.0102714,0.0323722) (0.0105501,0.0323722) (0.0108365,0.0323722) (0.0111306,0.0323722) (0.0114326,0.0323722) (0.0117429,0.0323722) (0.0120616,0.0323722) (0.0123889,0.0323722) (0.0127251,0.0323722) (0.0130705,0.0323722) (0.0134252,0.0323722) (0.0137895,0.0323722) (0.0141638,0.0323722) (0.0145482,0.0323722) (0.014943,0.0323722) (0.0153485,0.0323722) (0.0157651,0.0323722) (0.0161929,0.0323722) (0.0166324,0.0323722) (0.0170838,0.0323722) (0.0175474,0.0323722) (0.0180236,0.0323722) (0.0185127,0.0323722) (0.0190152,0.0323722) (0.0195312,0.0323722) (0.0200613,0.0323722) (0.0206057,0.0323722) (0.0211649,0.0323722) (0.0217393,0.0323722) (0.0223293,0.0323722) (0.0229353,0.0323722) (0.0235577,0.0323722) (0.0241971,0.0323722) (0.0248537,0.0323722) (0.0255282,0.0323722) (0.026221,0.0323722) (0.0269327,0.0323722) (0.0276636,0.0323722) (0.0284143,0.0323722) (0.0291855,0.0323722) (0.0299775,0.0323722) (0.0307911,0.0323722) (0.0316267,0.0323722) (0.032485,0.0323722) (0.0333666,0.0323722) (0.0342722,0.0323722) (0.0352023,0.0323722) (0.0361576,0.0323722) (0.0371389,0.0323722) (0.0381468,0.0323722) (0.0391821,0.0323722) (0.0402454,0.0323722) (0.0413377,0.0323722) (0.0424595,0.0323722) (0.0436118,0.0323722) (0.0447954,0.0323722) (0.0460111,0.0323722) (0.0472598,0.0323722) (0.0485424,0.0323722) (0.0498597,0.0323722) (0.0512129,0.0323722) (0.0526027,0.0323722) (0.0540303,0.0323722) (0.0554966,0.0323722) (0.0570028,0.0323722) (0.0585498,0.0323722) (0.0601387,0.0323722) (0.0617708,0.0323722) (0.0634472,0.0323722) (0.0651691,0.0323722) (0.0669377,0.0323722) (0.0687543,0.0323722) (0.0706203,0.0323722) (0.0725368,0.0323722) (0.0745054,0.0323722) (0.0765274,0.0323722) (0.0786042,0.0323722) (0.0807375,0.0323722) (0.0829286,0.0323722) (0.0851792,0.0323722) (0.0874908,0.0323722) (0.0898652,0.0323722) (0.0923041,0.0323722) (0.0948091,0.0323722) (0.0973821,0.0323722) (0.100025,0.0672048) (0.10274,0.0672048) (0.105528,0.0672048) (0.108392,0.0672048) (0.111333,0.0672048) (0.114355,0.0672048) (0.117458,0.0672048) (0.120646,0.0672048) (0.12392,0.0672048) (0.127283,0.0672048) (0.130737,0.0672048) (0.134286,0.0672048) (0.13793,0.0672048) (0.141673,0.0672048) (0.145518,0.0672048) (0.149467,0.0672048) (0.153524,0.0672048) (0.15769,0.0672048) (0.16197,0.0672048) (0.166365,0.0672048) (0.17088,0.0672048) (0.175518,0.0672048) (0.180281,0.0672048) (0.185174,0.0672048) (0.190199,0.0672048) (0.195361,0.0672048) (0.200663,0.0883546) (0.206109,0.0883546) (0.211702,0.0883546) (0.217447,0.0883546) (0.223349,0.0883546) (0.22941,0.0883546) (0.235636,0.0883546) (0.242031,0.0883546) (0.248599,0.0883546) (0.255346,0.0883546) (0.262276,0.0883546) (0.269394,0.0883546) (0.276705,0.0883546) (0.284214,0.0883546) (0.291928,0.0883546) (0.29985,0.0883546) (0.307988,0.107044) (0.316346,0.107044) (0.324932,0.107044) (0.33375,0.107044) (0.342807,0.107044) (0.352111,0.107044) (0.361667,0.107044) (0.371482,0.107044) (0.381564,0.107044) (0.391919,0.107044) (0.402555,0.119044) (0.41348,0.119044) (0.424701,0.119044) (0.436227,0.119044) (0.448066,0.119044) (0.460226,0.119044) (0.472716,0.119044) (0.485545,0.119044) (0.498722,0.119044) (0.512257,0.127676) (0.526159,0.127676) (0.540438,0.127676) (0.555105,0.127676) (0.57017,0.127676) (0.585644,0.127676) (0.601538,0.134194) (0.617863,0.134194) (0.634631,0.134194) (0.651854,0.134194) (0.669544,0.134194) (0.687715,0.134194) (0.706379,0.145718) (0.725549,0.145718) (0.74524,0.145718) (0.765465,0.145718) (0.786239,0.145718) (0.807576,0.156941) (0.829493,0.156941) (0.852005,0.156941) (0.875127,0.156941) (0.898877,0.156941) (0.923271,0.16566) (0.948328,0.16566) (0.974065,0.16566) (1.0005,0.174594) (1.02765,0.174594) (1.05554,0.174594) (1.08419,0.174594) (1.11361,0.184435) (1.14383,0.184435) (1.17488,0.184435) (1.20676,0.1935) (1.23951,0.1935) (1.27315,0.1935) (1.3077,0.202478) (1.34319,0.202478) (1.37964,0.202478) (1.41709,0.212146) (1.45554,0.212146) (1.49505,0.212146) (1.53562,0.22121) (1.57729,0.22121) (1.6201,0.232476) (1.66407,0.232476) (1.70923,0.241583) (1.75562,0.241583) (1.80326,0.249741) (1.8522,0.249741) (1.90247,0.258762) (1.9541,0.258762) (2.00713,0.26528) (2.0616,0.26528) (2.11755,0.272229) (2.17502,0.272229) (2.23405,0.278617) (2.29467,0.278617) (2.35695,0.285135) (2.42091,0.291911) (2.48662,0.291911) (2.5541,0.298602) (2.62341,0.304299) (2.69461,0.304299) (2.76774,0.309953) (2.84285,0.317723) (2.92001,0.323248) (2.99925,0.323248) (3.08065,0.328772) (3.16425,0.334211) (3.25013,0.340901) (3.33833,0.349404) (3.42893,0.356958) (3.52199,0.364166) (3.61757,0.371115) (3.71575,0.379359) (3.81659,0.387647) (3.92017,0.394078) (4.02656,0.399689) (4.13583,0.407156) (4.24807,0.41553) (4.36336,0.420925) (4.48178,0.425544) (4.60341,0.432968) (4.72834,0.435471) (4.85666,0.438838) (4.98847,0.442291) (5.12385,0.449111) (5.2629,0.451657) (5.40573,0.456664) (5.55244,0.459643) (5.70312,0.466549) (5.8579,0.470088) (6.01688,0.476649) (6.18017,0.479972) (6.34789,0.488519) (6.52017,0.495425) (6.69712,0.497842) (6.87887,0.501511) (7.06555,0.504834) (7.2573,0.50928) (7.45426,0.512819) (7.65656,0.516661) (7.86435,0.521064) (8.07778,0.525164) (8.297,0.529394) (8.52217,0.541048) (8.75346,0.550501) (8.99101,0.560083) (9.23502,0.574845) (9.48565,0.58434) (9.74308,0.601519) (10.0075,0.628065) (10.2791,0.658106) (10.558,0.696003) (10.8446,0.709988) (11.1389,0.725872) (11.4412,0.73908) (11.7517,0.753583) (12.0706,0.768949) (12.3982,0.783192) (12.7347,0.799119) (13.0803,0.813277) (13.4353,0.828427) (13.7999,0.836844) (14.1744,0.84634) (14.5591,0.852383) (14.9542,0.856742) (15.36,0.86067) (15.7769,0.864468) (16.2051,0.868612) (16.6448,0.871504) (17.0966,0.87582) (17.5605,0.881777) (18.0371,0.886611) (18.5266,0.891531) (19.0294,0.895373) (19.5459,0.899517) (20.0763,0.905257) (20.6212,0.912681) (21.1808,0.917213) (21.7556,0.922479) (22.346,0.92809) (22.9525,0.934349) (23.5754,0.938838) (24.2152,0.944277) (24.8724,0.949845) (25.5474,0.957528) (26.2407,0.963959) (26.9528,0.970218) (27.6843,0.977728) (28.4356,0.984979) (29.2073,0.991454) (30,1)};
\addlegendentry{LocalMIP}
\addplot[const plot, black!55, line width=1pt, dash pattern=on 4pt off 2pt on 1pt off 2pt] coordinates {(0.01,0.0927152) (0.0102714,0.0927152) (0.0105501,0.0927152) (0.0108365,0.0927152) (0.0111306,0.0927152) (0.0114326,0.0927152) (0.0117429,0.0927152) (0.0120616,0.0927152) (0.0123889,0.0927152) (0.0127251,0.0927152) (0.0130705,0.0927152) (0.0134252,0.0927152) (0.0137895,0.0927152) (0.0141638,0.0927152) (0.0145482,0.0927152) (0.014943,0.0927152) (0.0153485,0.0927152) (0.0157651,0.0927152) (0.0161929,0.0927152) (0.0166324,0.0927152) (0.0170838,0.0927152) (0.0175474,0.0927152) (0.0180236,0.0927152) (0.0185127,0.0927152) (0.0190152,0.0927152) (0.0195312,0.0927152) (0.0200613,0.0927152) (0.0206057,0.0927152) (0.0211649,0.0927152) (0.0217393,0.0927152) (0.0223293,0.0927152) (0.0229353,0.0927152) (0.0235577,0.0927152) (0.0241971,0.0927152) (0.0248537,0.0927152) (0.0255282,0.0927152) (0.026221,0.0927152) (0.0269327,0.0927152) (0.0276636,0.0927152) (0.0284143,0.0927152) (0.0291855,0.0927152) (0.0299775,0.0927152) (0.0307911,0.0927152) (0.0316267,0.0927152) (0.032485,0.0927152) (0.0333666,0.0927152) (0.0342722,0.0927152) (0.0352023,0.0927152) (0.0361576,0.0927152) (0.0371389,0.0927152) (0.0381468,0.0927152) (0.0391821,0.0927152) (0.0402454,0.0927152) (0.0413377,0.0927152) (0.0424595,0.0927152) (0.0436118,0.0927152) (0.0447954,0.0927152) (0.0460111,0.0927152) (0.0472598,0.0927152) (0.0485424,0.0927152) (0.0498597,0.0927152) (0.0512129,0.0927152) (0.0526027,0.0927152) (0.0540303,0.0927152) (0.0554966,0.0927152) (0.0570028,0.0927152) (0.0585498,0.0927152) (0.0601387,0.0927152) (0.0617708,0.0927152) (0.0634472,0.0927152) (0.0651691,0.0927152) (0.0669377,0.0927152) (0.0687543,0.0927152) (0.0706203,0.0927152) (0.0725368,0.0927152) (0.0745054,0.0927152) (0.0765274,0.0927152) (0.0786042,0.0927152) (0.0807375,0.0927152) (0.0829286,0.0927152) (0.0851792,0.0927152) (0.0874908,0.0927152) (0.0898652,0.0927152) (0.0923041,0.0927152) (0.0948091,0.0927152) (0.0973821,0.0927152) (0.100025,0.167219) (0.10274,0.167219) (0.105528,0.167219) (0.108392,0.167219) (0.111333,0.167219) (0.114355,0.167219) (0.117458,0.167219) (0.120646,0.167219) (0.12392,0.167219) (0.127283,0.167219) (0.130737,0.167219) (0.134286,0.167219) (0.13793,0.167219) (0.141673,0.167219) (0.145518,0.167219) (0.149467,0.167219) (0.153524,0.167219) (0.15769,0.167219) (0.16197,0.167219) (0.166365,0.167219) (0.17088,0.167219) (0.175518,0.167219) (0.180281,0.167219) (0.185174,0.167219) (0.190199,0.167219) (0.195361,0.167219) (0.200663,0.213576) (0.206109,0.213576) (0.211702,0.213576) (0.217447,0.213576) (0.223349,0.213576) (0.22941,0.213576) (0.235636,0.213576) (0.242031,0.213576) (0.248599,0.213576) (0.255346,0.213576) (0.262276,0.213576) (0.269394,0.213576) (0.276705,0.213576) (0.284214,0.213576) (0.291928,0.213576) (0.29985,0.213576) (0.307988,0.241722) (0.316346,0.241722) (0.324932,0.241722) (0.33375,0.241722) (0.342807,0.241722) (0.352111,0.241722) (0.361667,0.241722) (0.371482,0.241722) (0.381564,0.241722) (0.391919,0.241722) (0.402555,0.261589) (0.41348,0.261589) (0.424701,0.261589) (0.436227,0.261589) (0.448066,0.261589) (0.460226,0.261589) (0.472716,0.261589) (0.485545,0.261589) (0.498722,0.261589) (0.512257,0.279801) (0.526159,0.279801) (0.540438,0.279801) (0.555105,0.279801) (0.57017,0.279801) (0.585644,0.279801) (0.601538,0.299669) (0.617863,0.299669) (0.634631,0.299669) (0.651854,0.299669) (0.669544,0.299669) (0.687715,0.299669) (0.706379,0.324503) (0.725549,0.324503) (0.74524,0.324503) (0.765465,0.324503) (0.786239,0.324503) (0.807576,0.337748) (0.829493,0.337748) (0.852005,0.337748) (0.875127,0.337748) (0.898877,0.337748) (0.923271,0.349338) (0.948328,0.349338) (0.974065,0.349338) (1.0005,0.359272) (1.02765,0.359272) (1.05554,0.359272) (1.08419,0.359272) (1.11361,0.36755) (1.14383,0.36755) (1.17488,0.36755) (1.20676,0.374172) (1.23951,0.374172) (1.27315,0.374172) (1.3077,0.392384) (1.34319,0.392384) (1.37964,0.392384) (1.41709,0.402318) (1.45554,0.402318) (1.49505,0.402318) (1.53562,0.418874) (1.57729,0.418874) (1.6201,0.425497) (1.66407,0.425497) (1.70923,0.432119) (1.75562,0.432119) (1.80326,0.44702) (1.8522,0.44702) (1.90247,0.455298) (1.9541,0.455298) (2.00713,0.466887) (2.0616,0.466887) (2.11755,0.471854) (2.17502,0.471854) (2.23405,0.481788) (2.29467,0.481788) (2.35695,0.485099) (2.42091,0.491722) (2.48662,0.491722) (2.5541,0.498344) (2.62341,0.511589) (2.69461,0.511589) (2.76774,0.516556) (2.84285,0.521523) (2.92001,0.528146) (2.99925,0.528146) (3.08065,0.536424) (3.16425,0.546358) (3.25013,0.551325) (3.33833,0.556291) (3.42893,0.557947) (3.52199,0.562914) (3.61757,0.571192) (3.71575,0.576159) (3.81659,0.57947) (3.92017,0.582781) (4.02656,0.586093) (4.13583,0.592715) (4.24807,0.597682) (4.36336,0.597682) (4.48178,0.599338) (4.60341,0.60596) (4.72834,0.612583) (4.85666,0.615894) (4.98847,0.61755) (5.12385,0.630795) (5.2629,0.640728) (5.40573,0.655629) (5.55244,0.660596) (5.70312,0.663907) (5.8579,0.667219) (6.01688,0.673841) (6.18017,0.678808) (6.34789,0.68543) (6.52017,0.693709) (6.69712,0.69702) (6.87887,0.700331) (7.06555,0.708609) (7.2573,0.715232) (7.45426,0.720199) (7.65656,0.726821) (7.86435,0.735099) (8.07778,0.736755) (8.297,0.740066) (8.52217,0.751656) (8.75346,0.753311) (8.99101,0.756623) (9.23502,0.764901) (9.48565,0.766556) (9.74308,0.769868) (10.0075,0.773179) (10.2791,0.774834) (10.558,0.779801) (10.8446,0.786424) (11.1389,0.791391) (11.4412,0.796358) (11.7517,0.799669) (12.0706,0.801325) (12.3982,0.804636) (12.7347,0.806291) (13.0803,0.812914) (13.4353,0.816225) (13.7999,0.817881) (14.1744,0.822848) (14.5591,0.826159) (14.9542,0.82947) (15.36,0.836093) (15.7769,0.836093) (16.2051,0.84106) (16.6448,0.842715) (17.0966,0.849338) (17.5605,0.854305) (18.0371,0.859272) (18.5266,0.86755) (19.0294,0.872517) (19.5459,0.877483) (20.0763,0.880795) (20.6212,0.895695) (21.1808,0.900662) (21.7556,0.903974) (22.346,0.90894) (22.9525,0.910596) (23.5754,0.918874) (24.2152,0.927152) (24.8724,0.945364) (25.5474,0.950331) (26.2407,0.961921) (26.9528,0.970199) (27.6843,0.975166) (28.4356,0.980132) (29.2073,0.986755) (30,1)};
\addlegendentry{Scylla}
\end{groupplot}
\end{tikzpicture}
\\[2pt]
\pgfplotslegendfromname{probetimelegend}
\caption{When the probe's heuristics improve, each running alone and presolve-only under a $30$\,s wall cap. Left: the mean primal gap $p(t)$ over the instances on which that heuristic found a feasible solution, so each instance is scaled by its own reference objective and $p=2$ before the first incumbent. Right: the share of that heuristic's incumbents that have arrived by time $t$. Only incumbents carrying the heuristic's own source tag are counted.}
\label{fig:probe-time}
\end{figure}

%% file: arms.tex

\begin{table}[t]
\centering
\small
\setlength{\tabcolsep}{4pt}
\resizebox{\ifdim\width>\linewidth\linewidth\else\width\fi}{!}{%
\begin{tabular}{lllrrrrrl}
\toprule
Arm & What it is & Race ($\lambda$, id) & $e_{\text{fj}}$ & $e_{\text{fpr}}$ & $e_{\text{lm}}$ & $e_{\text{scylla}}$ & Total & Stages \\
\midrule
\texttt{A-fj-only} & FJ alone, the pick at $\lambda=1/600$ & free, $1/600$, \#335 & 0.1619 & 0 & 0 & 0 & 0.1619 & confirmation \\
\texttt{B-mix-lambda300} & the generous mix & free, $1/300$, \#351 & 4.027 & 2.66 & 4.088 & 0 & 10.78 & confirmation, held-out check \\
\texttt{B'-mix-cheapest} & the cheap mix, selected and shipped & free, $1/600$, \#337 & 0.3317 & 3.161 & 3.287 & 0 & 6.779 & confirmation, held-out check \\
\texttt{D-shipped} & the measured vector, the control & not raced & 0.5665 & 12.26 & 13.96 & 3.068 & 29.85 & confirmation, held-out check \\
\texttt{D'-tuned-all4} & winner of the race forced to run all four & all four, $1/600$, \#288 & 0.3047 & 2.89 & 6.869 & 0.4907 & 10.55 & confirmation \\
\bottomrule
\end{tabular}
}
\par\vspace{7pt}
\caption{The configurations Stage~3 ran at the full limit, with the effort budget $e_h$ each gave every heuristic, their total, and the stages it ran in. \emph{Race}: the race the arm is a candidate of, free or with all four heuristics forced on, its cost weight $\lambda$ and the candidate's irace id. The names are those of \texttt{bench/finalists.json} in the code repository; the prose uses the description beside them.}
\label{tab:arms}
\end{table}

%% file: fig/forest.tex

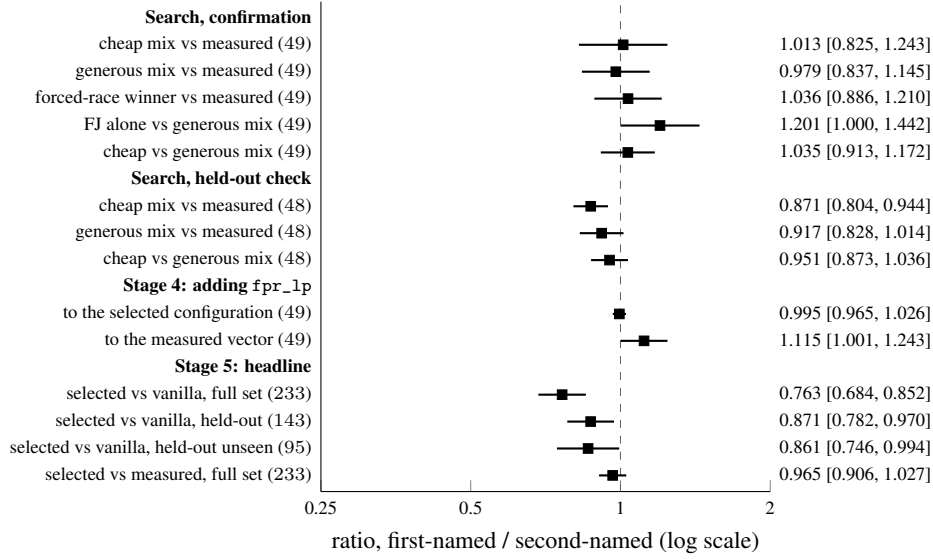
\begin{figure}[t]
\centering
\begin{tikzpicture}
\begin{axis}[
  scale only axis, width=0.36\linewidth, height=6.48cm,
  xmode=log, xmin=0.25, xmax=2,
  ymin=0.4, ymax=18.6,
  ytick={18,17,16,15,14,13,12,11,10,9,8,7,6,5,4,3,2,1}, yticklabels={{\textbf{Search, confirmation}},{cheap mix vs measured ($49$)},{generous mix vs measured ($49$)},{forced-race winner vs measured ($49$)},{FJ alone vs generous mix ($49$)},{cheap vs generous mix ($49$)},{\textbf{Search, held-out check}},{cheap mix vs measured ($48$)},{generous mix vs measured ($48$)},{cheap vs generous mix ($48$)},{\textbf{Stage 4: adding \texttt{fpr\_lp}}},{to the selected configuration ($49$)},{to the measured vector ($49$)},{\textbf{Stage 5: headline}},{selected vs vanilla, full set ($233$)},{selected vs vanilla, held-out ($143$)},{selected vs vanilla, held-out unseen ($95$)},{selected vs measured, full set ($233$)}},
  yticklabel style={font=\scriptsize, align=right}, ytick style={draw=none},
  xtick={0.25,0.5,1,2}, xticklabels={0.25,0.5,1,2},
  tick label style={font=\scriptsize}, label style={font=\small},
  xlabel={ratio, first-named / second-named (log scale)},
  clip=false, unbounded coords=jump, axis y line*=left, axis x line*=bottom,
]
\draw[dashed, black!60] (axis cs:1,0.4) -- (axis cs:1,18.6);
\addplot[black, line width=0.8pt] coordinates {(0.825126,17) (1.24267,17) (nan,nan) (0.836655,16) (1.14485,16) (nan,nan) (0.886245,15) (1.21043,15) (nan,nan) (0.999842,14) (1.44194,14) (nan,nan) (0.913138,13) (1.17232,13) (nan,nan) (0.804428,11) (0.944024,11) (nan,nan) (0.828415,10) (1.01413,10) (nan,nan) (0.872853,9) (1.03558,9) (nan,nan) (0.96475,7) (1.02599,7) (nan,nan) (1.00065,6) (1.24312,6) (nan,nan) (0.683956,4) (0.852052,4) (nan,nan) (0.7819,3) (0.970438,3) (nan,nan) (0.745513,2) (0.994125,2) (nan,nan) (0.905807,1) (1.0271,1) (nan,nan)};
\addplot[only marks, mark=square*, mark size=1.8pt, black] coordinates {(1.0126,17) (0.978694,16) (1.03573,15) (1.20071,14) (1.03464,13) (0.871435,11) (0.916583,10) (0.950743,9) (0.994898,7) (1.11531,6) (0.763391,4) (0.871083,3) (0.860891,2) (0.964549,1)};
\node[anchor=west, font=\scriptsize] at (axis cs:2,17) {1.013 [0.825, 1.243]};
\node[anchor=west, font=\scriptsize] at (axis cs:2,16) {0.979 [0.837, 1.145]};
\node[anchor=west, font=\scriptsize] at (axis cs:2,15) {1.036 [0.886, 1.210]};
\node[anchor=west, font=\scriptsize] at (axis cs:2,14) {1.201 [1.000, 1.442]};
\node[anchor=west, font=\scriptsize] at (axis cs:2,13) {1.035 [0.913, 1.172]};
\node[anchor=west, font=\scriptsize] at (axis cs:2,11) {0.871 [0.804, 0.944]};
\node[anchor=west, font=\scriptsize] at (axis cs:2,10) {0.917 [0.828, 1.014]};
\node[anchor=west, font=\scriptsize] at (axis cs:2,9) {0.951 [0.873, 1.036]};
\node[anchor=west, font=\scriptsize] at (axis cs:2,7) {0.995 [0.965, 1.026]};
\node[anchor=west, font=\scriptsize] at (axis cs:2,6) {1.115 [1.001, 1.243]};
\node[anchor=west, font=\scriptsize] at (axis cs:2,4) {0.763 [0.684, 0.852]};
\node[anchor=west, font=\scriptsize] at (axis cs:2,3) {0.871 [0.782, 0.970]};
\node[anchor=west, font=\scriptsize] at (axis cs:2,2) {0.861 [0.746, 0.994]};
\node[anchor=west, font=\scriptsize] at (axis cs:2,1) {0.965 [0.906, 1.027]};
\end{axis}
\end{tikzpicture}
\caption{The paired comparisons of the search (Stage~3), the LP-guided dive (Stage~4) and the headline (Stage~5), all in full-limit solves, on the primal integral: ratio of the shifted geometric means with its $95\%$ confidence interval, number of instances in parentheses. Below $1$ favours the first-named arm; an interval that crosses the dashed line does not separate the two. The cheap mix is the selected configuration, the measured vector is the probe's, used as the control, and the forced-race winner is the best of the race that had to run all four heuristics. Held-out unseen is the held-out complement without its held-out-check instances, the instances no stage saw before the headline.}
\label{fig:forest}
\end{figure}

%% file: headline.tex

\begin{table}[t]
\centering
\small
\setlength{\tabcolsep}{4pt}
\resizebox{\ifdim\width>\linewidth\linewidth\else\width\fi}{!}{%
\begin{tabular}{llrrrrrlrr}
\toprule
Set ($n$) & Arm & \#Feas. & \#Wins & T1st ($s{=}1$) & Gap ($s{=}0.001$) & \textbf{PI} ($s{=}0.001$) & PI ratio [95\% CI] & $p$ & B/T/W \\
\midrule
\multirow{2}{*}{Full set (233)} & selected & 214 & 49 & 5.57 & 0.0054 & 0.0484 & 0.763 [0.684, 0.852] & <0.001 & 114/24/95 \\
 & vanilla & 211 & 36 & 7.08 & 0.0066 & 0.0637 &  &  &  \\
\midrule
\multirow{2}{*}{Tuning set (90)} & selected & 86 & 19 & 3.06 & 0.0043 & 0.0344 & 0.619 [0.496, 0.772] & <0.001 & 51/5/34 \\
 & vanilla & 86 & 15 & 4.32 & 0.0056 & 0.0562 &  &  &  \\
\midrule
\multirow{2}{*}{Held-out (143)} & selected & 128 & 30 & 7.89 & 0.0062 & 0.0600 & 0.871 [0.782, 0.970] & 0.012 & 63/19/61 \\
 & vanilla & 125 & 21 & 9.52 & 0.0072 & 0.0690 &  &  &  \\
\bottomrule
\end{tabular}
}
\par\vspace{7pt}
\caption{The selected configuration versus a separately built unpatched HiGHS, HiGHS \texttt{v1.15.1}, \texttt{mipfeas} at $600$\,s, 16 workers, on the full set, the tuning set and the held-out complement. Feasible and wins (strictly best objective at the limit) are counts. Every SGM is over all instances: a run with no solution by the limit counts at the limit in T1st and scores $2$ in gap and PI. The last three columns are paired per instance on the primal integral: the ratio selected/vanilla with its 95\% confidence interval, the $p$-value, and the better/tied/worse counts. Lower is better throughout.}
\label{tab:headline}
\end{table}

%% file: attribution.tex

\begin{table}[t]
\centering
\small
\setlength{\tabcolsep}{4pt}
\resizebox{\ifdim\width>\linewidth\linewidth\else\width\fi}{!}{%
\begin{tabular}{lrr}
\toprule
Source & \#First feasible & \#Best at 600\,s \\
\midrule
FPR (\texttt{A}) & 16 & 2 \\
LocalMIP (\texttt{M}) & 5 & 4 \\
FJ (\texttt{J}) & 114 & 20 \\
HiGHS/other & 79 & 188 \\
\bottomrule
\end{tabular}
}
\par\smallskip
{\footnotesize Our FJ and HiGHS's own Feasibility Jump both emit the source tag \texttt{J}. HiGHS's own dispatch is patched off in this arm, so every \texttt{J} counted here is ours.}
\par\vspace{7pt}
\caption{Per-heuristic attribution over the instances the selected configuration made feasible, of 233, one seed, read from the solution-source tags of the log: the source of the first feasible solution and of the best solution at the limit, per instance.}
\label{tab:attribution}
\end{table}

%% file: fig/gap-time.tex

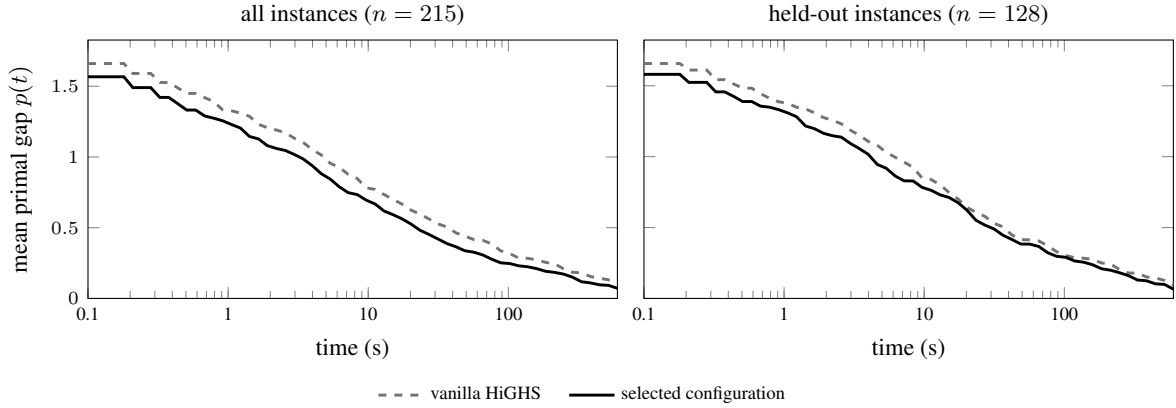
\begin{figure}[t]
\centering
\begin{tikzpicture}
\begin{groupplot}[
  group style={group size=2 by 1, horizontal sep=10pt,
    y descriptions at=edge left},
  width=0.52\linewidth, height=5cm,
  xmode=log, xmin=0.1, xmax=600, ymin=0,
  xtick={0.1,1,10,100}, xticklabels={0.1,1,10,100},
  tick label style={font=\scriptsize}, label style={font=\small},
  title style={font=\small, yshift=-4pt}, xlabel={time (s)},
]
\nextgroupplot[title={all instances ($n=215$)}, ylabel={mean primal gap $p(t)$}]
\addplot[black!55, line width=1.1pt, dashed] coordinates {(0.1,1.66062) (0.115887,1.66062) (0.134299,1.66062) (0.155636,1.66062) (0.180362,1.66062) (0.209017,1.5904) (0.242225,1.5904) (0.280708,1.5904) (0.325306,1.5267) (0.376989,1.5267) (0.436882,1.48871) (0.506292,1.44924) (0.586729,1.44924) (0.679945,1.41976) (0.787971,1.39223) (0.91316,1.33374) (1.05824,1.32648) (1.22637,1.30992) (1.4212,1.2872) (1.647,1.23088) (1.90866,1.20672) (2.2119,1.19063) (2.56332,1.17122) (2.97056,1.13266) (3.44251,1.10136) (3.98944,1.0523) (4.62326,1.01393) (5.35777,0.955206) (6.20899,0.927394) (7.19544,0.877121) (8.33861,0.84801) (9.66341,0.782084) (11.1987,0.770841) (12.9779,0.734116) (15.0397,0.700886) (17.4291,0.666321) (20.1982,0.624454) (23.4072,0.59594) (27.126,0.561899) (31.4356,0.527313) (36.4299,0.504887) (42.2177,0.463834) (48.9251,0.440932) (56.698,0.415979) (65.7059,0.410449) (76.1449,0.385499) (88.2424,0.33642) (102.262,0.317575) (118.509,0.285619) (137.337,0.283584) (159.156,0.273164) (184.442,0.254881) (213.745,0.24482) (247.704,0.209606) (287.058,0.185206) (332.664,0.181035) (385.516,0.154843) (446.764,0.143331) (517.744,0.131608) (600,0.110633)};
\addplot[black, line width=1.1pt] coordinates {(0.1,1.56699) (0.115887,1.56699) (0.134299,1.56699) (0.155636,1.56699) (0.180362,1.56699) (0.209017,1.4904) (0.242225,1.4904) (0.280708,1.4904) (0.325306,1.42108) (0.376989,1.42108) (0.436882,1.37694) (0.506292,1.33184) (0.586729,1.33184) (0.679945,1.29016) (0.787971,1.27438) (0.91316,1.25615) (1.05824,1.22923) (1.22637,1.20415) (1.4212,1.14483) (1.647,1.12764) (1.90866,1.07993) (2.2119,1.06015) (2.56332,1.04589) (2.97056,1.01753) (3.44251,0.985442) (3.98944,0.938364) (4.62326,0.88148) (5.35777,0.843784) (6.20899,0.789467) (7.19544,0.747766) (8.33861,0.734205) (9.66341,0.69509) (11.1987,0.667705) (12.9779,0.617727) (15.0397,0.591911) (17.4291,0.562599) (20.1982,0.525974) (23.4072,0.481926) (27.126,0.452428) (31.4356,0.420402) (36.4299,0.387911) (42.2177,0.365312) (48.9251,0.337226) (56.698,0.326864) (65.7059,0.307395) (76.1449,0.277031) (88.2424,0.251618) (102.262,0.246748) (118.509,0.230672) (137.337,0.224011) (159.156,0.210744) (184.442,0.191999) (213.745,0.184148) (247.704,0.17207) (287.058,0.150439) (332.664,0.117784) (385.516,0.109923) (446.764,0.0962019) (517.744,0.0927809) (600,0.0715937)};
\nextgroupplot[title={held-out instances ($n=128$)}, legend to name=gaptimelegend, legend columns=2, legend style={font=\scriptsize, draw=none, /tikz/every even column/.append style={column sep=8pt}}]
\addplot[black!55, line width=1.1pt, dashed] coordinates {(0.1,1.65673) (0.115887,1.65673) (0.134299,1.65673) (0.155636,1.65673) (0.180362,1.65673) (0.209017,1.61068) (0.242225,1.61068) (0.280708,1.61068) (0.325306,1.5429) (0.376989,1.5429) (0.436882,1.51224) (0.506292,1.48296) (0.586729,1.48296) (0.679945,1.44381) (0.787971,1.41136) (0.91316,1.38724) (1.05824,1.37623) (1.22637,1.34859) (1.4212,1.33717) (1.647,1.30866) (1.90866,1.27604) (2.2119,1.25692) (2.56332,1.2313) (2.97056,1.18695) (3.44251,1.15095) (3.98944,1.1074) (4.62326,1.08) (5.35777,1.03084) (6.20899,0.995589) (7.19544,0.961117) (8.33861,0.924467) (9.66341,0.849319) (11.1987,0.838643) (12.9779,0.802176) (15.0397,0.747188) (17.4291,0.701449) (20.1982,0.642263) (23.4072,0.607315) (27.126,0.578758) (31.4356,0.529486) (36.4299,0.50568) (42.2177,0.447555) (48.9251,0.415116) (56.698,0.413035) (65.7059,0.405555) (76.1449,0.371217) (88.2424,0.328934) (102.262,0.303715) (118.509,0.286486) (137.337,0.283884) (159.156,0.275065) (184.442,0.249995) (213.745,0.242703) (247.704,0.204148) (287.058,0.179606) (332.664,0.178466) (385.516,0.151422) (446.764,0.145332) (517.744,0.127167) (600,0.0932476)};
\addlegendentry{vanilla HiGHS}
\addplot[black, line width=1.1pt] coordinates {(0.1,1.5805) (0.115887,1.5805) (0.134299,1.5805) (0.155636,1.5805) (0.180362,1.5805) (0.209017,1.52386) (0.242225,1.52386) (0.280708,1.52386) (0.325306,1.45688) (0.376989,1.45688) (0.436882,1.42471) (0.506292,1.38881) (0.586729,1.38881) (0.679945,1.35624) (0.787971,1.34879) (0.91316,1.33176) (1.05824,1.30907) (1.22637,1.2819) (1.4212,1.21628) (1.647,1.19733) (1.90866,1.16473) (2.2119,1.14767) (2.56332,1.13733) (2.97056,1.09334) (3.44251,1.05691) (3.98944,1.01633) (4.62326,0.94518) (5.35777,0.918031) (6.20899,0.863508) (7.19544,0.829462) (8.33861,0.828588) (9.66341,0.783191) (11.1987,0.762778) (12.9779,0.730889) (15.0397,0.711952) (17.4291,0.675456) (20.1982,0.622284) (23.4072,0.551357) (27.126,0.516988) (31.4356,0.490454) (36.4299,0.446145) (42.2177,0.41606) (48.9251,0.383723) (56.698,0.382967) (65.7059,0.366623) (76.1449,0.322327) (88.2424,0.296989) (102.262,0.290342) (118.509,0.264153) (137.337,0.255498) (159.156,0.23766) (184.442,0.209311) (213.745,0.199199) (247.704,0.181235) (287.058,0.160821) (332.664,0.12987) (385.516,0.12458) (446.764,0.103081) (517.744,0.0991288) (600,0.0653052)};
\addlegendentry{selected configuration}
\end{groupplot}
\end{tikzpicture}
\\[2pt]
\pgfplotslegendfromname{gaptimelegend}
\caption{Mean primal gap $p(t)$ over time, log time axis, selected configuration against vanilla HiGHS; $p=2$ before the first incumbent, and the time average of a curve over the limit is the mean primal integral. Instances on which neither arm found a solution are left out (18 in all, 15 of them held-out).}
\label{fig:gap-time}
\end{figure}

%% file: appendix.tex
\appendix

\section{Per-instance headline}
\label{app:per-instance}

Table~\ref{tab:per-instance} lists every \texttt{mipfeas} instance with the selected
configuration and vanilla HiGHS side by side: time to first feasible solution, final gap and
primal integral for each arm, and their ratio on the primal integral.

\input{appendix_per_instance}

%% file: appendix_per_instance.tex

{\scriptsize
\setlength{\tabcolsep}{1.5pt}
\begin{longtable}{llrrrrrrr}
\toprule
Instance & Set & T1st (sel.) & T1st (van.) & Gap (sel.) & Gap (van.) & PI (sel.) & PI (van.) & Ratio \\
\midrule
\endfirsthead
\toprule
Instance & Set & T1st (sel.) & T1st (van.) & Gap (sel.) & Gap (van.) & PI (sel.) & PI (van.) & Ratio \\
\midrule
\endhead
\midrule \multicolumn{9}{r}{\emph{continued}}\\
\endfoot
\bottomrule
\caption{Per-instance headline, selected configuration (sel.) and vanilla HiGHS (van.): all 233 \texttt{mipfeas} instances at $600$\,s. T1st in seconds, gap and ratio as fractions, \textsc{n/a} where the arm never became feasible. The ratio is selected/vanilla on the primal integral (PI); below $1$ favours the selected configuration. Set: T tuning set, H held-out complement; C also in the confirmation set, K also in the held-out check.}\label{tab:per-instance}
\endlastfoot
\texttt{30n20b8} & H & 39.6 & 37.3 & 0.0000 & 0.0000 & 0.1764 & 0.1754 & 1.006 \\
\texttt{50v-10} & H & 0.0 & 0.0 & 0.0065 & 0.0064 & 0.0208 & 0.0274 & 0.770 \\
\texttt{CMS750\_4} & TC & 0.2 & 0.2 & 0.0079 & 0.0345 & 0.0425 & 0.1527 & 0.283 \\
\texttt{academictimetablesmall} & H & 1.4 & 1.6 & 1.0000 & 1.0000 & 1.0023 & 1.0027 & 1.000 \\
\texttt{air05} & T & 0.4 & 0.4 & 0.0000 & 0.0000 & 0.0072 & 0.0079 & 0.927 \\
\texttt{app1-1} & H & 0.2 & 0.2 & 0.0000 & 0.0000 & 0.0352 & 0.0292 & 1.199 \\
\texttt{app1-2} & H & n/a & n/a & n/a & n/a & 2.0000 & 2.0000 & 1.000 \\
\texttt{assign1-5-8} & TC & 0.0 & 0.0 & 0.0000 & 0.0000 & 0.0000 & 0.0002 & 0.843 \\
\texttt{atlanta-ip} & TC & 295.7 & 1.5 & 0.0323 & 0.0425 & 1.0151 & 0.1668 & 6.055 \\
\texttt{b1c1s1} & TC & 0.3 & 0.3 & 0.0170 & 0.0281 & 0.0777 & 0.0798 & 0.974 \\
\texttt{bab2} & H & 287.4 & 263.2 & 0.0031 & 0.0025 & 0.9663 & 0.8848 & 1.092 \\
\texttt{bab6} & H & 259.5 & 239.7 & 0.0131 & 0.0131 & 0.8724 & 0.8069 & 1.081 \\
\texttt{beasleyC3} & TC & 0.0 & 0.0 & 0.0000 & 0.0000 & 0.0005 & 0.0004 & 1.018 \\
\texttt{binkar10\_1} & HK & 0.0 & 0.0 & 0.0000 & 0.0000 & 0.0001 & 0.0001 & 1.008 \\
\texttt{blp-ar98} & H & 1.3 & 490.1 & 0.0058 & 0.0107 & 0.0870 & 1.6357 & 0.054 \\
\texttt{blp-ic98} & H & 1.1 & 72.8 & 0.0075 & 0.0138 & 0.0402 & 0.2666 & 0.154 \\
\texttt{bnatt400} & HK & 560.6 & 542.8 & 0.0000 & 0.0000 & 1.8687 & 1.8093 & 1.033 \\
\texttt{bppc4-08} & TC & 0.1 & 0.1 & 0.0000 & 0.0000 & 0.0012 & 0.0036 & 0.465 \\
\texttt{brazil3} & H & 0.3 & 0.3 & 0.8961 & 0.8723 & 0.9125 & 0.9197 & 0.992 \\
\texttt{buildingenergy} & H & 168.2 & 162.1 & 0.0204 & 0.0204 & 0.5753 & 0.5552 & 1.036 \\
\texttt{cbs-cta} & TC & 0.4 & 0.3 & 0.0000 & 0.0000 & 0.0737 & 0.0867 & 0.852 \\
\texttt{chromaticindex1024-7} & H & 5.0 & 5.1 & 0.0000 & 0.0000 & 0.0167 & 0.0170 & 0.981 \\
\texttt{chromaticindex512-7} & H & 1.6 & 1.8 & 0.0000 & 0.0000 & 0.0053 & 0.0060 & 0.905 \\
\texttt{cmflsp50-24-8-8} & H & 17.2 & 15.9 & 0.0193 & 0.0245 & 0.0930 & 0.0899 & 1.034 \\
\texttt{co-100} & H & 411.7 & 358.5 & 0.3113 & 0.3373 & 1.5070 & 1.4069 & 1.071 \\
\texttt{cod105} & H & 0.8 & 0.9 & 0.0000 & 0.0000 & 0.0027 & 0.0457 & 0.078 \\
\texttt{comp07-2idx} & T & 0.9 & 0.9 & 0.0000 & 0.7692 & 0.0482 & 0.8350 & 0.059 \\
\texttt{comp21-2idx} & T & 0.6 & 0.6 & 0.3782 & 0.5698 & 0.4491 & 0.6652 & 0.676 \\
\texttt{cost266-UUE} & TC & 0.4 & 0.3 & 0.0152 & 0.0037 & 0.0431 & 0.0414 & 1.041 \\
\texttt{cryptanalysiskb128n5obj16} & HK & n/a & n/a & n/a & n/a & 2.0000 & 2.0000 & 1.000 \\
\texttt{csched007} & H & 44.1 & 40.3 & 0.0028 & 0.0223 & 0.2207 & 0.2321 & 0.951 \\
\texttt{csched008} & HK & 20.5 & 17.6 & 0.0000 & 0.0000 & 0.0728 & 0.0626 & 1.160 \\
\texttt{cvs16r128-89} & TC & 0.1 & 0.1 & 0.0000 & 0.0000 & 0.0075 & 0.1084 & 0.078 \\
\texttt{dano3\_3} & HK & 19.2 & 17.8 & 0.0000 & 0.0000 & 0.0640 & 0.0593 & 1.077 \\
\texttt{dano3\_5} & H & 17.5 & 15.4 & 0.0000 & 0.0000 & 0.0591 & 0.0521 & 1.132 \\
\texttt{decomp2} & T & 0.3 & 0.3 & 0.0000 & 0.0000 & 0.0015 & 0.0065 & 0.341 \\
\texttt{drayage-100-23} & T & 0.1 & 0.1 & 0.0000 & 0.0000 & 0.0008 & 0.0008 & 1.027 \\
\texttt{drayage-25-23} & HK & 0.1 & 0.1 & 0.0000 & 0.0000 & 0.0016 & 0.0013 & 1.135 \\
\texttt{dws008-01} & HK & 43.5 & 37.7 & 0.1698 & 0.0925 & 0.4404 & 0.4166 & 1.057 \\
\texttt{eil33-2} & H & 1.7 & 1.5 & 0.0000 & 0.0000 & 0.0105 & 0.0131 & 0.814 \\
\texttt{eilA101-2} & HK & 84.0 & 83.5 & 0.0000 & 0.0000 & 0.3338 & 0.4503 & 0.742 \\
\texttt{enlight\_hard} & T & 1.8 & 1.6 & 0.0000 & 0.0000 & 0.0060 & 0.0053 & 1.105 \\
\texttt{ex10} & H & 26.6 & 26.9 & 0.0000 & 0.0000 & 0.0887 & 0.0897 & 0.989 \\
\texttt{ex9} & H & 3.5 & 3.2 & 0.0000 & 0.0000 & 0.0117 & 0.0107 & 1.086 \\
\texttt{exp-1-500-5-5} & HK & 0.1 & 0.1 & 0.0000 & 0.0000 & 0.0016 & 0.0016 & 1.000 \\
\texttt{fast0507} & H & 1.3 & 1.2 & 0.0057 & 0.0057 & 0.2379 & 0.0223 & 10.257 \\
\texttt{fastxgemm-n2r6s0t2} & H & 0.0 & 0.0 & 0.0000 & 0.0000 & 0.2117 & 0.0791 & 2.656 \\
\texttt{fhnw-binpack4-48} & T & n/a & n/a & n/a & n/a & 2.0000 & 2.0000 & 1.000 \\
\texttt{fiball} & H & 1.0 & 2.1 & 0.0000 & 0.0000 & 0.0058 & 0.0128 & 0.489 \\
\texttt{gen-ip002} & HK & 0.0 & 0.0 & 0.0016 & 0.0000 & 0.0021 & 0.0012 & 1.407 \\
\texttt{gen-ip054} & T & 0.0 & 0.0 & 0.0000 & 0.0000 & 0.0009 & 0.0011 & 0.949 \\
\texttt{germanrr} & TC & n/a & n/a & n/a & n/a & 2.0000 & 2.0000 & 1.000 \\
\texttt{gfd-schedulen180f7d50m30k18} & H & n/a & n/a & n/a & n/a & 2.0000 & 2.0000 & 1.000 \\
\texttt{glass-sc} & H & 0.1 & 0.0 & 0.0000 & 0.0000 & 0.0006 & 0.0082 & 0.175 \\
\texttt{glass4} & H & 0.0 & 0.6 & 0.2500 & 0.0000 & 0.2896 & 0.2082 & 1.389 \\
\texttt{gmu-35-40} & T & 0.0 & 0.0 & 0.0001 & 0.0001 & 0.0003 & 0.0003 & 0.967 \\
\texttt{gmu-35-50} & H & 0.1 & 0.1 & 0.0002 & 0.0002 & 0.0012 & 0.0016 & 0.849 \\
\texttt{graph20-20-1rand} & T & 0.1 & 0.1 & 0.0000 & 0.0000 & 0.0003 & 0.0727 & 0.018 \\
\texttt{graphdraw-domain} & TC & 0.0 & 0.0 & 0.0000 & 0.0000 & 0.0032 & 0.0195 & 0.204 \\
\texttt{h80x6320d} & TC & 0.1 & 0.2 & 0.0000 & 0.0000 & 0.0067 & 0.0070 & 0.951 \\
\texttt{highschool1-aigio} & T & 8.3 & 8.0 & 1.0000 & 1.0000 & 1.0138 & 1.0133 & 1.000 \\
\texttt{hypothyroid-k1} & HK & 2.3 & 2.3 & 0.0000 & 0.0000 & 0.0109 & 0.0129 & 0.858 \\
\texttt{ic97\_potential} & T & 4.1 & 4.0 & 0.0000 & 0.0000 & 0.0161 & 0.0151 & 1.064 \\
\texttt{icir97\_tension} & HK & 5.6 & 5.3 & 0.0022 & 0.0022 & 0.0209 & 0.0199 & 1.048 \\
\texttt{irish-electricity} & H & 71.3 & 70.1 & 0.0079 & 0.0079 & 0.2446 & 0.2406 & 1.016 \\
\texttt{irp} & H & 4.6 & 4.6 & 0.0000 & 0.0000 & 0.0158 & 0.0180 & 0.885 \\
\texttt{istanbul-no-cutoff} & TC & 0.5 & 1.4 & 0.0000 & 0.0000 & 0.0138 & 0.0193 & 0.729 \\
\texttt{k1mushroom} & HK & 194.0 & 169.3 & 0.0000 & 0.0000 & 0.9144 & 0.7911 & 1.156 \\
\texttt{lectsched-5-obj} & T & 1.6 & 412.6 & 0.4667 & 0.7647 & 0.6259 & 1.6179 & 0.387 \\
\texttt{leo1} & H & 0.2 & 0.3 & 0.0236 & 0.0105 & 0.0392 & 0.0365 & 1.071 \\
\texttt{leo2} & T & 0.3 & 0.5 & 0.0293 & 0.0288 & 0.0481 & 0.0474 & 1.015 \\
\texttt{lotsize} & H & 7.0 & 7.1 & 0.0071 & 0.0044 & 0.0338 & 0.0324 & 1.044 \\
\texttt{mad} & T & 0.0 & 0.0 & 0.0176 & 0.0184 & 0.0457 & 0.0412 & 1.106 \\
\texttt{map10} & T & 1.0 & 0.9 & 0.0000 & 0.0081 & 0.0940 & 0.0920 & 1.021 \\
\texttt{map16715-04} & TC & 1.0 & 0.9 & 0.0000 & 0.0000 & 0.2022 & 0.2542 & 0.796 \\
\texttt{markshare2} & TC & 0.0 & 0.0 & 0.9630 & 0.9375 & 0.9654 & 0.9564 & 1.009 \\
\texttt{markshare\_4\_0} & HK & 0.0 & 0.0 & 0.0000 & 0.0000 & 0.0239 & 0.0235 & 1.014 \\
\texttt{mas74} & H & 0.0 & 0.0 & 0.0078 & 0.0078 & 0.0150 & 0.0135 & 1.102 \\
\texttt{mas76} & TC & 0.0 & 0.0 & 0.0000 & 0.0000 & 0.0000 & 0.0000 & 1.000 \\
\texttt{mc11} & T & 0.1 & 0.0 & 0.0000 & 0.0000 & 0.0012 & 0.0009 & 1.164 \\
\texttt{mcsched} & H & 0.1 & 0.1 & 0.0000 & 0.0000 & 0.0057 & 0.0073 & 0.810 \\
\texttt{mik-250-20-75-4} & TC & 0.0 & 0.0 & 0.0000 & 0.0000 & 0.0005 & 0.0006 & 0.895 \\
\texttt{milo-v12-6-r2-40-1} & H & 4.1 & 3.9 & 0.0010 & 0.0010 & 0.0183 & 0.0177 & 1.035 \\
\texttt{momentum1} & H & n/a & n/a & n/a & n/a & 2.0000 & 2.0000 & 1.000 \\
\texttt{mushroom-best} & TC & 0.7 & 0.5 & 0.0000 & 0.0000 & 0.0143 & 0.0025 & 4.345 \\
\texttt{mzzv11} & HK & 23.9 & 24.4 & 0.0000 & 0.0000 & 0.1009 & 0.1203 & 0.840 \\
\texttt{mzzv42z} & H & 13.1 & 13.6 & 0.0000 & 0.0000 & 0.0456 & 0.0530 & 0.863 \\
\texttt{n2seq36q} & HK & 1.3 & 1.4 & 0.0000 & 0.0000 & 0.0230 & 0.0379 & 0.618 \\
\texttt{n3div36} & HK & 0.6 & 0.6 & 0.0015 & 0.0015 & 0.0243 & 0.0129 & 1.818 \\
\texttt{n5-3} & HK & 0.1 & 0.0 & 0.0000 & 0.0000 & 0.0053 & 0.0048 & 1.098 \\
\texttt{neos-1122047} & T & 6.2 & 4.5 & 0.0000 & 0.0000 & 0.0207 & 0.0150 & 1.354 \\
\texttt{neos-1171448} & TC & 0.2 & 0.2 & 0.0000 & 0.0000 & 0.0017 & 0.0287 & 0.090 \\
\texttt{neos-1171737} & H & 0.1 & 0.1 & 0.0115 & 0.0205 & 0.0193 & 0.0266 & 0.736 \\
\texttt{neos-1354092} & H & n/a & n/a & n/a & n/a & 2.0000 & 2.0000 & 1.000 \\
\texttt{neos-1445765} & T & 0.2 & 0.2 & 0.0000 & 0.0000 & 0.0078 & 0.0193 & 0.431 \\
\texttt{neos-1456979} & TC & 9.8 & 9.2 & 0.0000 & 0.0000 & 0.0593 & 0.0564 & 1.051 \\
\texttt{neos-1582420} & TC & 0.1 & 5.3 & 0.0000 & 0.0000 & 0.0028 & 0.0198 & 0.182 \\
\texttt{neos-2657525-crna} & T & 55.6 & 50.6 & 0.0000 & 0.0000 & 0.2006 & 0.2312 & 0.868 \\
\texttt{neos-2746589-doon} & HK & n/a & n/a & n/a & n/a & 2.0000 & 2.0000 & 1.000 \\
\texttt{neos-2978193-inde} & HK & 0.1 & 0.2 & 0.0000 & 0.0000 & 0.0006 & 0.0012 & 0.702 \\
\texttt{neos-2987310-joes} & T & 3.6 & 3.6 & 0.0000 & 0.0000 & 0.0165 & 0.0147 & 1.116 \\
\texttt{neos-3004026-krka} & TC & 265.1 & 268.2 & 0.0000 & 0.0000 & 0.8837 & 0.8940 & 0.988 \\
\texttt{neos-3024952-loue} & HK & 547.7 & 541.7 & 0.0280 & 0.0280 & 1.8300 & 1.8102 & 1.011 \\
\texttt{neos-3046615-murg} & TC & 0.6 & 0.4 & 0.0025 & 0.0000 & 0.0060 & 0.0043 & 1.305 \\
\texttt{neos-3083819-nubu} & TC & 0.1 & 0.9 & 0.0000 & 0.0000 & 0.0005 & 0.0030 & 0.364 \\
\texttt{neos-3216931-puriri} & TC & 58.8 & 55.5 & 0.0000 & 0.0000 & 0.2712 & 0.2577 & 1.052 \\
\texttt{neos-3381206-awhea} & H & 0.4 & 0.4 & 0.0000 & 0.0000 & 0.0014 & 0.0014 & 1.000 \\
\texttt{neos-3402294-bobin} & H & 2.2 & 3.4 & 0.0000 & 0.0000 & 0.0075 & 0.0372 & 0.223 \\
\texttt{neos-3555904-turama} & H & 5.8 & 6.0 & 0.0000 & 0.0432 & 0.0193 & 0.0628 & 0.319 \\
\texttt{neos-3627168-kasai} & H & 1.8 & 1.7 & 0.0000 & 0.0000 & 0.0060 & 0.0057 & 1.049 \\
\texttt{neos-3656078-kumeu} & H & n/a & n/a & n/a & n/a & 2.0000 & 2.0000 & 1.000 \\
\texttt{neos-3754480-nidda} & T & 0.0 & 0.0 & 0.0742 & 0.0585 & 0.0770 & 0.0616 & 1.246 \\
\texttt{neos-4300652-rahue} & T & 0.7 & 0.9 & 0.2700 & 0.4569 & 0.3201 & 0.5782 & 0.554 \\
\texttt{neos-4338804-snowy} & H & 0.0 & 0.0 & 0.0041 & 0.0054 & 0.0084 & 0.0092 & 0.925 \\
\texttt{neos-4387871-tavua} & HK & 8.8 & 8.4 & 0.0453 & 0.0406 & 0.0795 & 0.0742 & 1.070 \\
\texttt{neos-4413714-turia} & HK & 63.3 & 8.6 & 0.0000 & 0.0000 & 0.2113 & 0.0354 & 5.828 \\
\texttt{neos-4532248-waihi} & T & 12.6 & n/a & 0.1324 & n/a & 0.1717 & 2.0000 & 0.086 \\
\texttt{neos-4647030-tutaki} & H & 31.3 & 48.4 & 0.0002 & 0.0002 & 0.9346 & 0.1615 & 5.757 \\
\texttt{neos-4722843-widden} & HK & 22.3 & 18.5 & 0.0000 & 0.0000 & 0.0834 & 0.0706 & 1.179 \\
\texttt{neos-4738912-atrato} & T & 0.1 & 0.1 & 0.0000 & 0.0000 & 0.0058 & 0.0095 & 0.650 \\
\texttt{neos-4763324-toguru} & H & 4.1 & 9.4 & 0.1878 & 0.1861 & 0.4102 & 0.4116 & 0.997 \\
\texttt{neos-4954672-berkel} & H & 0.0 & 0.0 & 0.0060 & 0.0125 & 0.0158 & 0.0158 & 1.002 \\
\texttt{neos-5049753-cuanza} & T & 3.5 & 7.1 & 0.0226 & 0.4732 & 0.0920 & 0.6595 & 0.141 \\
\texttt{neos-5052403-cygnet} & T & 5.3 & 6.3 & 0.0808 & 0.4033 & 0.0989 & 0.4200 & 0.237 \\
\texttt{neos-5093327-huahum} & H & 17.0 & 14.4 & 0.1161 & 0.1161 & 0.1694 & 0.1613 & 1.050 \\
\texttt{neos-5104907-jarama} & TC & n/a & n/a & n/a & n/a & 2.0000 & 2.0000 & 1.000 \\
\texttt{neos-5107597-kakapo} & H & 403.0 & 378.1 & 0.8966 & 0.8979 & 1.6382 & 1.5924 & 1.029 \\
\texttt{neos-5114902-kasavu} & T & 11.3 & 23.9 & 0.1301 & 0.6984 & 0.1658 & 0.7503 & 0.222 \\
\texttt{neos-5188808-nattai} & TC & 13.3 & 12.3 & 0.0000 & 0.0000 & 0.0567 & 0.0539 & 1.051 \\
\texttt{neos-5195221-niemur} & T & 1.3 & 227.3 & 0.0005 & 0.0077 & 0.0048 & 0.7634 & 0.008 \\
\texttt{neos-631710} & HK & 35.2 & 35.8 & 0.0000 & 0.6336 & 0.1177 & 0.7151 & 0.166 \\
\texttt{neos-662469} & HK & 1.2 & 1.1 & 0.0001 & 0.0001 & 0.1667 & 0.1781 & 0.936 \\
\texttt{neos-787933} & TC & 4.3 & 3.9 & 0.0000 & 0.0323 & 0.0706 & 0.0829 & 0.853 \\
\texttt{neos-827175} & HK & 0.5 & 0.5 & 0.0000 & 0.0000 & 0.0017 & 0.0019 & 0.929 \\
\texttt{neos-848589} & H & 20.3 & 23.0 & 0.0256 & 0.0089 & 0.3981 & 0.4214 & 0.945 \\
\texttt{neos-860300} & H & 3.4 & 2.8 & 0.0000 & 0.0000 & 0.0159 & 0.0136 & 1.154 \\
\texttt{neos-873061} & TC & 9.4 & 4.8 & 0.0001 & 0.0001 & 0.0349 & 0.0193 & 1.766 \\
\texttt{neos-911970} & H & 0.0 & 0.0 & 0.0000 & 0.0000 & 0.0015 & 0.0013 & 1.068 \\
\texttt{neos-933966} & TC & 0.6 & 0.7 & 0.0000 & 0.0000 & 0.0544 & 0.0656 & 0.832 \\
\texttt{neos-950242} & HK & 0.3 & 0.3 & 0.0000 & 0.0000 & 0.0010 & 0.0073 & 0.242 \\
\texttt{neos-957323} & H & 5.7 & 5.6 & 0.0000 & 0.0000 & 0.0200 & 0.0466 & 0.442 \\
\texttt{neos-960392} & TC & 1.7 & 1.7 & 0.0000 & 0.0000 & 0.0058 & 0.0225 & 0.291 \\
\texttt{neos17} & HK & 0.0 & 0.0 & 0.0000 & 0.0000 & 0.0001 & 0.0001 & 1.010 \\
\texttt{neos5} & HK & 0.0 & 0.0 & 0.0000 & 0.0000 & 0.0000 & 0.0001 & 0.975 \\
\texttt{neos8} & T & 3.4 & 3.0 & 0.0000 & 0.0000 & 0.0115 & 0.0102 & 1.124 \\
\texttt{net12} & H & 0.4 & 0.3 & 0.0000 & 0.0000 & 0.0361 & 0.0109 & 3.122 \\
\texttt{netdiversion} & TC & 3.3 & 3.8 & 0.0000 & 0.0000 & 0.1277 & 0.3611 & 0.355 \\
\texttt{nexp-150-20-8-5} & HK & 0.1 & 0.2 & 0.0086 & 0.0043 & 0.0606 & 0.1105 & 0.552 \\
\texttt{ns1116954} & TC & 1.4 & 1.6 & 0.0000 & 0.0000 & 0.0047 & 0.0053 & 0.895 \\
\texttt{ns1208400} & TC & 86.3 & 81.7 & 0.0000 & 0.0000 & 0.2877 & 0.2723 & 1.056 \\
\texttt{ns1644855} & TC & n/a & 600.0 & n/a & 0.0776 & 2.0000 & 2.0000 & 1.000 \\
\texttt{ns1760995} & HK & n/a & n/a & n/a & n/a & 2.0000 & 2.0000 & 1.000 \\
\texttt{ns1830653} & H & 0.2 & 0.2 & 0.0000 & 0.0000 & 0.0683 & 0.0729 & 0.938 \\
\texttt{ns1952667} & H & n/a & n/a & n/a & n/a & 2.0000 & 2.0000 & 1.000 \\
\texttt{nu25-pr12} & H & 0.2 & 0.1 & 0.0000 & 0.0000 & 0.0007 & 0.0004 & 1.241 \\
\texttt{nursesched-medium-hint03} & TC & 28.9 & 31.4 & 0.2123 & 0.2123 & 0.7452 & 0.7356 & 1.013 \\
\texttt{nursesched-sprint02} & TC & 12.5 & 25.0 & 0.0000 & 0.0000 & 0.0608 & 0.0833 & 0.733 \\
\texttt{nw04} & H & 30.6 & 31.2 & 0.0000 & 0.0000 & 0.1299 & 0.1115 & 1.163 \\
\texttt{opm2-z10-s4} & HK & 3.1 & 2.8 & 0.0506 & 0.0508 & 0.1894 & 0.2262 & 0.838 \\
\texttt{p200x1188c} & H & 0.0 & 0.0 & 0.0000 & 0.0000 & 0.0000 & 0.0000 & 1.000 \\
\texttt{peg-solitaire-a3} & HK & n/a & n/a & n/a & n/a & 2.0000 & 2.0000 & 1.000 \\
\texttt{pg} & H & 0.0 & 0.0 & 0.0000 & 0.0000 & 0.0020 & 0.0020 & 0.989 \\
\texttt{pg5\_34} & H & 0.0 & 0.0 & 0.0000 & 0.0000 & 0.0004 & 0.0003 & 1.050 \\
\texttt{physiciansched3-3} & H & n/a & n/a & n/a & n/a & 2.0000 & 2.0000 & 1.000 \\
\texttt{physiciansched6-2} & HK & 1.3 & 9.1 & 0.0000 & 0.0000 & 0.0046 & 0.0303 & 0.179 \\
\texttt{piperout-08} & H & 22.2 & 22.8 & 0.0000 & 0.0000 & 0.0740 & 0.0760 & 0.974 \\
\texttt{piperout-27} & H & 37.1 & 36.8 & 0.0000 & 0.0000 & 0.1237 & 0.1227 & 1.008 \\
\texttt{pk1} & T & 0.0 & 0.0 & 0.0000 & 0.0000 & 0.0768 & 0.0314 & 2.401 \\
\texttt{proteindesign121hz512p9} & TC & 6.4 & 118.2 & 0.0121 & 0.0061 & 0.0385 & 0.4343 & 0.091 \\
\texttt{proteindesign122trx11p8} & TC & 6.2 & 77.6 & 0.0068 & 0.0068 & 0.0288 & 0.2688 & 0.111 \\
\texttt{qap10} & T & 0.5 & 0.5 & 0.0000 & 0.0000 & 0.0128 & 0.0238 & 0.555 \\
\texttt{radiationm18-12-05} & HK & 0.3 & 29.8 & 0.0002 & 0.1000 & 0.1950 & 0.3274 & 0.597 \\
\texttt{radiationm40-10-02} & H & 1.0 & n/a & 0.9172 & n/a & 0.9306 & 2.0000 & 0.466 \\
\texttt{rail01} & HK & n/a & n/a & n/a & n/a & 2.0000 & 2.0000 & 1.000 \\
\texttt{rail02} & HK & n/a & n/a & n/a & n/a & 2.0000 & 2.0000 & 1.000 \\
\texttt{rail507} & TC & 1.2 & 1.4 & 0.0057 & 0.0057 & 0.0145 & 0.0229 & 0.645 \\
\texttt{ran14x18-disj-8} & H & 0.0 & 0.0 & 0.0005 & 0.0000 & 0.0083 & 0.0054 & 1.440 \\
\texttt{rd-rplusc-21} & H & 34.2 & 30.9 & 0.0020 & 0.0026 & 0.1170 & 0.1060 & 1.102 \\
\texttt{reblock115} & H & 0.3 & 0.2 & 0.0001 & 0.0000 & 0.0148 & 0.0265 & 0.574 \\
\texttt{rmatr100-p10} & H & 0.0 & 0.3 & 0.0000 & 0.0000 & 0.0011 & 0.0100 & 0.191 \\
\texttt{rmatr200-p5} & TC & 0.2 & 7.0 & 0.0311 & 0.0856 & 0.0927 & 0.1468 & 0.634 \\
\texttt{rocI-4-11} & H & 12.8 & 12.0 & 0.0000 & 0.1611 & 0.1959 & 0.2038 & 0.962 \\
\texttt{rocII-5-11} & T & 16.2 & 15.4 & 0.1500 & 0.1500 & 0.2189 & 0.2166 & 1.011 \\
\texttt{rococoB10-011000} & HK & 0.1 & 0.1 & 0.0302 & 0.0308 & 0.0816 & 0.0912 & 0.895 \\
\texttt{rococoC10-001000} & H & 0.0 & 0.0 & 0.0000 & 0.0000 & 0.0077 & 0.0070 & 1.082 \\
\texttt{roi2alpha3n4} & H & 6.9 & 6.5 & 0.0000 & 0.0000 & 0.0647 & 0.0507 & 1.271 \\
\texttt{roi5alpha10n8} & T & 23.7 & 21.9 & 0.1026 & 0.0902 & 0.2860 & 0.3462 & 0.826 \\
\texttt{roll3000} & HK & 2.9 & 2.7 & 0.0000 & 0.0000 & 0.0104 & 0.0095 & 1.088 \\
\texttt{s100} & HK & 75.9 & 84.7 & 0.0192 & 0.0307 & 0.2698 & 0.3087 & 0.874 \\
\texttt{s250r10} & T & 29.6 & 106.3 & 0.0000 & 0.0000 & 0.1231 & 0.3674 & 0.337 \\
\texttt{satellites2-40} & H & 4.7 & 4.9 & 0.0000 & 0.0000 & 0.3643 & 0.5995 & 0.608 \\
\texttt{satellites2-60-fs} & HK & 3.0 & 116.6 & 0.0000 & 0.0000 & 0.2291 & 0.5657 & 0.406 \\
\texttt{savsched1} & TC & 11.7 & 336.1 & 0.7596 & 0.7596 & 0.8721 & 1.4544 & 0.600 \\
\texttt{sct2} & TC & 0.1 & 0.2 & 0.0000 & 0.0000 & 0.0011 & 0.0022 & 0.640 \\
\texttt{seymour} & H & 0.1 & 0.1 & 0.0000 & 0.0047 & 0.0005 & 0.0097 & 0.141 \\
\texttt{seymour1} & HK & 0.0 & 0.1 & 0.0000 & 0.0000 & 0.0007 & 0.0015 & 0.680 \\
\texttt{sing326} & HK & 8.6 & 7.0 & 0.1578 & 0.1578 & 0.3785 & 0.3788 & 0.999 \\
\texttt{sing44} & H & 3.2 & 9.0 & 0.0464 & 0.0464 & 0.1194 & 0.1258 & 0.950 \\
\texttt{snp-02-004-104} & HK & 20.1 & 13.3 & 0.0000 & 0.0000 & 0.1263 & 0.1013 & 1.244 \\
\texttt{sorrell3} & TC & 2.0 & 2.1 & 0.0000 & 0.1250 & 0.0070 & 0.2736 & 0.029 \\
\texttt{sp150x300d} & HK & 0.0 & 0.0 & 0.0000 & 0.0000 & 0.0000 & 0.0000 & 1.000 \\
\texttt{sp97ar} & H & 0.5 & 0.7 & 0.0234 & 0.0237 & 0.0413 & 0.0432 & 0.959 \\
\texttt{sp98ar} & H & 0.6 & 0.9 & 0.0097 & 0.0104 & 0.0263 & 0.0338 & 0.787 \\
\texttt{splice1k1} & H & 22.6 & 20.4 & 0.1447 & 0.1447 & 0.2801 & 0.3329 & 0.842 \\
\texttt{square41} & H & 155.4 & n/a & 0.8905 & n/a & 1.1779 & 2.0000 & 0.589 \\
\texttt{square47} & H & 331.4 & n/a & 0.9030 & n/a & 1.5090 & 2.0000 & 0.755 \\
\texttt{supportcase10} & H & 3.9 & 3.6 & 0.1250 & 0.6667 & 0.1375 & 0.7255 & 0.191 \\
\texttt{supportcase12} & H & 2.9 & 2.7 & 0.0018 & 0.0018 & 0.0364 & 0.0691 & 0.533 \\
\texttt{supportcase18} & TC & 0.2 & 0.2 & 0.0204 & 0.0204 & 0.0218 & 0.0230 & 0.952 \\
\texttt{supportcase19} & H & n/a & n/a & n/a & n/a & 2.0000 & 2.0000 & 1.000 \\
\texttt{supportcase22} & H & n/a & n/a & n/a & n/a & 2.0000 & 2.0000 & 1.000 \\
\texttt{supportcase26} & T & 0.0 & 0.0 & 0.0011 & 0.0000 & 0.0088 & 0.0060 & 1.402 \\
\texttt{supportcase33} & H & 10.5 & 10.8 & 0.0000 & 0.0000 & 0.0645 & 0.0750 & 0.862 \\
\texttt{supportcase40} & H & 0.1 & 0.5 & 0.0067 & 0.0084 & 0.0244 & 0.0474 & 0.526 \\
\texttt{supportcase42} & H & 0.3 & 0.4 & 0.0302 & 0.0234 & 0.2249 & 0.1446 & 1.552 \\
\texttt{supportcase6} & HK & 6.6 & 7.4 & 0.0000 & 0.0005 & 0.0828 & 0.0956 & 0.867 \\
\texttt{supportcase7} & T & 5.3 & 4.7 & 0.0000 & 0.0000 & 0.0704 & 0.0710 & 0.991 \\
\texttt{swath1} & T & 0.5 & 0.9 & 0.0000 & 0.0000 & 0.0058 & 0.0033 & 1.563 \\
\texttt{swath3} & H & 0.6 & 0.7 & 0.0000 & 0.0000 & 0.0074 & 0.0047 & 1.484 \\
\texttt{tbfp-network} & T & 2.3 & 2.6 & 0.0327 & 0.0000 & 0.2003 & 0.1906 & 1.050 \\
\texttt{thor50dday} & H & 0.5 & 0.6 & 0.0214 & 0.0083 & 0.0741 & 0.0447 & 1.641 \\
\texttt{timtab1} & TC & 1.8 & 1.6 & 0.0000 & 0.0000 & 0.0082 & 0.0075 & 1.078 \\
\texttt{tr12-30} & H & 0.1 & 0.1 & 0.0000 & 0.0000 & 0.0009 & 0.0009 & 1.051 \\
\texttt{traininstance2} & H & 11.5 & 11.4 & 0.0086 & 0.0077 & 0.0799 & 0.0791 & 1.010 \\
\texttt{traininstance6} & H & 2.8 & 2.4 & 0.0000 & 0.0000 & 0.0168 & 0.0154 & 1.083 \\
\texttt{trento1} & T & 0.4 & 2.9 & 0.0001 & 0.0002 & 0.0311 & 0.0308 & 1.009 \\
\texttt{triptim1} & H & 103.2 & 88.3 & 0.0000 & 0.0000 & 0.3440 & 0.2943 & 1.168 \\
\texttt{uccase12} & TC & 5.2 & 3.6 & 0.0000 & 0.0000 & 0.0185 & 0.0129 & 1.404 \\
\texttt{uccase9} & T & 9.9 & 9.3 & 0.0391 & 0.0391 & 0.3556 & 0.3542 & 1.004 \\
\texttt{uct-subprob} & TC & 0.1 & 0.1 & 0.0032 & 0.0000 & 0.0204 & 0.0282 & 0.733 \\
\texttt{unitcal\_7} & HK & 4.0 & 3.9 & 0.0000 & 0.0000 & 0.0146 & 0.0144 & 1.017 \\
\texttt{var-smallemery-m6j6} & TC & 8.2 & 5.4 & 0.0165 & 0.0096 & 0.0638 & 0.0447 & 1.418 \\
\texttt{wachplan} & H & 8.6 & 7.8 & 0.0000 & 0.0000 & 0.0287 & 0.0260 & 1.099 \\
\end{longtable}
}